\documentclass[12pt, a4article]{amsart}   	
\usepackage[top=1.15in, bottom=1.15in, left=1.17in, right=1.17in]{geometry}
\usepackage[dvipsnames,svgnames,x11names]{xcolor}
\usepackage[colorlinks=true, pdfstartview=FitV, linkcolor=blue, citecolor=blue, urlcolor=blue, breaklinks=true]{hyperref}
\usepackage{stmaryrd}
\usepackage{DotArrow}
\usepackage{amsmath, amsthm}
\usepackage{xypic}
\usepackage[markup=underlined]{changes}
\usepackage{epsf}
\usepackage[symbols,nogroupskip]{glossaries-extra}
\usepackage{glossary-mcols}
\usepackage{tikz-cd}
\usepackage[T1]{fontenc}

\renewcommand*{\glossarymark}[1]{}
\makenoidxglossaries

\glsxtrnewsymbol[description={extremal functions;}]{fp}{\ensuremath{f_p}}
\glsxtrnewsymbol[description={strata in Seshadri stratification;}]{Xp}{\ensuremath{X_p}}
\glsxtrnewsymbol[description={poset indexing strata in Seshadri stratification;}]{A}{\ensuremath{A}}
\glsxtrnewsymbol[description={fixed embedded projective variety;}]{X}{\ensuremath{X}}
\glsxtrnewsymbol[description={fixed vector space;}]{V}{\ensuremath{V}}
\glsxtrnewsymbol[description={subposet of $A$;}]{Ap}{\ensuremath{A_p}}
\glsxtrnewsymbol[description={length function on $A$;}]{l}{\ensuremath{\ell}}
\glsxtrnewsymbol[description={Grassmann variety}]{Grdn}{\ensuremath{\mathrm{Gr}_d\mathbb{K}^n}}
\glsxtrnewsymbol[description={affine cone over $X_p$;}]{Xphat}{\ensuremath{\hat{X}_p}}
\glsxtrnewsymbol[description={bonds between $p$ and $q$;}]{bpq}{\ensuremath{b_{p,q}}}
\glsxtrnewsymbol[description={normalization of $X_p$;}]{Xptilde}{\ensuremath{\tilde{X}_p}}
\glsxtrnewsymbol[description={homogeneous coordinate ring of $X_p$;}]{Rp}{\ensuremath{R_p}}
\glsxtrnewsymbol[description={valuation from vanishing multiplicity;}]{vpq}{\ensuremath{\nu_{p,q}}}
\glsxtrnewsymbol[description={principal ideal in $R_p$ generated by $f_p$;}]{Ip}{\ensuremath{I_p}}
\glsxtrnewsymbol[description={quasi-valuation associated to $I_p$;}]{nuIp}{\ensuremath{\nu_{I_p}}}
\glsxtrnewsymbol[description={homogenized quasi-valuation;}]{nuIpbar}{\ensuremath{\overline{\nu}_{I_p}}}
\glsxtrnewsymbol[description={maximal chain in $A$;}]{C}{\ensuremath{\mathfrak{C}}}
\glsxtrnewsymbol[description={sequence of rational function associated to $g$ along $\mathfrak{C}$;}]{gc}{\ensuremath{g_{\mathfrak{C}}}}
\glsxtrnewsymbol[description={normalized valuation associated to a maximal chain;}]{Vc}{\ensuremath{\mathcal{V}_{\mathfrak{C}}}}
\glsxtrnewsymbol[description={a variation of the valuation $\mathcal{V}_{\mathfrak{C}}$;}]{Vctilde}{\ensuremath{\widetilde{\mathcal{V}}_{\mathfrak{C}}}}
\glsxtrnewsymbol[description={a lattice;}]{Lc}{\ensuremath{L^{\mathfrak{C}}}}
\glsxtrnewsymbol[description={valuation monoid of $\mathcal{V}_{\mathfrak{C}}$;}]{VCX}{\ensuremath{\mathbb{V}_{\mathfrak{C}}}}
\glsxtrnewsymbol[description={submonoid of lattice $L^{\mathfrak{C}}$;}]{Lcdag}{\ensuremath{L^{\mathfrak{C},\dag}}}
\glsxtrnewsymbol[description={associated graded algebra of $\mathcal{V}_{\mathfrak{C}}$;}]{grcR}{\ensuremath{\mathrm{gr}_{\mathfrak{C}}R}}
\glsxtrnewsymbol[description={core of valuation monoid $\mathcal{V}_{\mathfrak{C}}(X)$;}]{PCX}{\ensuremath{P_{\mathfrak{C}}(X)}}
\glsxtrnewsymbol[description={an open affine subset of $\hat{X}$;}]{UC}{\ensuremath{U_{\mathfrak{C}}}}
\glsxtrnewsymbol[description={valuation monoid on $U_{\mathfrak{C}}$;}]{VCUC}{\ensuremath{\mathbb{V}_{\mathfrak C}(U_{\mathfrak{C}})}}
\glsxtrnewsymbol[description={fixed total order on $A$;}]{geqt}{\ensuremath{>^t}}
\glsxtrnewsymbol[description={set of all maximal chains in $A$;}]{CC}{\ensuremath{\mathcal{C}}}
\glsxtrnewsymbol[description={minimized quasi-valuation of function $g$;}]{Vg}{\ensuremath{\mathcal{V}(g)}}
\glsxtrnewsymbol[description={support of function g;}]{suppg}{\ensuremath{\mathrm{supp}\mathcal{V}(g)}}
\glsxtrnewsymbol[description={the fan of monoid associated to $\mathcal{V}$;}]{Gamma}{\ensuremath{\Gamma}}
\glsxtrnewsymbol[description={maximal chains hitting minimum;}]{Cg}{\ensuremath{\mathcal{C}(g)}}
\glsxtrnewsymbol[description={cone associated to chain $C$ in $A$;}]{Kc}{\ensuremath{K_C}}
\glsxtrnewsymbol[description={monoid in $\Gamma$;}]{GammaC}{\ensuremath{\Gamma_{\mathfrak{C}}}}
\glsxtrnewsymbol[description={the fan algebra associated to $\Gamma$;}]{KGamma}{\ensuremath{\mathbb{K}[\Gamma]}}
\glsxtrnewsymbol[description={associated graded algebra of quasi-valuation $\mathcal{V}$;}]{grvR}{\ensuremath{\mathrm{gr}_{\mathcal{V}}R}}
\glsxtrnewsymbol[description={subalgebra of $\mathrm{gr}_{\mathcal{V}}R$;}]{grvcR}{\ensuremath{\mathrm{gr}_{\mathcal{V},\mathfrak{C}}R}}
\glsxtrnewsymbol[description={annihilating ideal associated to $\mathfrak{C}$;}]{Ic}{\ensuremath{I_{\mathfrak{C}}}}
\glsxtrnewsymbol[description={Hasse graph with bonds;}]{GA}{\ensuremath{\mathcal{G}_A}}
\glsxtrnewsymbol[description={l.c.m of all bonds in $\mathcal{G}_A$;}]{N}{\ensuremath{N}}
\glsxtrnewsymbol[description={set of all chains in $A$;}]{K}{\ensuremath{\mathcal{K}}}
\glsxtrnewsymbol[description={set of all saturated chains in $A$;}]{Ks}{\ensuremath{\mathcal{K}^s}}
\glsxtrnewsymbol[description={affine variety associated to saturated set $\mathcal{M}$;}]{ZM}{\ensuremath{Z_\mathcal{M}}}
\glsxtrnewsymbol[description={Newton-Okounkov complex associated tp Seshadri stratification;}]{DV}{\ensuremath{\Delta_{\mathcal{V}}}}
\glsxtrnewsymbol[description={lattice generated by $\Gamma_{\mathfrak{C}}$;}]{LUC}{\ensuremath{\mathcal{L}^{\mathfrak{C}}}}
\glsxtrnewsymbol[description={extended poset;}]{Ahat}{\ensuremath{\hat{A}}}

\def\mfC{\mathfrak C}
\DeclareMathOperator{\supp}{supp}

\theoremstyle{plain}
\newtheorem{IntroThm}{Theorem}[]
\newtheorem{IntroProp}{Proposition}[]
\newtheorem{IntroCor}{Corollary}[]
\newtheorem{theorem}{Theorem}[section]
\newtheorem{coro}[theorem]{Corollary}
\newtheorem{proposition}[theorem]{Proposition}
\newtheorem{lemma}[theorem]{Lemma}

\theoremstyle{definition}
\newtheorem{algo}[theorem]{Algorithm}
\newtheorem{rem}[theorem]{Remark}
\newtheorem{example}[theorem]{Example}
\newtheorem{definition}[theorem]{Definition}

\title{Seshadri stratifications and standard monomial theory}
\author{Rocco Chiriv\`i}
\address{Dipartimento di Matematica e Fisica ``Ennio De Giorgi'', Universit\`a del Salento, Lecce, Italy}
\email{rocco.chirivi@unisalento.it}

\author{Xin Fang} 
\address{Lehrstuhl f\"ur Algebra und Darstellungstheorie, RWTH Aachen University, Pontdriesch 10-16, 52062 Aachen, Germany}
\email{xinfang.math@gmail.com}

\author{Peter Littelmann}
\address{Department Mathematik/Informatik, Universit\"at zu K\"oln, 50931, Cologne, Germany}
\email{peter.littelmann@math.uni-koeln.de}

\begin{document}
\maketitle

\begin{abstract}
We introduce the notion of a Seshadri stratification on an embedded projective variety. Such a structure enables us to construct a Newton-Okounkov simplicial complex and a flat degeneration of the projective variety into a union of toric varieties. We show that the Seshadri stratification provides a geometric setup for a standard monomial theory. In this framework, Lakshmibai-Seshadri paths for Schubert varieties  get a geometric interpretation as successive vanishing orders of regular functions.
\end{abstract}

\section{Introduction}

We fix throughout the article an algebraically closed field $\mathbb{K}$.

The aim of the article is to develop a theory parallel to that of Newton-Okounkov bodies, built on a web rather than a flag of subvarieties. The other ingredient making our approach different from that in Newton-Okounkov theory is a finite collection of functions with a prescribed set-theoretical vanishing behavior, leading to the notion of a Seshadri stratification. 
Compared to the Newton-Okounkov theory: instead of a valuation we have a quasi-valuation, but with values in the positive orthant; instead of a single monoid we obtain a fan of monoids, but the monoids in this fan are always finitely generated; instead of a body in an Euclidean space we get a simplicial complex with a rational structure.

As an example of this comparison, in the case of flag varieties, in the same way in which the string cones and their associated monoids \cite{BZ01, Lit98} show up in the application of Newton-Okounkov theory \cite{K1}, our setup leads to a polyhedral geometric version of the Lakshmibai-Seshadri path model \cite{L2}.
 
Before going into the details, we review a few aspects of the various versions of standard monomial theory, which have been part of the motivating background for this article.

\subsection{Theories of standard monomials}

One of the motivations for theories of standard monomials is the computation of the Hilbert function of a graded finitely generated algebra $R$ over a field $\mathbb{K}$. As their name suggests, there are two choices to make: what are the generators of the algebra $R$, and which monomials in the generators are chosen to be ``standard''.

\subsubsection{Algebraic setting}
To the best of our knowledge, the first work in this direction is by Macaulay \cite{Mac}. He considered the case $R=\mathbb{K}[x_1,\ldots,x_n]/I$ where $I$ is a homogeneous ideal. The important idea of Macaulay is to mix the structure of an order into the algebraic structure, transforming $R$ to a ``simpler'' algebra sharing the same Hilbert function as $R$.

The theory of Gr\"obner basis, introduced by Buchberger \cite{Buch}, associates a unique reduced Gr\"obner basis $\mathrm{GB}(I,>)$ to the ideal $I$ and a fixed monomial order $>$. Monomials in $x_1,\ldots,x_n$, which are not contained in the initial ideal $\mathrm{in}_>(I)$, are chosen to be \emph{standard}. The standard monomials form a basis of $\mathbb{K}[x_1,\ldots,x_n]/\mathrm{in}_>(I)$, which share the same Hilbert polynomial as $R$. Determining the Hilbert function of $R$ is thus reduced to a purely combinatorial problem of counting standard monomials. 

The reduced Gr\"obner basis contains further information: each element in $\mathrm{GB}(I,>)$ has the form 
$$\text{a non-standard monomial }+\text{ a linear combination of standard monomials}.$$
Not only does it tell which monomials are standard, but also how to rewrite a non-standard one as a linear combination of standard monomials. Elements in a reduced Gr\"obner basis are called straightening laws.

\subsubsection{Algebro-geometric setting}

Hodge \cite{Hodge} studied this problem when $R$ is the homogeneous coordinate ring of a Grassmann variety or a Schubert subvariety in the Pl\"ucker embedding. The monomials are those in Pl\"ucker coordinates; a monomial is standard if the associated Young tableau is standard. He proved that the standard monomials form a basis of $R$, which allows him to deduce the postulation formula describing the Hilbert function. The Pl\"ucker relations have been used to write the non-standard monomials as linear combinations of standard ones. 

The idea of Hodge is extracted in the work of De Concini, Eisenbud and Procesi \cite{DEP} (see also \cite{Eis}), where they coined the name ``Hodge algebra'' (a.k.a. algebra with straightening laws). Such an algebra $R$ is defined with the following choices: a generating set of $R$ indexed by a partially ordered set (poset), a linear basis consisting of \emph{standard} monomials, i.e. those supported on a maximal chain in the poset, and a rule how to write non-standard monomials into a linear combination of standard ones. Verifying several geometric properties (Gorenstein property, Cohen-Macaulayness, etc) of the algebra $R$ can be reduced to combinatorial problems. Later De Concini and Lakshmibai \cite{DL} generalized Hodge algebras to doset algebras.

\subsubsection{Geometric setting}

Seshadri \cite{Ses} generalizes the work of Hodge from Grassmann varieties to a partial flag variety $G/P$ with $G$ a reductive algebraic group and $P$ a minuscule maximal parabolic subgroup in $G$. In this work he took a geometric approach to the standard monomials in order to avoid applying the explicit straightening relations, and set up the paradigm of deducing geometric properties, such as vanishing of higher cohomology, normality and singular locus of Schubert varieties in $G/P$ from the existence of a \emph{standard monomial theory}. Motivated by the work of De Concini and Procesi \cite{DCP}, in collaboration with Lakshmibai and Musili \cite{LSII,LSIII,LSIV,LS}, Seshadri succeeded in generalizing the results in \cite{Ses} to Schubert varieties in $G/Q$, where $Q$ is a parabolic subgroup of classical type in $G$, by introducing the notion of admissible pairs. This case corresponds to the doset algebras above.

Going beyond classical type, the definition of admissible pairs becomes involved. Lakshmibai (see \cite{LS2}, or Appendix C of \cite{SMT} for a reprint) made a conjecture on a possible index system of a basis of $\mathrm{H}^0(X(\tau),\mathcal{L}_\lambda)$ where $\tau\in W$ is an element in the Weyl group $W$ of $G$, $X(\tau)\subseteq G/Q$ is the corresponding Schubert variety and $\mathcal{L}_\lambda$ is an ample line bundle on $G/Q$ associated to a dominant weight $\lambda$. Such an index system consists of a chain of elements in the Bruhat graph of $W/W_Q$ below $\tau$, with $W_Q$ the Weyl group of $Q$, together with a sequence of rational numbers. It was meanwhile asked to associate an explicit global section to each element in the index system.

\subsection{Path models and standard monomial bases}\label{Sec:IntroPathSMT}

The conjecture of Lakshmibai on the indexing systems is established by the third author in \cite{L2,L3} as a special path model consisting of Lakshmibai-Seshadri (LS)-paths of shape $\lambda$. The LS-paths are piece-wise linear paths starting from the origin in the dual space of a fixed Cartan subalgebra in the Lie algebra of $G$, and their endpoints coincide with the weights appearing in $V(\lambda)$, the Weyl module of $G$ of highest weight $\lambda$, counted with multiplicities. This gives a type-free positive character formula of $V(\lambda)$, i.e. without cancellations like in the Weyl character formula.

Later in \cite{L1}, the third author solved the question about the construction of global sections. For each LS-paths $\pi$ of shape $\lambda$, using the Frobenius map of Lusztig in quantum groups at roots of unity, he constructed a path vector $p_\pi\in\mathrm{H}^0(G/Q,\mathcal{L}_\lambda)$ such that when $\pi$ runs over all LS-paths of shape $\lambda$, the path vectors $p_\pi$ form a basis of the space of global sections. Such constructions are compatible with Schubert varieties (in the sense of \emph{loc.cit}).

The first author \cite{Chi} introduced LS-algebras, which further generalized the above-mentioned work on Hodge algebras and doset algebras, to establish an algebro-geometric setting of the constructions in \cite{L2,L3,L1}. An LS-algebra can be degenerated to a much simpler LS-algebra (called discrete LS-algebra). Together with results in \cite{L1}, this gives a proof of normality and Koszul property of the Schubert varieties by transforming these properties to combinatorics of the discrete LS-algebra (see also \cite{Chi2}).

The LS-paths were defined in a combinatorial way, and their geometric interpretation was missing. We quote the following observation/question by Seshadri from \cite{Ses2}:
\vskip 0.35cm
$\quad$
\begin{minipage}[c]{0,85\textwidth}
{\small\emph{``The character formula via paths or standard diagrams is a formula which involves only the cellular decomposition and its topological properties. It leads one to suspect that there could be a \guillemotleft cellular Riemann-Roch\guillemotright\,\,which could also explain the character formula.''}}
\end{minipage}
\vskip 0.35cm
One of the goals of this article is to build up a framework to provide a geometric interpretation of LS-paths, and at the same time, generalizing them from Schubert varieties to projective varieties with a Seshadri stratification (see below). An algebraic approach has been already taken in \cite{CFL2} by establishing a connection between LS-algebras and valuation theory, generalizing results in \cite{Chi}.

\subsection{Newton-Okounkov theory}\label{Sec:IntroNOB}

Newton-Okounkov bodies first appeared in the work of Okounkov \cite{O1}. His construction has been systemized by Kaveh-Khovanskii \cite{KK} and Lazarsfeld-Musta\c{t}$\breve{\rm a}$ \cite{LM} into the theory of Newton-Okounkov bodies. These discrete geometric objects received great attention in the past ten years. 

To be more precise, the inputs of this machinery are an embedded projective variety $Y$, a flag $Y_\bullet:=(Y=Y_r\supseteq Y_{r-1}\supseteq\ldots\supseteq Y_0=\{\mathrm{pt}\})$ of normal subvarieties (we assume the normality only for simplicity), a collection of rational functions $u_r,\ldots,u_1\in\mathbb{K}(Y)$ such that the restriction of $u_k$ to $Y_k$ is a uniformizer in $\mathcal{O}_{Y_{k},Y_{k-1}}$, and a total order on $\mathbb{Z}^r$. For a non-zero rational function $f$ in $\mathbb{K}(Y)$, one first looks at the vanishing order $a_r$ of $f$ at $Y_{r-1}$, then considers the function $f_{r-1}:=fu_r^{-a_r}\vert_{Y_{r-1}}$ to eliminate the zero or pole, and repeats this procedure for $f_{r-1}$ and the flag starting from $Y_{r-1}$. The outcome is a point $\nu_{Y_\bullet}(f):=(a_r,\ldots,a_1)\in\mathbb{Z}^r$; taking into account the total order on $\mathbb{Z}^r$, one obtains a valuation $\nu_{Y_\bullet}:\mathbb{K}(Y)\setminus\{0\}\to \mathbb{Z}^r$. 

Extending the valuation to the homogeneous coordinate ring $\mathbb{K}[Y]$ of $Y$ by sending a homogeneous function $f\in\mathbb{K}[Y]$ of degree $m$ to $(m,\nu_{Y_\bullet}(f))$ yields a valuation $\tilde{\nu}_{Y_\bullet}:\mathbb{K}[Y]\setminus\{0\}\to\mathbb{Z}\times \mathbb{Z}^r$ . The image of $\tilde{\nu}_{Y_\bullet}$ is a monoid. One of the most important questions in Newton-Okounkov theory is to determine when this monoid is finitely generated. If it happens to be so, Anderson \cite{And} obtains a toric degeneration\footnote{In this article, toric varieties are irreducible but not necessarily normal.} of $Y$ to the toric variety associated to this monoid. 

Motivated by seeking for an interpretation of the LS-paths in the above setup as vanishing order of functions, the second and the third author in \cite{FL} studied the case of Grassmann varieties. Instead of a flag of subvarieties, a web of subvarieties consisting of Schubert varieties is fixed. They constructed in \emph{loc.cit.} a quasi-valuation by choosing the minimum of all possible vanishing orders at each step. The graded algebra associated to the filtration arising from this quasi-valuation coincides with the discrete Hodge algebra in \cite{DEP} (a.k.a. the Stanley-Reisner algebra of the poset arising from the web). It was asked in \cite{FL, FaFoL} how to generalize this construction to Schubert varieties in a partial flag variety.

\subsection{Seshadri stratification, Newton-Okounkov bodies and standard monomial theory}

In this article we introduce the notion of a Seshadri stratification. In such a framework we construct a Newton-Okounkov simplicial complex and the associated semi-toric degeneration: this enables us to prove a formula on the degree of $X$ and to give a new geometric setup for standard monomial theory.

\subsubsection{Semi-toric degeneration from Seshadri stratification}

The geometric setting in the entire article is encoded in the concept of a Seshadri stratification (Definition \ref{Defn:SS}) of an embedded projective variety $X\subseteq\mathbb{P}(V)$\footnote{Such a setup is more suitable to Standard Monomial Theory, line bundles or linear systems will be discussed in a future work.}, where $V$ is a finite dimensional vector space. Such a stratification consists of a collection of subvarieties $X_p$ in $X$ which are smooth in codimension one, together with homogeneous functions $f_p$ on $V$, both indexed by a finite set $A$, i.e. $p\in A$. The set $A$ is naturally endowed with a poset structure from the inclusion of subvarieties, such that covering relation $q<p$\footnote{A relation $q<p$ in $A$ is called a \emph{covering relation}, if there is no $r\in A$ satisfying $q<r<p$.} means that $X_q$ is a divisor in $X_p$. The poset $A$ is assumed to have a maximal element $p_{\mathrm{max}}$ with $X_{p_{\mathrm{max}}}=X$. These subvarieties and the functions are compatible in the following sense: 
\begin{itemize}
\item[-] the vanishing set of the restriction of $f_p$ to $X_p$ is the union of all divisors in $X_p$ which are of form $X_q$;
\item[-] $f_p$ vanishes on $X_r$ for $p\not\leq r$.
\end{itemize}

Typical examples of this setting are Schubert varieties in a flag variety and extremal weight functions (Section \ref{Schubert}). More examples, varying from quadrics and elliptic curves to Grassmann varieties and group compactifications, will be discussed in Section \ref{backexamples}.

The requirements above seem to be restrictive. It is natural to ask for the existence and the uniqueness of Seshadri stratifications on an embedded projective variety. 

Concerning the existence: the definition of a Seshadri stratification demands the variety to be smooth in codimension one, and this is in fact sufficient:

\begin{IntroProp}[Proposition \ref{Prop:Generic}]\label{Prop1}
Every embedded projective variety $X\subseteq\mathbb{P}(V)$, smooth in codimension one, admits a Seshadri stratification.
\end{IntroProp}

The Seshadri stratification is far away from being unique: examples will be discussed in the article (Example \ref{examplep3}, Remark \ref{Rmk:Positroid}). 

One of the purposes of this article is to prove the following theorem, constructing semi-toric\footnote{Different to the terminology in symplectic geometry, a semi-toric variety is a variety whose irreducible components are toric varieties \cite{Cal}.} degenerations of $X$ from a Seshadri stratification on $X$.
\begin{IntroThm}[Theorem \ref{Thm:Degeneration}]\label{Thm1}
Let $X\subseteq\mathbb{P}(V)$ be a projective variety and $X_p,f_p$, $p\in A$ defines a Seshadri stratification on $X$. There exists a flat degeneration of $X$ into a reduced union of projective toric varieties $X_0$. Moreover, $X_0$ is equidimensional, and its irreducible components are in bijection with maximal chains in $A$.
\end{IntroThm}

Combining with Proposition \ref{Prop1} gives the following

\begin{IntroCor}[Corollary \ref{Cor:SemiToric}]\label{Cor1}
Every embedded projective variety, which is smooth in codimension one, admits a flat degeneration into a reduced union of projective toric varieties, the number of irreducible components coincides with its degree.
\end{IntroCor}

An important problem in the study of toric degenerations is to construct degenerations of projective varieties into projective toric varieties. The above theorem does not go precisely in this direction: our aim is rather to seek for degenerations of a projective variety which are compatible with a prescribed collection of subvarieties. Generally speaking, such a degeneration can not be toric, as being pointed out by Olivier Mathieu already for Schubert varieties (see the introduction of \cite{Cal}). The above theorem provides an answer to this problem if the projective variety admits a Seshadri stratification. In other words, such a degeneration exists if there are regular functions with prescribed set-theoretic vanishing locus. This condition is in the same vein as the Riemann-Roch theorem: the geometry gets controlled by the existence of certain functions.

The proof of the above theorem occupies a large part of the article. The general idea is similar to the one in \cite{And}, as soon as an analogue of a Newton-Okounkov polytope (not just a body) can be associated to the Seshadri stratification. In fact, we will construct a Newton-Okounkov simplicial complex from a Seshadri stratification, and the semi-toric variety $X_0$ is determined by this simplicial complex together with a lattice in each simplex.

This article is influenced by the idea of Allen Knutson \cite{Kn06} to use Rees valuations. Later in the work of Alexeev and Knutson \cite{AK}, they suggested to apply this idea to recover the degenerations in \cite{Chi} for Schubert varieties. 

\subsubsection{Newton-Okounkov simplicial complex}\label{Sec:IntroNOSC}

In a Seshadri stratification, all maximal chains in $A$ have the same length, which is $\dim X$. We start with a na\"{i}ve idea. Fix a maximal chain $\mathfrak{C}:p_r>p_{r-1}>\ldots>p_1>p_0$ in $A$, we aim to produce a convex body as in \cite{KK,LM} from the flag of subvarieties
$$X=X_{p_r}\supseteq X_{p_{r-1}}\supseteq\ldots \supseteq X_{p_1}\supseteq X_{p_0}$$
and the functions $f_{p_r},\ldots,f_{p_0}$. 

We switch to the affine picture: let $\hat{X}_{p_k}$ be the affine cone over $X_{p_k}$. The problem is that the restriction of $f_{p_k}$ to $\hat{X}_{p_k}$ is not necessarily a uniformizer in the local ring $\mathcal{O}_{\hat{X}_{p_{k-1}},\hat{X}_{p_k}}$. As in the first step of the construction of a valuation associated to the flag, for a rational function $g\in\mathbb{K}(\hat{X})$, there is no reason why there exists $m\in\mathbb{Z}$ such that the function $gf_p^m$, when restricted to $\hat{X}_{p_{r-1}}$, yields a well-defined and non-zero rational function. One could think about choosing a uniformizer in $\mathcal{O}_{\hat{X}_{p_{k-1}},\hat{X}_{p_k}}$, however, we will not concentrate on only one maximal chain $\mathfrak{C}$ but take into account all of them, there is no control of this uniformizer on the flag associated to other maximal chains in $A$.

We adjust the construction of the Newton-Okounkov body by keeping track of the vanishing multiplicities of the functions $f_p$ along a divisor. We consider the Hasse graph of the poset $A$ and colour an edge arising from the covering relation $q<p$ by the vanishing order of $f_p$ on $X_q$ (see Section \ref{Hasse}). These colours are called bonds. 

We fix $N$ to be the l.c.m. of all bonds appearing in the coloured Hasse graph.

For a non-zero rational function $g_r:=g\in\mathbb{K}(\hat{X})$ with vanishing order $a_r$ along the divisor $\hat{X}_{p_{r-1}}$ in $\hat{X}=\hat{X}_{p_r}$, we define a rational function 
$$h:=\frac{g_r^N}{f_{p_r}^{N\frac{a_r}{b_r}}}\in\mathbb{K}(\hat{X}_{p_{r}}),$$
where $b_r$ is the vanishing order of $f_{p_r}$ along $\hat{X}_{p_{r-1}}$. The restriction of $h$ to $\hat{X}_{p_{r-1}}$, denoted by $g_{r-1}$, gives rise to a well-defined non-zero rational function in $\mathbb{K}(\hat{X}_{p_{r-1}})$ (Lemma \ref{one}). This procedure can be henceforth iterated, yielding a sequence of rational functions $g_{\mathfrak{C}}:=(g_r,g_{r-1},\ldots,g_0)$ with $g_k\in\mathbb{K}(\hat{X}_{p_k})\setminus\{0\}$. The vanishing order of $g_k$ (resp. $f_{p_k}$) on $\hat{X}_{p_{k-1}}$ will be denoted by $a_k$ (resp. $b_k$).

Similar to the Newton-Okounkov theory, the vanishing orders will be collected to define a valuation. In view of the $N$-th powers appearing in the sequence of rational functions, we define a map $\mathcal{V}_{\mathfrak{C}}:\mathbb{K}[\hat{X}]\setminus\{0\}\to\mathbb{Q}^\mathfrak{C}$ in the following way:
$$g\mapsto \frac{a_r}{b_r}e_{p_r}+\frac{1}{N}\frac{a_{r-1}}{b_{r-1}}e_{p_{r-1}}+\ldots+\frac{1}{N^r}\frac{a_{0}}{b_{0}}e_{p_0},$$
where $e_{p_k}$ is the coordinate function in $\mathbb{Q}^\mathfrak{C}$ corresponding to $p_k\in\mathfrak{C}$.
Such a map is indeed a valuation (Proposition \ref{VCvaluation}) having at most one-dimensional leaves (Theorem \ref{chainleavedimensiontwo}). As in the situation of Section \ref{Sec:IntroNOB}, we do not know whether the image of $\mathcal{V}_{\mathfrak{C}}$ is a finitely generated monoid. In general,
the finite generation property is not expected in general as the flag of subvarieties reveals rather the local geometry. 

In order to pass from local to global, we define a quasi-valuation (Definition \ref{quasidef}) $\mathcal{V}:\mathbb{K}[\hat{X}]\setminus\{0\}\to\mathbb{Q}^A$ by taking the minimum over all maximal chains in $A$. For this we choose a total order $>^t$ on $A$ refining the partial order (Equation \ref{totorder}), extend lexicographically to $\mathbb{Q}^A$, and define 
$$\mathcal{V}(g):=\min\{\mathcal{V}_{\mathfrak{C}}(g)\mid \mathfrak{C}\text{ is a maximal chain in }A\},$$
where $\mathbb{Q}^{\mathfrak{C}}$ is naturally embedded into $\mathbb{Q}^A$.

The first nice property of this quasi-valuation is its positivity: the image of $\mathcal{V}$ is contained in $\mathbb{Q}^A_{\geq 0}$ (Proposition \ref{positivity}). Such a property is guaranteed by the Valuation Theorem of Rees (Theorem \ref{ValuationTheoremRees}), whose spirit is already incorporated as part of the Seshadri stratification, as well as the valuation $\mathcal{V}_{\mathfrak{C}}$. This positivity encodes in fact the regularity: for a non-zero regular function $g\in\mathbb{K}[\hat{X}]$ and a maximal chain $\mathfrak{C}$ on which the minimum of $\mathcal{V}(g)$ is attained, the function $g_k$ in the sequence of rational functions $g_{\mathfrak{C}}$ is regular in the normalization of $\hat{X}_{p_{k}}$ for all $k=0,1,\ldots,r$.

For such a function $g\in\mathbb{K}[\hat{X}]$, there could be many maximal chains on which the minimum $\mathcal{V}(g)$ is attained. We will prove (Proposition \ref{supportcondition}) that these maximal chains are precisely those containing the support of $\mathcal{V}(g)$, defined as the set of elements in $A$ on which $\mathcal{V}(g)$ takes non-zero (hence positive) value. As a consequence of this characterization via supports, we are able to decompose the image $\Gamma$ of the quasi-valuation $\mathcal{V}$ into a finite union of (finitely generated) monoids $\Gamma_{\mathfrak{C}}$ (Corollary \ref{finiteunionmonoid}) where $\mathfrak{C}$ runs over all maximal chains in $A$ and $\Gamma_{\mathfrak{C}}$ consists of elements in $\Gamma$ supported on $\mathfrak{C}$.

The set $\Gamma$ encodes rather the global aspects of $X$: first, for a given regular function, it tells how to smooth out its zeros using the functions $f_p$ and keeping the regularity simultaneously; secondly, the quasi-valuation $\mathcal{V}$ has at most one-dimensional leaves (Lemma \ref{NuLeaves}), hence $\mathbb{K}[\hat{X}]$ and $\Gamma$ have the same ``size''; moreover, the monoids $\Gamma_{\mathfrak{C}}$ are finitely generated (Lemma \ref{coreminimaldifference}).
 
The finite generation of $\Gamma_{\mathfrak{C}}$ allows us to investigate the geometry of $\Gamma$. We will define a fan algebra $\mathbb{K}[\Gamma]$ by gluing different $\Gamma_{\mathfrak{C}}$ in a Stanley-Reisner way (Definition \ref{Defn:FanAlgebra}). The affine variety $\mathrm{Spec}(\mathbb{K}[\Gamma])$ associated to the fan algebra is an irredundant union of affine toric varieties $\mathrm{Spec}(\mathbb{K}[\Gamma_{\mathfrak{C}}])$ where $\mathfrak{C}$ runs over all maximal chains in $A$, each of dimension $\dim \hat{X}$ (Proposition \ref{Prop:DecompsitionIrrComp}). 

In order to prove Theorem \ref{Thm1}, we need to construct a flat family over $\mathbb{A}^1$ with special fibre $\mathrm{Proj}(\mathbb{K}[\Gamma])$. The quasi-valuation $\mathcal{V}$ induces an algebra filtration on $R:=\mathbb{K}[\hat{X}]$. The associated graded algebra $\mathrm{gr}_{\mathcal{V}}R$ is finitely generated and reduced (Corollary \ref{all_finite-generated}). Different to the toric case as in \cite{And}, some work is needed in proving that the fan algebra $\mathbb{K}[\Gamma]$ and the associated graded algebra $\mathrm{gr}_{\mathcal{V}}R$ are indeed isomorphic as algebras (Theorem \ref{fanAndDegeneratetheorem}). Once this isomorphism is established, the machinery of the Rees algebra associated to a filtration can be applied to construct the flat family in Theorem \ref{Thm1}. 

As an application to Theorem \ref{Thm1}, we show that if the poset $A$ is Cohen-Macaulay over $\mathbb{K}$ and the monoids $\Gamma_{\mathfrak{C}}$ are saturated, then the embedded projective variety is projectively normal (Theorem \ref{Thm:ProjNormal}).

Out of the set $\Gamma$ we define the associated Newton-Okounkov simplicial complex $\Delta_{\mathcal{V}}$ (Definition \ref{Defn:NOSC}), where each monoid $\Gamma_{\mathfrak{C}}$ for a maximal chain $\mathfrak{C}$ in $A$ contributes a simplex. Different to the Newton-Okounkov theory \cite{LM,KK}, where the collection of rational functions are uniformizers, to each simplex we associate a natural lattice $\mathcal{L}^{\mathfrak{C}}$ from the quasi-valuation. The degree $\mathrm{deg}(X)$ of the embedded projective variety $X\hookrightarrow\mathbb{P}(V)$ can be computed as a sum of volumes of simplexes:

\begin{IntroThm}[Theorem \ref{volumetheorem1}]
For each maximal chain $\mathfrak{C}$, we provide an $r$-dimensional simplex with rational vertices $D_{\mathfrak{C}}$ such that
$$\mathrm{deg}(X)=r!\sum\mathrm{vol}(D_\mathfrak{C}),$$
where the sum runs over all maximal chains $\mathfrak{C}$ in $A$.
\end{IntroThm}

For Schubert varieties, such a formula was obtained by Knutson \cite{Kn99} in the symplectic geometric setting, and the first author \cite{Chi} in the algebro-geometric setting. 

When the monoids $\Gamma_{\mathfrak{C}}$ are saturated, the Hilbert function can be calculated in the same way as Ehrhart functions of simplexes (Corollary \ref{volumetheorem2}). 

\subsubsection{The case of a totally ordered poset} 

In order to help the reader compare our approach with the usual Netwon-Okounkov context, we want to shortly discuss the case of a Seshadri stratification with a totally ordered poset $A = \{p_r > p_{r-1} > \cdots > p_1 > p_0\}$. The quasi-valuation $\mathcal{V}$ coincides with the valuation $\mathcal{V}_A$ for the unique maximal chain $A$.

Further, the zero locus of the extremal function $f_k$ in $\hat{X}_{p_k}$ is the divisor $\hat{X}_{p_{k-1}}$. It is then clear that $f_k|_{\hat{X}_{p_k}}$ is the $b_k$--th power of a uniformizer $u_k$ in $\mathcal{O}_{\hat{X}_{p_{k}},\hat{X}_{p_{k-1}}}$. In particular the first $r$ components of $\mathcal{V}(g)\in\mathbb{Q}^{r+1}$, for homogeneous $g\in\mathbb{K}[\hat{X}]\setminus\{0\}$, are just renormalizations of the components of the usual valuation $\nu_{X_\bullet}(g)\in\mathbb{Z}^r$ associated to the flag of subvarieties $X_\bullet = (X=X_{p_r}\supseteq X_{p_{r-1}}\supseteq\cdots\supseteq X_{p_0})$ and uniformizer $u_k$, $k = 1,\ldots,r$, in the Newton-Okounkov theory. More precisely: if $\nu_{X_\bullet}(g) = (n_r,\ldots,n_1)$ then
\[
\mathcal{V}(g) = \left(\frac{n_r}{b_r}, \frac{n_{r-1}}{b_{r-1}},\ldots,\frac{n_0}{b_0}\right)
\]
where $n_0$ is such that $\sum_{j=0}^r\deg f_j \frac{n_j}{b_j} = \deg g$.

Finally, the Newton-Okounkov simplicial complex of $\mathcal{V}$ is just a simplex in this case.

If such a setup arises from a generic hyperplane stratification (see Section \ref{Sec:Generic}), the finite generation result in this article is closely related to \cite[Proposition 14]{AKL}. Indeed, the assumptions in that proposition allow to construct a Seshadri stratification with a linear poset, which is exactly the proof in \emph{loc.cit}.

\subsubsection{Standard monomial theory}\label{Sec:IntroSMT}

Seshadri stratifications provide geometric setups for standard monomial theories on $R=\mathbb{K}[\hat{X}]$.

In view of the calculation of Hilbert functions (Proposition \ref{volumetheorem2}), in order to have a standard monomial theory, the monoids $\Gamma_{\mathfrak{C}}$ should be assumed to be saturated. Under this assumption, each element $\underline{a}\in\Gamma_{\mathfrak{C}}$ can be uniquely decomposed into a sum of indecomposable elements (Definition \ref{Defn:Indec}) in $\Gamma_{\mathfrak{C}}$ (Proposition \ref{Prop:UniqueDec}). 

From the one-dimensional leaf property of the quasi-valuation $\mathcal{V}$, for each indecomposable element $\underline{a}\in\Gamma_{\mathfrak{C}}$ we choose a regular function $x_{\underline{a}}\in R$ with $\mathcal{V}(x_{\underline{a}})=\underline{a}$. The condition of being standard will be defined on monomials in these regular functions: a monomial $x_{\underline{a}_1}\cdots x_{\underline{a}_k}$ with $\underline{a}_1,\ldots,\underline{a}_k\in\Gamma_{\mathfrak{C}}$ is standard if for $i=1,\ldots,k-1$, $\min\supp\underline{a}_i\geq \max\supp\underline{a}_{i+1}$. This defines a standard monomial theory on $R$ (Proposition \ref{proposition_standard_monomial_basis}).

One of the tasks in Seshadri's paradigm is to construct standard monomial bases compatible with all strata $X_p$. We call a standard monomial $x_{\underline{a}_1}\cdots x_{\underline{a}_k}$ standard on $X_p$ if the maximal element in $\mathrm{supp}\,\underline{a}_1,\ldots,\mathrm{supp}\,\underline{a}_k$ is $\leq p$. This compatibility requires extra conditions on the independence to the choice of the total order $>^t$ in the definition of $\mathcal{V}$. To eliminate this dependency, we introduce the \emph{balanced} conditions on Seshadri stratifications (Definition \ref{definition_balanced_statification}); this extra structure allows us to show:

\begin{IntroThm}[Theorem \ref{prop:SMT:for:subvarieties}]\label{IntroThm:3}
In the above situation, the following hold:
\begin{enumerate}
    \item[i)] All standard monomials on $X$ which are standard on $X_p$ form a basis of $\mathbb{K}[\hat{X}_p]$.
    \item[ii)] Standard monomials on $X$ which are not standard on $X_p$ are precisely those vanishing on $X_p$. They form a linear basis of the defining ideal of $X_p$ in $X$.
    \item[iii)] For any $p,q\in A$, the scheme-theoretic intersection $X_p\cap X_q$ is a reduced union of strata contained in both. 
\end{enumerate}
\end{IntroThm}

\subsubsection{L-S paths as vanishing orders}

We explain to what extent the framework in this article answers Seshadri's question in Section \ref{Sec:IntroPathSMT}.

Let $G$ be a simple simply connected algebraic group, $B$ a Borel subgroup, $T$ a maximal torus and $W$ the Weyl group of $G$, viewed as a poset with the Bruhat order. The assumption on $G$ is made only to simplify the statements: the results will be proved in \cite{CFL} for Schubert varieties in a symmetrizable Kac-Moody group. The Schubert varieties $X(\sigma)$, together with the extremal weight functions $p_\sigma$ for $\sigma\in W$, form a Seshadri stratification of the flag variety $X:=G/B$, embedded in $\mathbb{P}(V(\lambda))$ for a regular dominant weight $\lambda$. 

To describe the image of the associated quasi-valuation $\mathcal{V}$, we introduce for each maximal chain $\mathfrak{C}$ an explicit lattice $L_{\mathfrak{C},\lambda}$ (Equation \ref{EqLSLattice}) and a fan of (saturated) monoids $L_\lambda^+$, which is in an easy bijection with the Lakshmibai-Seshadri paths (LS-paths) of shape $\lambda$ (Lemma \ref{LSlattice}). 

\begin{IntroThm}[Theorem \ref{flag:standard:monomial}, Section \ref{Sec:Schubert}, Proposition \ref{Prop:degreebonds}]\par\noindent
\begin{enumerate}
    \item[i)] The image of $\mathcal{V}$ coincides with $L_\lambda^+$.
    \item[ii)] The degree of $X(\sigma)$ is a sum of products of bonds.
    \item[iii)] The Seshadri stratification is balanced, hence all statements in Section \ref{Sec:IntroSMT} hold for Schubert varieties too.
    \item[iv)] The Schubert varieties $X(\sigma)\hookrightarrow\mathbb{P}(V(\lambda))$ are projectively normal.
\end{enumerate}
\end{IntroThm}

The difficult part is i), and the key point to the proof is Theorem \ref{valuationofpathvector}: for any element $\pi$ in $L_\lambda^+$ (looked as an LS-path) we seek for a regular function $p_\pi$ with $\mathcal{V}(p_\pi)=\pi$. A candidate for $p_\pi$ has already been constructed by the third author in \cite{L1}, but to prove the desired vanishing properties requires results and techniques from representation theory of algebraic groups and quantum groups at roots of unity. The complete proof is given in a separate article \cite{CFL}.

As a consequence of the theorem, the LS-paths get interpreted as vanishing orders of regular functions. This fits perfectly into Seshadri's expectation of ``cellular Riemann-Roch''.

\subsection{Outline of the article}
In Section~\ref{A_partially_ordered_collection} we introduce the Seshadri stratifications of an embedded projective variety and the associated Hasse graph coloured by bonds. A few quick examples are discussed therein as running examples for the article. In Section~\ref{ValuationsAndquasi-valuations} we collect a few standard facts about valuations and quasi-valuations, and, in particular, we recall the homogenized quasi-valuation arising from ideal filtration, and the Rees valuation theorem investigating their structures.
In Section~\ref{Section:ValuationCodimOne} and \ref{chainsandfunctions} we prepare the procedure used in 
Section~\ref{value_valuations} to define a valuation associated to a maximal chain in the poset $A$. Further studies regarding these valuations are carried out in Section~\ref{Sec:FiniteGen}.
In Section~\ref{nonegativequasival} we introduce the main point of this article: a quasi-valuation, defined
as the minimum over the collection of valuations introduced in Section~\ref{value_valuations}. 

We introduce in Section~\ref{fanandcone} the notion of fan monoids and fan algebras to describe the associated graded algebra. In Section~\ref{leavesandgr} we discuss some nice properties of this quasi-valuation. Section~\ref{Sec:FanDeg} is devoted to proving that the associated graded algebra and the fan algebra are isomorphic as algebras, which is the crucial step in the construction of the semi-toric degeneration in Section~\ref{flatdegen}.

In Section~\ref{Nocomplex} we associate to the Seshadri stratification a Newton-Okounkov simplicial complex to investigate the discrete geometry behind the semi-toric variety. This complex, unlike the usual Newton-Okounkov setting, is not necessarily a convex body. We compensate this difference by endowing the complex with a rational or integral structure. As an application, we prove a criterion on the projective normality in Section \ref{Sec:ProjNormal}. Under certain hypothesis on the Seshadri stratification, we define a standard monomial theory for the homogeneous coordinate ring in Section \ref{section_standard_monomial_theory}.

Examples such as Schubert varieties in partial flag varieties, compactification of torus and $\mathrm{PSL}_2(\mathbb{C})$, quadrics and elliptic curves are discussed in Section~\ref{backexamples}.

The structure of the article is rather linear, the notations used throughout the article are gathered in a list of notations after Section~\ref{backexamples}.

\subsection{Recent development}
The Seshadri stratification of a Schubert variety consisting of its Schubert subvarieties is studied in \cite{CFL}, where results announced in Section~\ref{Schubert} are proved with the help of quantum groups at roots of unity. A different approach without using quantum groups is given in \cite{CFL4}. The algebraic counterpart of this article is studied in \cite{CFL2} in the framework of valuations on LS-algebras, the connection to the current article is made clear in \cite{CFL3}. More results on normal Seshadri stratification, especially its connection to Gr\"obner theory, are topics of \emph{loc.cit}.

\subsection{Acknowledgements}
The work of R.C. is partially supported by PRIN \\2017YRA3LK\_005 ``Moduli and Lie Theory''. The work of X.F. and P.L. is partially supported by DFG SFB/Transregio 191 ``Symplektische Strukturen in Geometrie, Algebra und Dynamik''. We thank Michel Brion, Ghislain Fourier, Kiumars Kaveh, Andrea Maffei, Luca Migliorini and Christopher Manon for helpful discussions in different stages of this work. We would like to thank the referee, whose comments help us to simplify various proofs.

\section{Seshadri stratifications}\label{A_partially_ordered_collection}

\subsection{Conventions}

Let $V$ be a finite dimensional $\mathbb{K}$-vector space. For a homogeneous polynomial function $f\in\mathrm{Sym}(V^*)$, we denote its vanishing set $\mathcal{H}_f:=\{[v]\in\mathbb{P}(V)\mid f(v)=0\}$. For a partially ordered set (poset) $(A,\leq)$ and $p\in A$, we denote $\gls{Ap}:=\{q\in A\mid q\leq p\}$: $(A_p,\leq)$ is a poset. A relation $q<p$ in $A$ is called a covering relation, if there is no $r\in A$ such that $q<r<p$.

In this article, varieties are assumed to be irreducible. Toric varieties are not necessarily normal. When it is necessary to emphasize on the normality, we use the term ``normal toric variety''.

\subsection{Definition and examples}\label{setup}

Let $X\subseteq \mathbb P(V)$ be an embedded projective variety with graded
homogeneous coordinate ring $R=\mathbb K[X]$. The degree $k$ component of $R$ will be denoted by $R(k)$: $R=\bigoplus_{k\ge 0} R(k)$.

Let $\{\gls{Xp}\mid p\in A\}$ be a collection of projective subvarieties $X_p$ in $X$ indexed by a finite set $\gls{A}$. The set $\gls{A}$ is naturally endowed with a partial order $\leq$ by: for $p,q\in A$, $p\leq q$ if and only if $X_p\subseteq X_q$. We assume that there exists a unique maximal element $p_{\max}\in A$ with $X_{p_{\max}}=X$.

For each $p\in A$, we fix a homogeneous function $\gls{fp}$ on $V$ of degree larger or equal to $1$.

\begin{definition}\label{Defn:SS}
The collection of subvarieties $\gls{Xp}$ and homogeneous functions $\gls{fp}$ for $p\in A$ is called a \emph{Seshadri stratification}, if the following conditions are fulfilled:
\begin{enumerate}
\item[(S1)] the projective varieties $X_p$, $p\in A$, are all smooth in codimension one; for each covering relation $q<p$ in $A$, $X_q\subseteq X_p$ is a codimension one subvariety; 
\item[(S2)] for any $p\in A$ and any $q\not\leq p$, $f_q$ vanishes on $X_p$;
\item[(S3)] for $p\in A$, it holds: set theoretically 
$$\mathcal{H}_{f_p}\cap X_p=\bigcup_{q\text{ covered by }p} X_q.$$
\end{enumerate}
The subvarieties $X_p$ will be called \emph{strata}, and the functions $\gls{fp}$ are called \emph{extremal functions}.
\end{definition}

The following lemma will be used throughout the article often without mention.

\begin{lemma}\label{Lem:SS}
\begin{enumerate}
    \item[i)] The function $f_p$ does not identically vanish on $X_p$.
    \item[ii)] All maximal chains in $A$ have the same length, which coincides with $\dim X$. In particular, the poset $A$ is graded.
    \item[iii)] The intersection of two strata is a union of strata; in particular for each $p,q\in A$ we have
\[
X_p \cap X_q = \bigcup_{t\leq p,q} X_t.
\]
\end{enumerate}
\end{lemma}

\begin{proof}
If $X_p$ is just a point, then (S3) implies the intersection $\mathcal{H}_{f_p}\cap X_p$ is empty, which implies i) in this case.
If $X_p$ is not just a point, then the intersection $\mathcal{H}_{f_p}\cap X_p$ is not empty, and (S1) and (S3) enforce the intersection to be a union of divisors. In particular: $X_p\not\subseteq \mathcal{H}_{f_p}$, which implies i),
and there must exist elements $q$ in $A$ covered by $p$, which implies ii).

From the definition of the partial order on $A$ it follows that $\bigcup_{t\leq p,q} X_t\subseteq X_p\cap X_q$. We prove by induction on the length of a maximal chain joining $p$ with a minimal element in $A$ that the intersection $X_p\cap X_q$ is a union of strata. Such a length is well-defined by the part ii).

When $p = p_0$ is a minimal element in $A$, it follows that either $p_0 \leq q$ and $X_{p_0} \cap X_q = X_{p_0}$, or $p_0 \not\leq q$ and $X_{p_0} \cap X_q = \emptyset$; in both cases the claim is proved.

For an arbitrary $p\in A$, if $p \leq q$ then $X_p \cap X_q = X_p$, so we can assume that $p \not\leq q$; hence $f_p|_{X_q} = 0$ by (S2). In particular $f_p|_{X_p\cap X_q} = 0$. But, for $x\in X_p$, (S3) implies that $f_p(x) = 0$ if and only if $x\in \bigcup_{p'\leq p}X_{p'}$, which gives the inclusion $X_p\cap X_q \subseteq \bigcup_{p'<p}(X_{p'}\cap X_q)$. Since the reverse inclusion clearly holds, we have proved
\[
X_p\cap X_q = \bigcup_{p'<p}(X_p' \cap X_q).
\]
By induction, each intersection $X_{p'}\cap X_q$ is a union of strata, hence $X_p\cap X_q$ is a union of strata.
\end{proof}

Thanks to the part ii) of the lemma we can define the length of an element in $A$.

\begin{definition}\label{length}
Let $p\in A$. The \emph{length} $\gls{l}(p)$ of $p$ 
is the length of a (hence any) maximal chain joining $p$ with a minimal element in $A$.
\end{definition} 

It is clear that $\ell(p)=\dim X_p$.

\begin{rem}\label{induction}
For a fixed $p\in A$, by Lemma~\ref{Lem:SS} (ii), the poset $\gls{Ap}$ has a unique maximal element, and all maximal chains have the same length.
The collection of varieties $X_q$, $q\in A_p$, and the extremal functions 
$f_q$, $q\in A_p$ satisfy the conditions (S1)-(S3), and hence defines a Seshadri stratification for $X_p\hookrightarrow  \mathbb P(V)$.
\end{rem}

Before going further we look at some examples of Seshadri stratifications. Further examples will be given in Section~\ref{backexamples}.

\begin{example}\label{grasstwofour1}
Let $\{e_1,e_2,e_3,e_4\}$ be the standard basis of $\mathbb K^4$. The wedge products $e_i\wedge e_j$, $1\le i<j\le 4$, form a basis of $\bigwedge^2 \mathbb K^4$. Denote the indexing set of the basis by $I_{2,4}$, it consists of pairs of positive numbers $(i,j)$, strictly increasing, and 
smaller or equal to $4$. The corresponding elements $\{x_{i,j}\mid  (i,j)\in I_{2,4}\}$ of the dual basis are called Pl\"ucker coordinates. 

The set $I_{2,4}$ is endowed with a partial order: $(i,j)\le(k,\ell)$ if and only if $i\le k$ and $j\le \ell$. Let $X:=\mathrm{Gr}_2\mathbb K^4\subseteq \mathbb P(\bigwedge^2 \mathbb K^4)$ be the Grassmann variety
of $2$-planes in $\mathbb K^4$, emdedded into $\mathbb P(\bigwedge^2 \mathbb K^4)$ via the Pl\"ucker embedding. 
Set theoretically, the Schubert varieties $X(i,j)\subseteq \mathrm{Gr}_2\mathbb K^4$ for $(i,j)\in I_{2,4}$ are defined by
$$
X(i,j):=\{[v]\in \mathrm{Gr}_2\mathbb K^4\mid x_{k,\ell}([v])=0\ \text{for all}\ (k,\ell)\in I_{2,4}\ \text{such that}\ (k,\ell)\not\leq (i,j)\}.
$$
The collection of subvarieties $X(i,j)$, $(i,j)\in I_{2,4}$, together with the functions $f_{(i,j)}:=x_{i,j}$, $(i,j)\in  I_{2,4}$, define a Seshadri stratification on $\mathrm{Gr}_{2}\mathbb{K}^4$.

Below the Hasse diagram showing the inclusion relations between the Schubert varieties, here
$X(i,j)\rightarrow X{(k,\ell)}$ means $X{(i,j)}$ is contained in $X{(k,\ell)}$ of codimension one. 
It depicts meanwhile the Hasse diagram of the partial order on $I_{2,4}$.
$${ \tiny
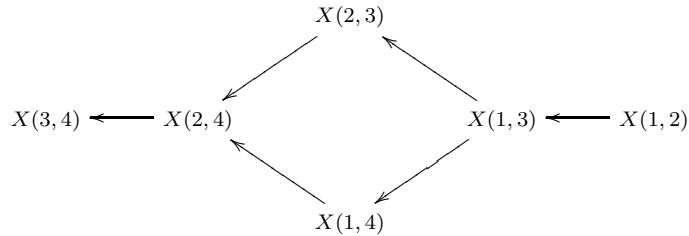
\xymatrix{
&&X(2,3)  \ar[dl]&\\
X(3,4)&\ar[l]X(2,4) &&X(1,3)\ar[ul] \ar[dl] &X(1,2)\ar[l]\\
&&X(1,4) \ar[ul]&}}
$$
\end{example}

A Seshadri stratification of $\mathrm{Gr}_2\mathbb{C}^4$ is not necessarily given by Schubert varieties, see Remark \ref{Rmk:Positroid}.

\begin{example}\label{GmodBexample}
More generally, consider the Grassmann variety $X:=\gls{Grdn}\subseteq \mathbb P(\bigwedge^d \mathbb K^n)$ of $d$-dimensional subspaces in $\mathbb K^n$, embedded into $ \mathbb P(\bigwedge^d \mathbb K^n)$ via the Pl\"ucker embedding. The Schubert varieties in $\mathrm{Gr}_d \mathbb{K}^n$ are indexed by the set $I_{d,n}:=\{\underline{i}=(i_1,\ldots,i_d)\mid 1\le i_1<\ldots<i_d\le n\}$ of strictly increasing sequences of length $d$. For $\underline{i}\in I_{d,n}$ let $X(\underline{i})$ denote the corresponding Schubert variety in $\mathrm{Gr}_d\mathbb{K}^n$. The partial order on $I_{d,n}$ induced by the inclusion of subvarieties coincides with the usual partial order on $I_{d,n}$: $\underline i\le \underline j$ if and only if $i_1\le j_1$, $\ldots$, $i_d\le j_d$. For $\underline{i}\in I_{d,n}$ let $f_{\underline{i}}:=x_{\underline{i}}$ be the Pl\"ucker coordinate. The collection of Schubert varieties $X(\underline{i})$, $\underline{i}\in I_{d,n}$ and the Pl\"ucker coordinates $x_{\underline{i}}$, $\underline{i}\in I_{d,n}$ define a Seshadri stratification on $X$.
\end{example} 

\begin{example}\label{examplep3}
Even on simple varieties such as projective spaces, there may exist \textit{several}
non-trivial Seshadri stratifications. One has been given in Example \ref{GmodBexample} as the special case $d=1$. We define now a different stratification.

We consider the projective space $X:=\mathbb P(\mathbb K^3)$. By fixing the standard basis $\{e_1,e_2,e_3\}$ of $\mathbb{K}^3$, the homogeneous coordinate ring $\mathbb{K}[X]$ can be identified with $\mathbb{K}[x_1,x_2,x_3]$. 

As indexing set we take the power set $A:=\mathfrak P(\{1,2,3\})\setminus \emptyset$ by omitting the empty set. As the collection of subvarieties we set for a subset $p\in A$: $X_p:=\mathbb P(\langle e_j\mid j\in p\rangle_\mathbb{K})$. The poset structure on $A$ induced by the inclusion of subvarieties coincides with that arising from inclusion of sets. For $p\in A$, let $f_p=\prod_{j\in p}x_j$ be the extremal function.

We leave it to the reader to verify that the collection of subvarieties $X_p$, together with the monomials $f_p$, $p\in A$, define a Seshadri stratification on $X$. Below the inclusion diagram of the varieties, and, in the same scheme, the functions $f_p$ corresponding to the subvariety $X_p$. 
{\tiny
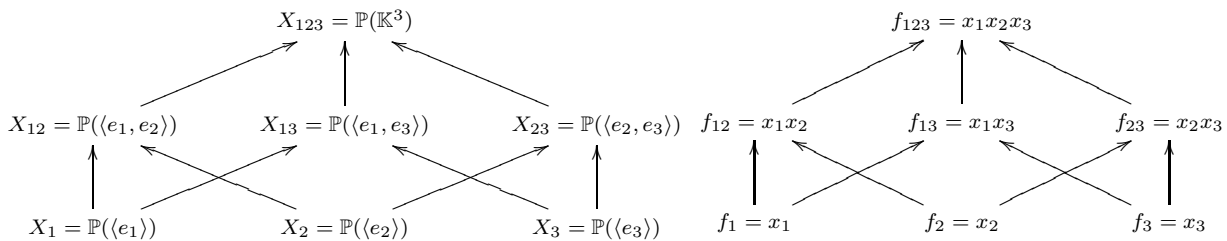
$$
\hskip -20pt
\xymatrix{ 
&X_{123}=\mathbb P(\mathbb K^3) \ar@{<-}[dl]\ar@{<-}[d] \ar@{<-}[dr]& \\
X_{12}=\mathbb P(\langle e_1,e_2\rangle) \ar@{<-}[d] \ar@{<-}[dr]&X_{13}=\mathbb P(\langle e_1,e_3\rangle ) \ar@{<-}[dr]\ar@{<-}[dl]
&X_{23}=\mathbb P(\langle e_2,e_3\rangle)\ar@{<-}[d] \ar@{<-}[dl]\\
X_1=\mathbb P(\langle e_1\rangle) 
&X_2=\mathbb P(\langle e_2\rangle) 
&X_3=\mathbb P(\langle e_3\rangle)\\
}
\xymatrix{
&f_{123}=x_1x_2x_3 \ar@{<-}[dl]\ar@{<-}[d] \ar@{<-}[dr]& \\
f_{12}=x_1x_2 \ar@{<-}[d] \ar@{<-}[dr]&f_{13}=x_1x_3 \ar@{<-}[dr]\ar@{<-}[dl]&f_{23}=x_2x_3\ar@{<-}[d] \ar@{<-}[dl]\\
f_{1}=x_1 
&f_{2}=x_2
&f_3=x_3 \\
}
$$}
\end{example}

\subsection{A Hasse diagram with bonds}\label{Hasse}

For a given Seshadri stratification of a projective variety $X$ consisting of subvarieties $X_p$ and extremal functions $f_p$ for $p\in A$, we associate to it in this section an edge-coloured directed graph.

Let $\gls{GA}$ be the Hasse diagram of the poset $A$. The edges in $\mathcal G_A$ are covering relations in $A$. If $p$ covers $q$, then the affine cone $\hat X_q$ of $X_q$ is a prime divisor in the affine cone $\gls{Xphat}$ of $X_p$. We denote by $\gls{bpq}\ge 1$ the vanishing multiplicity of $f_p$ at the prime divisor $\hat X_q$ (see Section~\ref{ValutaionDefn} for the definition of the vanishing multiplicity),
it is called the \textit{bond} between $p$ and $q$. The Hasse diagram with bonds is the diagram with edges coloured with the corresponding bonds: 
$q \stackrel{\gls{bpq}}{\longrightarrow} p$.

\begin{example}\label{grasstwofour}
For  $X=\mathrm{Gr}_2\mathbb K^4\subseteq \mathbb P(\bigwedge^2 \mathbb K^4)$ as in Example~\ref{grasstwofour1},
the corresponding Hasse diagram with bonds is:
$$ {\tiny
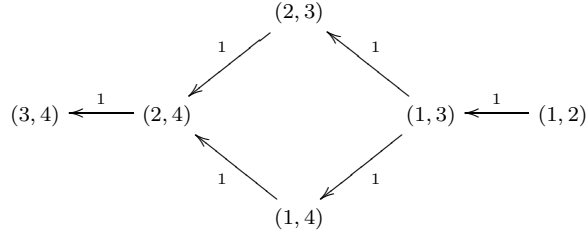
\xymatrix{
&&(2,3)  \ar@{<-}[dr]^1&\\
(3,4)\ar@{<-}[r]^1&(2,4) \ar@{<-}[dr]_1 \ar@{<-}[ur]^1&&(1,3) \ar@{<-}[r]^1&(1,2)\\
&&(1,4) \ar@{<-}[ur]_1&}}
$$
More generally, for the Grassmann variety $\mathrm{Gr}_d\mathbb{K}^n$ in Example~\ref{GmodBexample}, all the bonds are $1$.
\end{example}

As we will see in the following example, the bonds in a Seshadri stratification are not necessarily one.

\begin{example}\label{fullflag}
To avoid technical details in small characteristics, we assume in this example that the characteristic of $\mathbb{K}$ is zero or a large prime number.

Let $\mathrm{SL}_3$ be the group of $3\times 3$ matrices over $\mathbb{K}$ having determinant $1$. Its Lie algebra $\mathfrak{sl}_3$ consists of traceless $3\times 3$ matrices over $\mathbb{K}$. Let $B\subseteq\mathrm{SL}_3$ be the subgroup consisting of upper triangular matrices.

The group $\mathrm{SL}_3$ acts linearly on its Lie algebra $\mathfrak{sl}_3$ via the adjoint representation: for $g\in \mathrm{SL}_3$ and $M\in\mathfrak{sl}_3$, $g\cdot M:=gMg^{-1}$. 
By choosing $M$ to be a root vector for the highest root, this action induces an embedding $\mathrm{SL}_3/B\hookrightarrow \mathbb P(\mathfrak{sl}_3)$.

Let $\mathrm{S}_3$ be the Weyl group of $\mathrm{SL}_3$: it is the symmetric group acting on three letters. By abuse of notation we
identify $\sigma\in \mathrm{S}_3$ with an appropriately chosen representative $\sigma\in \mathrm{SL}_3$; it is, up to the sign of the entries, the corresponding permutation matrix.

In the Bruhat decomposition $\mathrm{SL}_3=\bigsqcup_{\sigma\in\mathrm{S}_3}B\sigma B$ of $\mathrm{SL}_3$, the class of the closure of each cell in $\mathrm{SL}_3/B$
$$X(\sigma):=\overline{B\sigma B}/ B\subseteq\mathrm{SL}_3/ B$$
is the Schubert variety associated to $\sigma\in\mathrm{S}_3$.
We fix a basis of $\mathfrak{sl}_3$ as follows:
$$v_{(13)}:=E_{3,1},\ v_{(132)}:=E_{3,2},\ v_{(123)}:=E_{2,1},\ v_{(23)}:=E_{1,2},\ v_{(12)}:=E_{2,3},\ v_{\mathrm{id}}=E_{1,3},$$
$$h_1:=E_{1,1}-E_{2,2},\ h_2:=E_{2,2}-E_{3,3},$$
where $E_{i,j}$ stands for the matrix whose only non-zero entry is a $1$ at the $i$-th row and $j$-th column and the indexes of $v$ are permutations. For $\sigma\in\mathrm{S}_3$, let $f_\sigma$ be the dual basis of $v_\sigma$. Then the Schubert varieties $X(\sigma)$ and the extremal functions $f_\sigma$ for $\sigma\in\mathrm{S}_3$ define a Seshadri stratification on the embedded projective variety $\mathrm{SL}_3/ B$.

The bonds can be determined using the Pieri-Chevalley formula \cite{Brion}. The Hasse diagram with bonds for this Seshadri stratification is depicted below:
$${\scriptsize
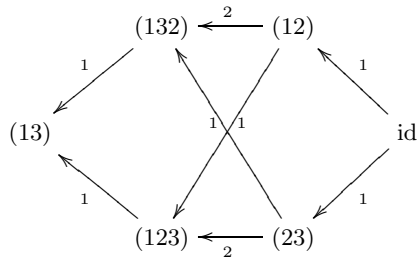
\xymatrix{
& (132)  \ar@{<-}[r]^2\ar@{<-}[rdd]^1&(12)\ar@{<-}[dr]^1&\\
(13)\ar@{<-}[ur]^1\ar@{<-}[dr]_1&&&\mathrm{id}\\
& (123) \ar@{<-}[r]_2\ar@{<-}[ruu]^1&(23)\ar@{<-}[ur]_1&\\
}
}
$$
\end{example}

For the example of a Seshadri stratification of a Schubert variety in a flag variety, see Section \ref{Schubert}.

\begin{rem}\label{extrabond}
Later in the article, we will mainly consider the affine cone $\hat{X}$ of $X$, so it is helpful to extend the stratification one step further. 

For a minimal element $p_0\in A$, the affine cone $\hat X_{p_{0}}\simeq \mathbb A^1$ is an affine line. Let $\hat{X}_{p_{-1}}$ denote the origin of $V$: it is contained in the affine cone of $X_p$ for any minimal element $p\in A$. The bond $b_{p_0,p_{-1}}$ is defined to be the vanishing multiplicity of $f_{p_0}$ at $\hat{X}_{p_{-1}}$, which coincides with the degree of $f_{p_0}$.

The set $\gls{Ahat}:=A\cup\{p_{-1}\}$ admits a poset structure by requiring $p_{-1}$ to be the unique minimal element. This poset structure is compatible with the containment relations between the affine cones $\hat{X}_p$, $p\in \hat{A}$. Similarly we have the Hasse diagram $\mathcal{G}_{\hat{A}}$, the bonds on it are described as above.
\end{rem}

\subsection{Generic hyperplane stratifications}\label{Sec:Generic}

According to the following proposition, the requirements in a Seshadri stratification are not as restrictive as it looks.

\begin{proposition}\label{Prop:Generic}
Every embedded projective variety $X \subseteq \mathbb{P}(V)$, smooth in codimension one, admits a Seshadri stratification.
\end{proposition}

\begin{proof}
Let $r$ be the dimension of $X$. For $r\geq 2$, applying the Bertini theorem \cite[Th\'eor\`eme 6.3]{Jou} to the smooth locus of $X_r:=X$, there (generically) exists $f_r\in V^*$ such that $X_{r-1}:=X\cap\mathcal{H}_{f_r}$ is an irreducible variety which is again smooth in codimension one and has the same degree as $X$. Repeating this construction gives irreducible varieties $X_{r-2},\ldots,X_{1}$ and $f_{r-1},\ldots,f_2\in V^*$ such that $X_k=X_{k+1}\cap\mathcal{H}_{f_{k+1}}$ is smooth in codimension one and the degree of $X_k$ is the same as the degree of $X$ for $k=1,\ldots,r-2$. Now $X_1$ is a smooth curve. A generic hyperplane $\mathcal{H}_{f_1}$ for $f_1\in V^*$ intersects $X_1$ at finitely many points $X_{0,1}=\{\varpi_1\}, \ldots, X_{0,s}:=\{\varpi_s\}$. It suffices to choose homogeneous functions $f_{0,k}$, $1\leq k\leq s$ on $V$ satisfying: $f_{0,k}$ vanishes on $\varpi_\ell$ for $\ell\neq k$, but $f_{0,k}$ is non-zero at $\varpi_k$.
\end{proof}

Since the hyperplanes are chosen generically (i.e. in an open set), the geometric definition of the degree of an embedded variety implies $\deg X=s$, the number of points we get in the last intersection. 

A Seshadri stratification arising from generic hyperplanes in this way will be termed a \emph{generic hyperplane stratification}. Let $A=
\{q_r,\ldots,q_1,q_{0,1},\ldots,q_{0,s}\}$ be the indexing poset with $X_{q_k}:=X_k$ and $X_{q_{0,\ell}}:=\{\varpi_\ell\}$.

\begin{example}\label{Ex:Bertini1}
The functions $f_{0,k}$ in the above proposition can be chosen in the following precise way. Let $h_i\in V^*$ be such that $h_i(\varpi_j)\neq 0$ for $1\leq j\leq s$ with $j\neq i$ and $h_i(\varpi_i)=0$. We define $g_{k}=\prod_{i\neq k}h_i$: it is homogeneous of degree $s-1$. These functions satisfy the requirement on $f_{0,k}$ in the above proof. With this choice, the extended Hasse diagram $\mathcal{G}_{\hat{A}}$ with bonds is depicted below:
{\tiny
$$
\xymatrix{
 & & & & q_{0,1}\ar[dl]_1 & \\
 q_r & q_{r-1} \ar[l]_1 & \cdots \ar[l]_1 & q_1 \ar@{}[r]|-{\raisebox{+2.2ex}{\text{\hskip -15pt}\vdots}} \ar[l]_1 & \vdots & q_{-1}. \ar[ul]_{s-1} \ar[dl]^{s-1} \ar@{}[l]|-{\raisebox{+2.2ex}{\text{\hskip 15pt}\vdots}}\\
 & & & & q_{0,s}\ar[ul]^1 & 
}
$$}
\end{example}

\begin{rem}\label{Rmk:Hibi}
In \cite{Hibi1} Hibi proved that every finitely generated positively graded ring admits a Hodge algebra structure. The construction in Proposition \ref{Prop:Generic} resembles a geometric version of the algebraic approach in \emph{loc.cit.} under the smooth in codimension one assumption. In our setup, this extra condition allows us to extract a convex geometric skeleton of $X$ (see Corollary \ref{Cor:SemiToric} and Section \ref{Sec:DegFormulaGHS}), and to deduce geometric properties of $X$ (see Corollary \ref{Cor:Degree2}).

In \cite{KMM} a similar Bertini-type argument is applied to construct flat degenerations of embedded projective varieties into complexity-one $T$-varieties.
\end{rem}

\section{Generalities on valuations}\label{ValuationsAndquasi-valuations}

From now on until Section \ref{section_standard_monomial_theory}, we fix a Seshadri stratification of $X\subseteq\mathbb{P}(V)$ with subvarieties $X_p$ and extremal functions $f_p$ for $p\in A$. 

\subsection{Definition and example}\label{ValutaionDefn}

We recall the definition and some basic properties of valuations and quasi-valuations.

\begin{definition}\label{quasidef}
Let $\mathcal R$ be a $\mathbb K$-algebra. A {\it quasi-valuation} on $\mathcal R$ with 
values in a totally ordered abelian group $\textsc{G}$ 
is a map $\nu: \mathcal R\setminus \{0\}\rightarrow \textsc{G}$ satisfying the following conditions:
\begin{itemize}
\item[{(a)}] $\nu(x + y) \ge \min\{\nu(x), \nu(y)\}$ for all $x, y\in \mathcal R \setminus \{0\}$ with $x+y\not=0$;
\item[{(b)}] $\nu(\lambda x) = \nu(x)$ for all $x\in \mathcal R  \setminus \{0\}$ and $\lambda \in \mathbb K^*$;
\item[{(c)}] $\nu(xy) \ge \nu(x) + \nu(y)$ for all $x,y \in \mathcal R \setminus \{0\}$ with $xy\not=0$.
\end{itemize}
The map $\nu$ is called a \emph{valuation} if  
the inequality in (c) can be replaced by an equality:
\begin{itemize}
\item[{(c')}] $\nu(xy)= \nu(x) + \nu(y)$ for all $x,y \in \mathcal R \setminus \{0\}$ with $xy\not=0$.
\end{itemize}
\end{definition}

\begin{rem}
Quasi-valuations on $\mathcal{R}$ can be thought of as synonyms of algebra filtrations on $\mathcal{R}$ (see Section 2.4 in \cite{KM}).
\end{rem}

The following properties are well-known for valuations (see \cite[Lemma~2.1.1]{MS}). Since the proof uses only axioms (a) and (b) of a valuation, it holds also for quasi-valuations.

\begin{lemma}\label{simpleproperties}
Let $\nu:\mathcal{R}\setminus\{0\}\to\textsc{G}$ be a quasi-valuation and $x, y\in \mathcal R \setminus \{0\}$.
\begin{itemize}
\item[{i)}] If $\nu(x)\not= \nu(y)$, then $\nu(x + y) = \min\{\nu(x), \nu(y)\}$.
\item[{ii)}] If $x+y\not=0$ and $\nu(x + y)>\nu(x)$, then  $\nu(x)=\nu(y)$.
\end{itemize}
\end{lemma}

\begin{lemma}[{\cite[Proposition 4.1]{FL}}]\label{minquasi}
Let $\nu_1,\ldots,\nu_k:\mathcal{R}\setminus\{0\}\to \textsc{G}$ be a family of quasi-valuations. The map $\nu:\mathcal{R}\setminus\{0\}\to\textsc{G}$ defined by
$$
g\mapsto\min\{\nu_j(g)\mid j=1,\ldots,k\}
$$
is a quasi-valuation on $\mathcal{R}$.
\end{lemma}

In algebraic geometry, valuations usually arise from vanishing orders of rational functions.

We come back to the setup in Section~\ref{setup}. Let $\gls{Rp}:=\mathbb K[X_p]$ be the homogeneous coordinate 
ring of $X_p$ with respect to the embedding $X_p\subseteq X\subseteq \mathbb P(V)$. In the following we consider 
$R_p$ often as the coordinate ring of the affine cone $\hat X_p\subseteq V$ over $X_p$.

If $p$ covers $q$ in $\hat{A}$, then $\hat X_q\subseteq \hat X_p$ is a prime divisor in $ \hat X_p$.
The local ring $\mathcal O_{\hat X_p,\hat X_q}$ is a discrete valuation ring because
$ \hat X_p$ is smooth in codimension one by (S1). 
Let $\gls{vpq}$ be the associated valuation. We refer to the value $\nu_{p,q}(f)$ for $f\in R_p\setminus\{0\}$ as 
the \textit{vanishing multiplicity} of $f$ in the divisor $\hat X_q$:
\begin{equation}\label{geometricvaluation}
\nu_{p,q}: R_p\setminus\{0\}\rightarrow \mathbb Z.
\end{equation} 

The valuation $\nu_{p,q}$ can be naturally extended to a valuation on 
$\mathbb K(\hat X_p)={\rm Quot\,} R_p$, 
the quotient field of $R_p$, by the rule:
$$\nu_{p,q}\left(\frac{g}{h}\right):=\nu_{p,q}(g)-\nu_{p,q}(h)\text{ for }g,h\in R_p\setminus\{0\}.$$

\begin{rem}\label{specialvalues}
As a continuation of Remark~\ref{extrabond}, for a minimal element $p_0\in A$, $R_{p_0}$ is a polynomial ring. We define $\nu_{p_0,p_{-1}}$ to be the vanishing multiplicity of a polynomial in $R_{p_0}$ in $\{0\}$.
\end{rem}

\subsection{Valuation under normalization}\label{normalization}
The following Lemma compares the prime divisors in the normalization $\gls{Xptilde}$ of $\hat X_p$ with those in $\hat X_p$. Note that $X_p$ is smooth in codimension one.

\begin{lemma}\label{BijectionDivisor}
The normalization map $\omega: \tilde X_p\rightarrow \hat X_p$ induces a bijection $\omega_D:\tilde Z\mapsto \omega(\tilde Z)$ between the set of prime 
divisors $ \tilde Z\subseteq \tilde X_p$ and the set of prime divisors 
$Z\subseteq \hat X_p$. In addition, for all non-zero $f\in\mathbb K(\hat X_p)=\mathbb K(\tilde X_p)$ and all prime divisors $\tilde Z\subseteq \tilde X_p$ holds:
$\nu_{\tilde Z}(f)=\nu_{\omega(\tilde Z)}(f)$.
\end{lemma}
\begin{proof}
Let $\hat U\subseteq \hat X_p$
be the open and dense subset of smooth points and denote by $\tilde U\subseteq \tilde X_p$ the 
open and dense subset obtained as preimage $\omega^{-1}(\hat U)$.

By axiom (S1), $\hat X_p\setminus \hat U$ is of codimension
greater or equal to $2$. Since the normalization map is finite, the same holds for 
$\tilde X_p\setminus \tilde U$. It follows that any prime divisor as well as the induced valuation in 
$\hat X_p$ respectively $\tilde X_p$ is completely determined by its intersection with $\hat U$ respectively
$\tilde U$. The claim follows now since the normalization map $\omega:\tilde X_p \rightarrow \hat X_p$ induces an isomorphism  
$\omega\vert_{\tilde U} :\tilde U \rightarrow \hat U$.
\end{proof}
\subsection{Rees quasi-valuations}
For an element $p\in A$ let $I_p=(f_p\vert_{\hat X_p})$ be the principal ideal in $R_p$ generated by the restriction of the extremal function $f_p$ (see Definition \ref{Defn:SS}) to $\hat X_p$.
By abuse of notation we write in the following often just $f_p$ instead of $f_p\vert_{\hat X_p}$. 
We define a map $\nu_{I_p}:R_p\setminus\{0\}\rightarrow\mathbb Z_{\ge 0}\cup\{\infty\}$ by:
\begin{equation}\label{poweridealquasi-valuation}
\text{for any}\ g\in R_p\setminus\{0\},\ \nu_{I_p}(g):=\max\left\{m\in \mathbb Z_{\ge 0}\cup\{\infty\}\left\vert\, g \in I_p^m\right.\right\}.
\end{equation}
In our situation the value $\infty$ is never attained.

\begin{lemma}[{\cite[Section 2.2]{R1}}]
$\nu_{I_p}$ defines a quasi-valuation $\nu_{I_p}:R_p\setminus\{0\}\rightarrow \mathbb Z$.
\end{lemma}

Such a quasi-valuation is called \emph{homogeneous}, if for any $g\in R_p\setminus\{0\}$ and $m\in\mathbb{N}$, $\nu_{I_p}(g^m)=m\nu_{I_p}(g)$. The quasi-valuation $\nu_{I_p}$ is not necessarily homogeneous. Samuel \cite{Sam} introduced a limit procedure to homogenize such a quasi-valuation.

We denote by $\bar\nu_{I_p}$ the corresponding {\it homogenized quasi-valuation} (\cite[Lemma 2.11]{R1}): by definition, for $g\in R_p\backslash\{0\}$, we set
\begin{equation}\label{homogenizedpoweridealquasi-valuation}
\bar\nu_{I_p}(g):=\lim_{k\rightarrow \infty}\frac{\nu_{I_p}(g^k)}{k} \in \mathbb Q_{\ge 0}\cup\{\infty\}.
\end{equation}

The limit exists and is a non-negative rational number (see \cite{N,R2}). The homogenized quasi-valuation has an interpretation as the minimum function
over some rescaled discrete valuations: 

\begin{theorem}[{\cite[Valuation Theorem]{R2}}]\label{ValuationTheoremRees}
There exists a finite set of discrete, integer-valued valuations $\eta_1,\ldots,\eta_k$ on $R_p$ and  integers $e_1,\ldots, e_k$ such that for all $g\in R_p\setminus\{0\}$:
$$
\bar{\nu}_{I_p}(g):=\min\left\{\frac{\eta_j(g)}{e_j}\mid j=1,\ldots,k\right\}.
$$ 
\end{theorem}

In our setting, the valuations $\eta_j$ get an interpretation as vanishing multiplicities at the prime divisors occurring in $X_p\cap \mathcal{H}_{f_p}$. Let $q_1,\ldots,q_k\in A$ be such that 
$$X_p\cap \mathcal{H}_{f_p}=\bigcup_{i=1,\ldots,k} X_{q_i}.$$ 
The affine cones $\hat X_{q_i}$ are prime divisors in $\hat X_p$ inducing valuations 
$\nu_{p,q_i}:R_p\setminus\{0\}\rightarrow \mathbb Z$ (see~\eqref{geometricvaluation}). Let $b_{p,q_i}$, $i=1,\ldots,k$, be the corresponding bonds, which are the vanishing multiplicities of $f_p$ at the prime divisors $\hat{X}_{q_i}$ (see Section~\ref{Hasse}). The Valuation Theorem gets the following interpretation:

\begin{proposition} \label{minimuvaluation}
Let $g\in R_p\setminus\{0\}$. We have
\begin{equation}\label{valuationtheorem}
\bar\nu_{I_p}(g) =\min\left\{\left.\frac{\nu_{p,q_i}(g)}{b_{p,q_i}}\,\right\vert\, i=1,\ldots,k\right\}.
\end{equation}
\end{proposition}

\begin{proof}
Let $\tilde R_p$ be the normalization of $R_p$ (in $\mathrm{Quot\,}R_p$), $\tilde I_p$ be the principal ideal in
$\tilde R_p$ generated by $f_p$ and set $\tilde{Y}_p:=\mathrm{Spec\,} \tilde R_p$, the normalization of $\hat X_p$. 
We can construct the quasi-valuation $\nu_{\tilde I_p}:\tilde R_p\setminus\{0\}\rightarrow\mathbb Z$
and its homogenized version $\bar \nu_{\tilde I_p}$ in the same way as in \eqref{poweridealquasi-valuation}
and \eqref{homogenizedpoweridealquasi-valuation}. The constructions are related by Lemma 2.2 in \cite{R2}: 
one has for $g\in R_p\setminus\{0\}$, $\bar \nu_{\tilde I_p}(g)=\bar \nu_{I_p}(g)$.

Let $\Sigma$ be the set of all valuations $\nu_Z:{\rm Quot\,}\tilde R_p\setminus\{0\}\rightarrow \mathbb Z$ 
induced by prime divisors $Z\subseteq \tilde{Y}_p$. Since $\tilde R_p$ is an integrally closed noetherian ring, this set of valuations
makes $\tilde R_p$ into a finite discrete principal order (a.k.a. Krull domain, see \cite{B} Chapter VII, \S 1.3), so 
$\tilde R_p$ satisfies the assumptions for Lemma 2.3 in \cite{R2}.
The proof of the lemma identifies the valuations in Theorem~\ref{ValuationTheoremRees} as the 
valuations $\nu_Z$ such that $\nu_Z(f_p\vert_{\hat X_p})>0$. Keeping in mind the bijection in Lemma~\ref{BijectionDivisor} and (S3) in Definition \ref{Defn:SS}, 
we see that in such a case $\omega(Z)=\hat X_{q_i}$ is an irreducible component of the vanishing set of $f_p\vert_{\hat X_p}$, and the 
scaling factor is $b_{p,q_i}^{-1}=(\nu_{p,q_i}(f_p))^{-1}$ (see the proof of Lemma 2.3 in \cite{R2}).
\end{proof}

\begin{rem}\label{remarkvaluationtheoremormalized}\rm
The appropriate reformulation of the proposition above holds for elements in the normalization too. 
The proof above together with Lemma~\ref{BijectionDivisor} shows for $g\in \tilde R_p\setminus\{0\}$:
\begin{equation}\label{valuationtheoremormalized}
\bar\nu_{\tilde I_p}(g) =\min\left\{\left.\frac{\nu_{p,q_i}(g)}{b_{p,q_i}}\,\right\vert\, i=1,\ldots,k\right\}.
\end{equation}
\end{rem}

\section{Codimension one}\label{Section:ValuationCodimOne}

In this section we modify the procedure in \cite{KK, LM} in codimension one to overcome the difficulty mentioned in the introduction. In the next two sections we will link the codimension one constructions to obtain a higher rank valuation.

Let $p,q\in A$ be such that $p$ covers $q$ and $g\in\mathbb{K}(\hat{X}_p)\setminus\{0\}$ be a rational function. Let $b_{p,q}=\nu_{p,q}(f_p)$ be the bond between $p$ and $q$ (see Section \ref{Hasse}).
Let $\gls{N}$ be the least common multiple of all bonds in $\mathcal G_{A}$, so the number $N\frac{\nu_{p,q}(g)}{b_{p,q}}\in\mathbb Z$.
We set
\begin{equation}\label{functionprocedure}
h:=\frac{g^N}{f_p^{N\frac{\nu_{p,q}(g)}{b_{p,q}}}}\in \mathbb K(\hat X_p).
\end{equation} 

\begin{lemma}\label{one}
The restriction $h\vert_{\hat X_q}$ is a well-defined, non-zero rational function on $\hat X_{q}$.
\end{lemma}

\begin{proof}
Note that
$$
\nu_{p,q}(h)=\nu_{p,q}\left(\frac{g^N}{f_p^{N\frac{\nu_{p,q}(g)}{b_{p,q}}}}\right) 
=N\nu_{p,q}(g)-N{\frac{\nu_{p,q}(g)}{b_{p,q}}}b_{p,q}=0,
$$
so this rational function is an element of the local ring $\mathcal O_{\hat X_p,\hat X_{q}}$ 
of the prime divisor $\hat X_{q}\subseteq \hat X_p$.
But it is not in its maximal ideal $\mathfrak m_{\hat X_q,\hat X_{p}}$ and
hence its restriction gives a non-zero element in the residue field 
$\mathcal O_{\hat X_q,\hat X_{p}}/\mathfrak m_{\hat X_q,\hat X_{p}}$, which is the field $\mathbb K(\hat X_{q})$ of rational functions on $\hat X_{q}$.  
\end{proof}

\begin{rem}\label{alternativH}
Instead of taking the function $h$ as above one could take $\tilde h:= g^{b_{p,q}}/f_p^{\nu_{p,q}(g)}$. 
Lemma~\ref{one} holds for $\tilde h$ with the same proof.  In fact, $h=\tilde h^{\frac{N}{b_{p,q}}}$. 
Since it is later more convenient to work uniformly with the $N$-th power instead of the $b_{p,q}$-th power,
 we will stick to this construction. For the valuation which will be defined in 
Section~\ref{chainvaluation} the choice makes no difference, one just has to rescale the values appropriately, see
Remark~\ref{alternativeEvaluation}.
\end{rem}

What can we say about $h$ if the starting function $g\in R_p$ is a regular function? 

For the inductive procedure we will use later it is necessary to take a slightly more general point of view. 
Consider again the setting in Proposition~\ref{minimuvaluation} 
respectively Remark~\ref{remarkvaluationtheoremormalized} and
let $q_1,\ldots,q_k\in A$ be such that $X_p\cap \mathcal{H}_{f_p}=\bigcup_{i=1,\ldots,k} X_{q_i}$.
Let $g\in \mathbb K(\hat X_p)\setminus\{0\}$ be a rational function which  
is integral over $\mathbb K[\hat X_p]$. By Lemma~\ref{BijectionDivisor}, this property is equivalent to 
$\nu_Z(g)\ge 0$ for all prime divisors $Z\subseteq \hat X_p$. 

\begin{proposition}\label{normularinherited}
Let $g\in \mathbb K(\hat X_p)\setminus\{0\}$ be a rational function which  
is integral over $\mathbb K[\hat X_p]$. We assume that the enumeration of the divisors $\hat X_{q_1},\ldots,\hat X_{q_k}$ is such that 
\begin{equation}\label{choiceOfq1}
\forall\, i=1,\ldots k,\ \bar\nu_{\tilde I_p}(g)=\frac{\nu_{p,q_1}(g)}{b_{p,q_1}}\le \frac{\nu_{p,q_i}(g)}{b_{p,q_i}}.
\end{equation}
Set 
$h=g^Nf_{p}^{-N\frac{\nu_{p,q_1}(g)}{b_{p,q_1}}}$ (as in \eqref{functionprocedure}). 
Then $h$ is integral over $\mathbb K[\hat X_p]$, and 
$h\vert_{\hat X_{q_1}}\in\mathbb K(\hat X_{q_1})$ is integral over $\mathbb K[\hat X_{q_1}]$.
\end{proposition}

\begin{proof}
Given a prime divisor $Z\subseteq \hat X_p$, we have $$\nu_Z(h)=N\left(\nu_Z(g)-\nu_Z(f_p)\frac{\nu_{p,q_1}(h)}{b_{p,q_1}}\right).$$
By assumption we have $\nu_Z(g)\ge 0$ for all prime divisors $Z\subseteq \hat X_p$ and 
$\nu_Z(f_p)=0$ for $Z\not= \hat X_{q_i}$, $j=1,\ldots,k$. It follows: if $Z\not= \hat X_{q_j}$, then $\nu_Z(h)\ge 0$.

For the prime divisors $\hat X_{q_j}$, $j=1,\ldots,k$, and the associated valuations $\nu_{p,q_j}$ we obtain:
$$
\begin{array}{rcl}
\nu_{p,q_j}(h)=\nu_{p,q_j}\left(\frac{g^N}{f_p^{N(\frac{\nu_{p,q_1}(g)}{b_{p,q_1}})}}\right) 
&=&N(\nu_{p,q_j}(g)-\frac{\nu_{p,q_1}(g)}{b_{p,q_1}}\nu_{p,q_j}(f_p))\\
&=& N\nu_{p,q_j}(f_p)\left( {\frac{\nu_{p,q_j}(g)}{b_{p,q_j}}}-\frac{\nu_{p,q_1}(g)}{b_{p,q_1}}\right) \ge 0
\end{array}
$$
by the choice of $q_1$ (see \eqref{choiceOfq1}). Hence
$\nu_Z(h)\ge 0$ for all prime divisors $Z\subseteq \hat X_p$, so $h$ is integral in $\mathbb K(\hat X_p)$
over $\mathbb K[\hat X_p]$. 
By Lemma~\ref{one}, $h\vert_{\hat X_{q_1}}$ is a well-defined, non-zero rational function, and hence $h\vert_{\hat X_{q_1}}\in\mathbb{K}(\hat{X}_{q_1})$ is also integral over $K[\hat X_{q_1}]$.
\end{proof}

\section{Maximal chains and sequences of rational functions}\label{chainsandfunctions}

We fix a maximal chain in $A$ joining $p_{\max}$ with a minimal element $p_0$:
\begin{equation}\label{fixachain}
\gls{C}: p_{\max}=p_r>\ldots >p_2>p_1>p_0.
\end{equation}

In particular, $r$ is the length of every maximal chain. To avoid double indexes as much as possible,
we use abbreviations.  Since we have fixed a maximal chain, an element $p\in\mathfrak C$ 
is either the minimal element, or there is a unique element in $\mathfrak C$, say $q$, covered by $p$. It makes sense to omit the second index and write $\nu_p$ and $b_p$ instead of $\nu_{p,q}$ and $b_{p,q}$. Moreover, when $p=p_i$ we will simplify the notation further by just writing
\begin{equation}\label{abbreviations}
\text{$\nu_i$ (resp. $b_i$, $f_i$, $X_i$) instead of $\nu_{p_i,p_{i-1}}$ (resp. $b_{p_i,p_{i-1}}$, $f_{p_i}$, $X_{p_i}$).}
\end{equation}

We use the procedure in \eqref{functionprocedure} to attach to a non-zero function $g\in R$ 
inductively a sequence $(g_r,\ldots,g_0)$ of non-zero rational functions $g_j\in\mathbb K(\hat X_{j})$, $j=0,1,\ldots,r$. 
For each $j=0,1,\ldots,r$, the function $g_j$ will depend on $g$ and $f_{{j+1}},\ldots,f_{r}$.
Starting with a regular function $g\in R\setminus\{0\}$, the following inductive procedure is well defined by Lemma~\ref{one}: we set
$$
g_r:=g\quad\textrm{and}\quad D_r=\frac{\nu_{r}(g_r)}{b_r},
$$
and then inductively for $j=r-1,\ldots,1,0$:
\begin{equation}\label{definitiocproc1}
g_j=\left.\frac{g_{j+1}^N}{f_{{j+1}}^{ND_{j+1}}}\right\vert_{\hat X_{j}} 
\quad\textrm{and}\quad 
D_j=\frac{\nu_{j}(g_j)}{b_{j}}.
\end{equation}

\begin{rem}\label{nearlyclosed}
We can provide a {\it nearly} closed formula for $g_i$: for $r\ge j\ge i+1$ let $D_j$ be as defined above in \eqref{definitiocproc1}, then 
\begin{equation}\label{closedfromula1}
g_i=g^{N^{r-i}}f_{r}^{-N^{r-i}D_r}
f_{{r-1}}^{-N^{r-i-1}D_{r-1}}
\cdots 
f_{{i+1}}^{-ND_{i+1}}\vert_{\hat X_{i}}.
\end{equation}
It is only {nearly} closed because the valuations $\nu_{j}(g_{j})$, $j>i$, show up in the formula
for the numbers $D_j$. But for our purpose this will be good enough.
\end{rem}

We give this sequence of functions a name:

\begin{definition}\label{sequence-of-rational-functions}
The tuple $\gls{gc}:=(g_r,\ldots,g_1,g_0)$ associated to $g\in R\setminus\{0\}$
is called the \textit{sequence of rational functions}
associated to $g$ along $\mathfrak C$.
\end{definition}

Before going further to define a higher rank valuation from this sequence of rational functions, we look at some concrete examples.

\begin{example}\label{constant_function}
The constant function $a$, $a\in \mathbb K^*$, vanishes nowhere, so $\nu_{j}(a\vert_{\hat X_{j}})=0$ for all $j=1,\ldots,r$, and the sequence associated to the constant function is $a_{\mathfrak C}=(a,a^N,\ldots,a^{N^r})$. Note that this holds independent of the choice of the maximal chain $\mathfrak C$.
\end{example}

\begin{example}\label{fpitupel}
We consider the case when $g$ is the extremal function $f_{i}$ for some $0\le i\le r$.
By Lemma~\ref{Lem:SS}, $f_{i}$ does not vanish identically on $\hat X_{j}$ for $j\ge i$, one determines inductively: 
$$
g_r=f_{i}, D_r=0,\textrm{\ and\ }g_{r-1}=f^{N}_{i}, D_{r-1}=0,\textrm{\ and\ }\ldots,\textrm{\ and\ }
g_{i+1}=f^{N^{r-i-1}}_{i}, D_{i+1}=0.
$$
Next consider the function $g_i=f_{i}^{N^{r-i}}$. This function vanishes on the divisor $\hat X_{{i-1}}$ and we have
$D_i=\frac{\nu_{i}(f_{i}^{N^{r-i}})} {b_i}=N^{r-i}$ by the definition of $b_{{i}}$.
It follows
$$
g_{i-1}=\frac{g_i^N}{(f_{i})^{ND_i}}=\frac{f_{i}^{N^{r-i+1}}}{f_{i}^{N^{r-i+1}}} =1\textrm{\ and\ }D_{i-1}=0.
$$
The procedure implies now $g_j=1$ for all $j<i$. Summarizing the above computation gives
$$
(f_{i})_{\mathfrak C}=(f_{i},f_{i}^N,f_{i}^{N^2},\ldots,f_{i}^{N^{r-i}},1,\ldots,1).
$$
This holds for all choices of a maximal chain $\mathfrak C$ as long as $p_i\in\mathfrak C$.
If $p\not\in \mathfrak C$, then $(f_{p})_{\mathfrak C}$ looks rather different, as the following example shows.
\end{example}

\begin{example}\label{grasstwofour2}
Let $X=\mathrm{Gr}_2\mathbb K^4$ with the Seshadri stratification defined in Example~\ref{grasstwofour1} and~\ref{grasstwofour}. 
The bonds are all equal to $1$ and hence $N=1$. We fix as maximal chain 
$$
\mathfrak C: (3,4)>(2,4)>(2,3)>(1,3)>(1,2).
$$
For the Pl\"ucker coordinates we get as sequences of rational functions:
$$
\begin{array}{cclccll}
(x_{3,4})_{\mathfrak C}&=&(x_{3,4},1,1,1,1),&&(x_{2,4})_{\mathfrak C}&=&(x_{2,4},x_{2,4},1,1,1),\\
(x_{2,3})_{\mathfrak C}&=&(x_{2,3},x_{2,3},x_{2,3},1,1),&&(x_{1,4})_{\mathfrak C}&=&(x_{1,4},x_{1,4},\frac{x_{1,4}}{x_{2,4}},\frac{x_{2,3}x_{1,4}}{x_{2,4}},1),\\
(x_{1,3})_{\mathfrak C}&=&(x_{1,3},x_{1,3},x_{1,3},x_{1,3},1),&&(x_{1,2})_{\mathfrak C}&=&(x_{1,2},x_{1,2},x_{1,2},x_{1,2},x_{1,2}).\\
\end{array}
$$
The Pl\"ucker coordinate $x_{1,4}$ is the only extremal function whose index is not in $\mathfrak C$. By Example~\ref{fpitupel},
the sequence $(x_{1,4})_{\mathfrak C}=(g_4,\ldots,g_0)$ is the only one which needs an explanation. 

Recall that for the sequence $(g_4,\ldots,g_0)$ of rational functions associated to $x_{1,4}$ along $\mathfrak{C}$, we denote $D_j=\nu_j(g_j)/b_j$ for $j=0,1,\ldots,4$. In this example all bonds are equal to $1$, hence $D_j=\nu_j(g_j)$.

The restriction of $x_{1,4}$ to $\hat X(2,4)$ does not vanish identically and hence $D_{4}=0$; this gives $g_4=g_3=x_{1,4}$. 

In order to compute the vanishing order, we introduce coordinates on an open subset of $\hat{X}(2,4)$. For this we choose a maximal torus and a Borel subgroup for $\mathrm{SL}_4(\mathbb{K})$ as in Example \ref{fullflag}. We use the standard enumeration of the simple roots, i.e. $\alpha_1=\epsilon_1-\epsilon_2$,
$\alpha_2=\epsilon_2-\epsilon_3$ and  $\alpha_3=\epsilon_3-\epsilon_4$.

For a root $\alpha$, let $U_{\alpha}\subseteq \mathrm{SL}_4(\mathbb{K})$ be the corresponding root subgroup. The root subgroup associated to the a (negative) simple root is given by 
$$U_{-\alpha_i}(t)=1\hskip -3pt\textrm{\rm I}+tE_{i+1,i},\ \ i=1,2,3.$$
We consider in $\hat X(2,4)$ the open subset $U_{-\alpha_3}(d)U_{-\alpha_1}(c)U_{-\alpha_2}(b)(a e_1\wedge e_2)$,
it is equal to
\begin{equation}\label{opensetGrass}
\left\{ae_1\wedge e_2 +abe_1\wedge e_3+abc e_2\wedge e_3+abd e_1\wedge e_4+abcd e_2\wedge e_4\mid a,b,c,d\in\mathbb K\right\}.
\end{equation} 
This open subset is compatible with the other Schubert varieties in the maximal chain $\mathfrak C$ contained in $\hat{X}(2,4)$, i.e. we get open subsets in the affine cones over these Schubert varieties by setting some of the coordinates equal to $0$. For instance, by setting $d=0$ we obtain an open subset in the affine cone $\hat{X}(2,3)$.

The restricted Pl\"ucker coordinate $g_3=x_{1,4}\vert_{\hat X(2,4)}$ 
vanishes on the divisor $\hat X(2,3)$ with multiplicity  $1$, so $D_3=1$. To get $g_2$ in the sequence $(x_{1,4})_{\mathfrak C}$ 
we have to divide $g_3$ by $x_{2,4}\vert_{\hat X(2,4)}$ and get $g_2=\frac{x_{1,4}}{x_{2,4}}\vert_{\hat X(2,3)}$.  The rational function $\frac{x_{1,4}}{x_{2,4}}$ takes on the open set in 
\eqref{opensetGrass} the value $\frac{1}{c}$, so $\frac{x_{1,4}}{x_{2,4}}\vert_{\hat X(2,3)}$ has a pole of order $1$ at the divisor $\hat X(1,3)\subseteq\hat X(2,3)$,
which implies $D_2=-1$.
Hence, for the term $g_1$ in the sequence $(x_{1,4})_{\mathfrak C}$, one has to multiply $g_2$ by $x_{2,3}$ 
and obtains $g_1=\frac{x_{1,4}x_{2,3}}{x_{2,4}}\vert_{\hat X(1,3)}$. 
This rational function takes on the open set in \eqref{opensetGrass} the same values as the Pl\"ucker coordinate $x_{1,3}$, and hence vanishes with multiplicity $1$ 
on the divisor $\hat X(1,2)\subseteq\hat X(1,3)$, so we have $D_1=1$. We have hence to divide $g_1$ by $x_{1,3}$ to get $g_0$, 
which is the constant function $1$.

The calculations may be simplified by using the Pl\"ucker relation $x_{1,2}x_{3,4}+x_{2,3}x_{1,4}-x_{1,3}x_{2,4}=0$.
\end{example}

Even if one starts with a regular function $g\in R\setminus\{0\}$, for a fixed maximal chain $\mathfrak C$
there is no reason why a function $g_j$ occurring in the sequence $g_{\mathfrak C}=(g_r,\ldots,g_0)$
should be a regular function on the respective subvariety $\hat X_{p_j}$, see Example~\ref{grasstwofour2} with $g=x_{1,4}$.
Nevertheless, Proposition~\ref{normularinherited} shows that for a fixed function, poles can be avoided if one chooses the maximal chain carefully: 

\begin{coro}\label{normularsequence}
For every regular function $g\in R\setminus\{0\}$ there exists a maximal chain $\mathfrak C$ so that 
the associated tuple $g_{\mathfrak C}=(g_r,\ldots,g_0)$ consists of rational functions $g_j\in\mathbb K(\hat{X}_{p_j})$ such that $\nu_{j}(g_j)\ge 0$ for all $j=1,\ldots,r$.
\end{coro}

\section{Valuations from maximal chains}\label{value_valuations}

Throughout this section, we will work with a fixed maximal chain in $A$: 
$$\mathfrak C:p_{\max}=p_r>p_{r-1}>\ldots >p_0.$$

As in the previous section, to avoid the usage of double indexes as much as possible, we keep the conventions in \eqref{abbreviations}.
Moreover, if not mentioned otherwise, we will write $e_i$ instead of $e_{p_i}$.

Let $\mathbb Q^{\mathfrak C}$ be the $\mathbb Q$-vector space with basis $\{e_{j}\mid j=0,\ldots,r\}$. We will write $v=(a_r,\ldots,a_0)$ for the vector $v=\sum_{j=0}^r a_je_{j}\in \mathbb Q^{\mathfrak C}$.

\begin{definition}\label{orderchain}
We endow $\mathbb Q^{\mathfrak C}$ with the \textit{lexicographic order}, i.e. 
$$(a_r,\ldots,a_0)\ge (b_r,\ldots,b_0)\  \text{if and only if}\ a_r>b_r,\ \text{or}\ a_r=b_r\ \text{and}\ a_{r-1}>b_{r-1},\ \text{or etc}.$$
\end{definition}

This total order is compatible with the addition of vectors: for any $u,v,w\in\mathbb{Q}^{\mathfrak{C}}$, if $u\geq v$, then $u+w\geq v+w$ holds.

\subsection{Linking up valuations}

We start with linking together the rank one valuations arising from the covering relations in the poset $A$.

\begin{definition}\label{chainvaluation}
For $g\in R\setminus\{0\}$ let $g_{\mathfrak C}=(g_r, \ldots, g_0)$ be the associated sequence of rational functions
(see Definition~\ref{sequence-of-rational-functions}).
We attach to the maximal chain $\mathfrak C$ the map
$$
\gls{Vc} :R\setminus\{0\}\rightarrow \mathbb Q^{\mathfrak C},\quad
g\mapsto \sum_{j=0}^r  \left(\frac{\nu_{j}(g_{j})}{N^{r-j} b_{j}}\right) e_{j}.
$$
\end{definition}

\begin{rem}\label{idealinequality2}
In terms of the numbers $D_j$, $j=0,\ldots,r$, defined in \eqref{definitiocproc1}, we have:
$$
\mathcal V_{\mathfrak C}(g)= \sum_{j=0}^r \frac{D_j}{N^{{r-j}}} e_{j}.
$$
\end{rem}

We will prove in Proposition \ref{VCvaluation} that $\mathcal V_{\mathfrak C}$ is a valuation. Before that we introduce a lattice containing the image of $\mathcal{V}_{\mathfrak{C}}$ and discuss some examples.

\begin{rem}\label{idealinequality}
The construction is similar to the one in Newton-Okounkov theory which associates
a valuation to a flag of subvarieties (see Introduction). The scaling factor $\frac{1}{N^{r-j}}$ in our construction shows up due to the 
fact that $g$ occurs with the power $N^{r-j}$ in $g_j$, see Remark~\ref{nearlyclosed}. The scaling factor 
$\frac{1}{b_{j}}$ occurs due to the fact that we divide by the function $f_{j}$ and not, as usual, 
by a (local) equation defining the prime divisor scheme theoretically.
\end{rem}

\begin{rem}\label{alternativeEvaluation}
For the construction of the functions $g_r,g_{r-1},\ldots$ we have made use of the procedure described in \eqref{functionprocedure}.
One could as well define these functions by using the algorithm described in Remark~\ref{alternativH}: $\tilde g_r=g$ and, 
inductively, $\tilde g_{i-1}=\tilde g_{i}^{b_i}/f_i^{\nu_i(\tilde g_i)}\vert_{\hat X_{i-1}}$ for $i=1,\ldots,r$. As map one would take instead:
$$
\gls{Vctilde}(g)= \sum_{j=0}^r  \left(\frac{\nu_{j}(\tilde g_{j})}{\prod_{k=j}^r b_{k}}\right) e_{j}.
$$
One sees easily from the definition that these two maps coincide: for all $g\in R\setminus\{0\}$ one has $\widetilde{\mathcal V}_{\mathfrak C}(g)={\mathcal V}_{\mathfrak C}(g)$.
\end{rem}

An immediate consequence of Remark~\ref{alternativeEvaluation} is:

\begin{lemma}\label{Lem:ImageValC}
The map ${\mathcal V}_{\mathfrak C}$ takes values in the lattice:
\begin{equation}\label{latticedef}
\gls{Lc}=\{\ell=(\ell_r,\ldots,\ell_0)\in \mathbb Q^{\mathfrak C}\,\vert\, b_r\cdots b_j\ell_j\in\mathbb{Z},\ 0\leq j\leq r\}.
\end{equation}
\end{lemma}

\begin{rem}\label{Rmk:Bond1}
 If all bonds in the extended Hasse diagram $\mathcal{G}_{\hat{A}}$ are equal to $1$, then $L^{\mathfrak C}\,\cong\,\mathbb Z^{r+1}$.
\end{rem}

\begin{example}\label{extremalvaluation}
Let $p_i$ be an element in the maximal chain $\mathfrak C$.
The renormalization is chosen so that $\mathcal V_{\mathfrak C}(f_{i})= e_{i}$.
 
Indeed, by Example~\ref{fpitupel} we know that 
$(f_{i})_{\mathfrak C}=(f_{i},f_{i}^N,f_{i}^{N^2},\ldots,f_{i}^{N^{r-i}},1,\ldots,1)$.
Let $\mathcal V_{\mathfrak C}(f_{i})=(a_r,\ldots,a_0)$.
Since $f_{i}\vert_{\hat X_{j}}\not\equiv 0$ for $j\ge i$, it follows that 
$a_j=0$ for $j>i$. Since $1$ is a nowhere vanishing function, one has $a_j=0$ for $j<i$.
It remains to determine $a_i$. Now 
$$
\nu_{i}(f_{i}^{N^{r-i}})=N^{r-i}\nu_{i}(f_{i})=N^{r-i}b_{i},
$$ 
so the renormalization implies $a_i=1$, and hence $\mathcal V_{\mathfrak C}(f_{i})= e_{i}$.
Note that this holds no matter which maximal chain one chooses as long as $p_i$ shows up in the chain.
The situation becomes different if $p\not\in\mathfrak C$, as the following example shows.
\end{example}

\begin{example}\label{grasstwofour3}
Let $X=\mathrm{Gr}_2\mathbb K^4$ be as in the Examples~\ref{grasstwofour1},~\ref{grasstwofour} and~\ref{grasstwofour2}. We fix the maximal chain $\mathfrak C: (3,4)>(2,4)>(2,3)>(1,3)>(1,2)$ and $N=1$ as in Example~\ref{grasstwofour}. The calculations in Example~\ref{grasstwofour2} and Example~\ref{extremalvaluation} imply that the values of $\mathcal V_{\mathfrak C}$ 
on the Pl\"ucker coordinates are given by:
$$
\begin{array}{clcllcll}
\mathcal V_{\mathfrak C}(x_{3,4})=(1,0,0,0,0),
&&\mathcal V_{\mathfrak C}(x_{2,4})=(0,1,0,0,0),
&&\mathcal V_{\mathfrak C}(x_{2,3})=(0,0,1,0,0),\\
\mathcal V_{\mathfrak C}(x_{1,4})=(0,1,-1,1,0), 
&&\mathcal V_{\mathfrak C}(x_{1,3})=(0,0,0,1,0),
&&\mathcal V_{\mathfrak C}(x_{1,2})=(0,0,0,0,1).\\
\end{array}
$$
\end{example}

\subsection{$\mathcal V_{\mathfrak C}$ is a valuation}

We extend the map $\mathcal{V}_{\mathfrak C}$ to $\mathbb{K}(\hat{X})=\mathrm{Quot} R$ by:
$$\mathcal{V}_{\mathfrak C}\left(f/g\right):=\mathcal{V}_{\mathfrak C}(f)-\mathcal{V}_{\mathfrak C}(g),\ \text{for}\ f,g\in R\setminus\{0\}.$$

The goal of this subsection is to show that the map is well defined and

\begin{proposition}\label{VCvaluation}
$\mathcal V_{\mathfrak C}$ is an $L^{\mathfrak C}$-valued valuation on $\mathbb K(\hat X)$.
\end{proposition}

\begin{definition}\label{valumonoid}
We denote by $\gls{VCX}(X)$ the \textit{valuation monoid} associated to $X$ by $\mathcal V_{\mathfrak C}$, i.e. 
$\mathbb V_{\mathfrak C}(X)=\{\mathcal V_{\mathfrak C}(g)\mid g\in R\setminus\{0\}\}\subseteq L^{\mathfrak C}$.
\end{definition} 

\begin{proof}
It suffices to consider elements in $R\setminus\{0\}$.
Given $g,h\in R\setminus\{0\}$, one verifies by induction (using Remark~\ref{nearlyclosed}) that if 
$g_{\mathfrak C}=(g_r,\ldots,g_0)$ and $h_{\mathfrak C}=(h_r,\ldots,h_0)$, then
$(gh)_{\mathfrak C}=(g_rh_r,\ldots,g_0h_0)$. Since the $\nu_{j}$, $j=0,\ldots,r$, 
are valuations, property (c')  in Definition~\ref{quasidef} follows.

If one replaces $g\not=0$ by a non-zero scalar 
multiple $\lambda g$, then the components $g_j$ in the tuple 
$g_{\mathfrak C}$ are replaced 
by some non-zero scalar multiples (see Remark~\ref{nearlyclosed}). Since the $\nu_{j}$, $j=0,\ldots,r$, 
are valuations, we see that property (b) in  Definition~\ref{quasidef} holds.

It remains to verify (a) in Definition~\ref{quasidef}. Let $-1\le j\le r$ be minimal such that
$\nu_{i}((g+h)_i)= \nu_{i}(g_i)=\nu_{i}(h_i)$ for all $i>j$. If $j=-1$, then obviously property (a) holds.

Suppose $j\ge 0$ and set $F_j=\prod_{\ell=j+1}^r f_{\ell}^{-N^{\ell-j}D_\ell}$,
where $F_r=1$ because we have an empty product.
Remark~\ref{nearlyclosed} implies that 
$(g+h)_j$ is a linear combination of functions of the form $g^th^{N^{r-j}-t}F_j$, $t=0,\ldots, N^{r-j}$.
A small calculation shows: for all $t=0,\ldots, N^{r-j}$, one has
\begin{equation}\label{firstformulaone}
\begin{array}{rcl}
\nu_{{j}}((g+h)_j)&\ge&\min\{\nu_{{j}}(g^th^{N^{r-j}-t}F_j)\mid t=0,\ldots, N^{r-j} \} \\
&=&\min\{ \frac{t}{N^{r-j}}\nu_{j}(g_j) +(1-\frac{t}{N^{r-j}})\nu_{j}(h_j)\mid t=0,\ldots, N^{r-j}\}.
\end{array}
\end{equation}
and hence:  $\nu_{{j}}((g+h)_j)\ge \min\{\nu_{{j}}(g_j),\nu_{{j}}(h_j)\}$. If the inequality is strict, then 
property (a) holds. 

If the inequality is not strict, then (by the assumption on $j$) we have $\nu_{{j}}(g_j)\not=\nu_{{j}}(h_j)$.
Without loss of generality assume $\nu_{{j}}(g_j)<\nu_{{j}}(h_j)$. But then:
$$
(g+h)_{j-1}=\left.\frac{(g+h)_{j}^{N}}{f_j^{ND_j}}\right\vert_{\hat X_{j-1}}=\left.\frac{g_{j}^{N}}{f_j^{ND_j}}\right\vert_{\hat X_{j-1}}=g_{j-1}
$$
because the restrictions of all the other terms in the expansion of  $(g+h)_{j-1}$  vanish. But this implies
$(g+h)_\ell=g_\ell$ for $0\le \ell <j$ and hence $\mathcal V_{\mathfrak C}(g+h)\ge \min\{\mathcal V_{\mathfrak C}(g), \mathcal V_{\mathfrak C}(h)\}$.
\end{proof}

\subsection{The lattice generated by the image of $\mathcal V_{\mathfrak C}$}\label{leaf:lattice}

The lattice $L^{\mathfrak C}$ introduced in Lemma~\ref{Lem:ImageValC} should be considered as a first approximation of the lattice generated by the 
image of $\mathcal{V}_{\mathfrak{C}}$. The valuation monoid $\mathbb{V}_{\mathfrak{C}}(X)=
\{\mathcal V_{\mathfrak C}(g)\mid g\in R\setminus\{0\}\}$ may be 
contained in a proper sublattice of $L^{\mathfrak C}$.

In this section we propose an approach to determine the sublattice $L^{\mathfrak C}_{\mathcal V}\subseteq L^{\mathfrak C}$ generated
by $\mathbb{V}_{\mathfrak{C}}(X)$. This strategy highlights the strong connection between the point of view in this article
and the usual procedure in the theory of Newton-Okounkov bodies. The difference between these approaches will only become evident in Section~\ref{nonegativequasival}.

We fix a maximal chain $\mathfrak{C}:p_r > p_{r-1} > \ldots > p_0$.

\begin{lemma}
There exist rational functions $F_r,\ldots,F_0\in\mathbb{K}(\hat{X})\setminus \{0\}$ such that
\[
\mathcal V_\mathfrak{C}(F_j) = (\underbrace{0,\ldots,0}_{r-j},1/b_j,*,\ldots,*),
\]
where the $*$ are certain numbers in $\mathbb{Q}$.
\end{lemma}

\begin{proof}
Suppose $r\geq j\geq 0$. By assumption, the variety $ \hat{X}_j$ is smooth in codimension 1, so the local ring
$\mathcal{O}_{\hat X_j,\hat X_{j-1}}$ is a discrete valuation ring. Let $\eta_j$ be a uniformizer in the maximal ideal.
It is a rational function on $\hat{X}_j$ with the property $\nu_j(\eta_j)=1$.

As a rational function on $\hat{X}_j$, $\eta_j$ can be represented as the restriction to $\hat X_j$ of a rational function on $V$: 
there exist $g,h\in \mathbb K[V]$ such that both do not identically vanish on $\hat X_j$, and $\eta_j=\frac{g}{h}\vert_{\hat X_j}$. In particular, $F_j:=\frac{g}{h}\vert_{\hat X}\in\mathbb K(\hat X)$ is a well defined rational function, which has neither poles nor vanishes on non-empty open subsets in $\hat X=\hat X_r$, $\hat X_{r-1}$, $\ldots$, $\hat X_j$. It follows that the first
$r-j$ entries in $\mathcal V_{\mathfrak C}(F_j)$ are equal to zero. One gets for the associated sequence of rational functions along $\mathfrak C$:
$$
(F_j)_{\mathfrak C}=(F_j,F_j^N\vert_{\hat X_{r-1}},\ldots, F_j^{N^{r-j}}\vert_{\hat X_{j}},\ldots),
$$
and hence the $(r-j+1)$-th entry is $\nu_j(F_j^{N^{r-j}}\vert_{\hat X_{j}})/(N^{r-j}b_j)=\frac{\nu_j(\eta_j)}{b_j}=\frac{1}{b_j}$.
\end{proof}

Let $L^{\mathfrak C}_{\mathcal V}\subseteq L^{\mathfrak C}$ be the sublattice generated by the valuation monoid
$\mathbb V_{\mathfrak C}(X)$.

\begin{proposition}
Let $F_r,\ldots,F_0\in\mathbb{K}(\hat{X})\setminus \{0\}$ be rational functions such that 
\[
\mathcal V_\mathfrak{C}(F_j) = (\underbrace{0,\ldots,0}_{r-j},1/b_j,*,\ldots,*),
\] 
where the $*$ are certain numbers in $\mathbb{Q}$. Then 
$L^{\mathfrak C}_{\mathcal V} = \langle \mathcal V_\mathfrak{C}(F_r),\ldots,\mathcal V_\mathfrak{C}(F_0)\rangle_\mathbb{Z}$.
\end{proposition}

\begin{proof}
Let $\bar L^{\mathfrak C}_{\mathcal V}=\langle \mathcal V_\mathfrak{C}(F_r),\ldots,
\mathcal V_\mathfrak{C}(F_0)\rangle_\mathbb{Z} \subseteq L^{\mathfrak C}$
be the lattice generated by the valuations of the rational functions $F_r,\ldots F_0$.
It is obvious that $\bar L^{\mathfrak C}_{\mathcal V}\subseteq L^{\mathfrak C}_{\mathcal V}$.

To prove the reverse inclusion, it suffices to show $\mathcal V_{\mathfrak C}(g)\in \bar L^{\mathfrak C}_{\mathcal V}$ for $g\in \mathbb K(\hat X)\setminus \{0\}$.
Fix a rational function $g$. Note that without loss of generality we 
can modify $g$ by multiplying it by a power of $F_j$ for some $j=0,\ldots,r$. Indeed, since $\mathcal V_\mathfrak{C}$ is a valuation, on has 
$\mathcal V_\mathfrak{C}(g)\in \bar L^{\mathfrak C}_{\mathcal V}$
if and only if  $\mathcal V_\mathfrak{C}(gF^a_j)\in \bar L^{\mathfrak C}_{\mathcal V}$ for some $a\in \mathbb Z$. 

We proceed by induction on the number of entries equal to zero at the beginning of $\mathcal V_\mathfrak{C}(g)$. Suppose the first entry
is non-zero: it is equal to $\frac{\nu_r(g)}{b_r}$. After replacing $g$ by $gF^{-\nu_r(g)}_r$, we can assume  
that $g\in\mathbb K(\hat{X})\setminus\{0\}$ is such that the first entry in $\mathcal V_\mathfrak{C}(g)$ is equal to zero.

Suppose $g$ is a non-zero rational function such that the first $r-j$ entries in $\mathcal V_\mathfrak{C}(g)$ are equal to zero. If $j=-1$, then
we are done. If $j\ge 0$, then $(g)_{\mathfrak C}=(g, g^N,\ldots, g^{N^{r-j}},\ldots)$, and the $(r-j+1)$-st entry in $\mathcal V_\mathfrak{C}(g)$
is $\nu_{j}(g^{N^{r-j}})/(N^{r-j}b_{j})$, which is equal to $\frac{\nu_{j}(g)}{b_{j}}$. So by multiplying $g$ with an appropriate power
of $F_{j}$, we get a new rational function having one more entry equal to zero in its valuation. 
\end{proof}

The rational functions $F_r,\ldots,F_0\in \mathbb K(\hat{X})$ used above are far from being unique. But all possible choices have one common 
feature. Let $A_{\mathfrak C}$ be the rational $(r+1)\times (r+1)$ matrix having as columns the valuations $\mathcal V_{\mathfrak C}(F_r),
\ldots \mathcal V_{\mathfrak C}(F_0)$, this is a lower triangular matrix. Let $B_{\mathfrak C}$ be its inverse. It is of the following form: 
\[
B_\mathfrak{C} = \left(
\begin{matrix}
b_r \\
* & b_{r-1}\\
\vdots & \ddots & \ddots\\
* & \cdots & * &  b_0\\
\end{matrix}
\right),
\]
where the $*$ are certain numbers in $\mathbb{Q}$.

\begin{proposition}
Let $v\in\mathbb Q^{\mathfrak C}$. Then $v\in L^{\mathfrak C}_{\mathcal V}$ if and only if $B_\mathfrak{C} \cdot v\in \mathbb{Z}^{r+1}$.
Moreover, the entries of $B_\mathfrak{C}$ are integers.
\end{proposition}

\begin{proof}
The first part is just a reformulation of the previous proposition. For the integrality property note that $e_j = \nu_\mathfrak{C}(f_j)$ is an element of 
$L^{\mathfrak C}_{\mathcal V}$, hence the $j$-th column of $B_\mathfrak{C}$ must have integral entries.
\end{proof}

\subsection{$\mathbb{K}^*$-action versus valuation}

The varieties $\hat X$ and $\hat X_p$ are all endowed with a $\mathbb K^*$-action, and the algebra $R$ as well as
the algebras $R_p$  are correspondingly graded. For $g\in R_p$ and $\lambda\in\mathbb K^*$ denote by $g^\lambda$ the 
function $g^\lambda(y):=g(\lambda y)$ for $y\in \hat X_p$. This $\mathbb{K}^*$-action on $R_p$ naturally extends to $\mathbb{K}(\hat{X}_p)$.
\begin{lemma}\label{homogeneousparts}
\begin{itemize}
\item[{i)}] For any $g\in R\setminus\{0\}$ and $\lambda\in \mathbb K^* : \mathcal V_{\mathfrak C}(g^\lambda)=\mathcal V_{\mathfrak C}(g)$.
\item[{ii)}] Let $h=h_{1}+\ldots +h_{t}\in R=\bigoplus_{i\ge 0} R(i)$ be a decomposition of $h\not=0$ into its homogeneous parts.
Then
$$
\mathcal V_{\mathfrak C}(h)=\min\{\mathcal V_{\mathfrak C}(h_j)\mid 1\leq j\leq t \text{ such that } h_j\not=0\}.
$$
\end{itemize}
\end{lemma}
\begin{proof}
The $\mathbb K^*$-action on $ \hat X_{{j}}$ stabilizes
the divisor $\hat X_{{j-1}}$. The associated algebra automorphisms of $R_{p_j}$ stabilize hence
the vanishing ideal of $\hat X_{{j-1}}$ as well as the local ring $\mathcal O_{\hat X_{{j}},\hat X_{{j-1}}}\subseteq 
\mathbb K(\hat X_{{j}})$ and its maximal ideal. So for all $g\in R_{p_j}$: $\nu_{j}(g^\lambda)=\nu_{j}(g)$.

For $g\in R\setminus\{0\}$ let $g_{\mathfrak C}=(g_r,\ldots,g_0)$ be the associated sequence of rational functions. 
We can use the $\mathbb K^*$-action to construct two new tuples: the $\lambda$-twisted
tuple ${}^\lambda g_{\mathfrak C}=(g^\lambda_r,\ldots,g^\lambda_0)$ obtained by twisting
component-wise each of the rational functions occurring in $g_{\mathfrak C}$. And we can
consider $g_{\mathfrak C}^\lambda=(g'_r,\ldots,g'_0)$, the tuple associated
to the function $g^\lambda$.

The functions $\{f_{p}\mid p\in A\}$ are homogeneous, so by Remark~\ref{nearlyclosed}
the $j$-th component $g^\lambda_j$ in the $\lambda$-twisted sequence ${}^\lambda g_{\mathfrak C}$
and the $j$-th component $g'_j$ in  $g_{\mathfrak C}^\lambda$ differ only by a scalar multiple.
It follows: $\mathcal V_{\mathfrak C}(g^\lambda)=\mathcal V_{\mathfrak C}(g)$.

Let $h=h_{1}+\ldots +h_{t}$ be a decomposition of $h\not=0$ into its homogeneous parts. We know:
$\mathcal V_{\mathfrak C}(h)\ge \min\{\mathcal V_{\mathfrak C}(h_j)\mid 1\leq j\leq t \text{ such that } h_j\not=0\}$
by the minimum property of a valuation (see Definition~\ref{quasidef} a)). To prove the equality, note that 
one can find pairwise distinct, non-zero scalars $\lambda_1,\ldots,\lambda_t\in\mathbb{K}^*$ 
such that the linear span of the following functions coincide:
$$
\langle h_1,\ldots,h_t\rangle_{\mathbb K} =\langle h^{\lambda_1},\ldots,h^{\lambda_t}\rangle_{\mathbb K}.
$$
So one can express the homogeneous function $h_i$ as a linear combination of the functions 
$h^{\lambda_1},\ldots,h^{\lambda_t}$, and hence by part i) of the lemma:
$$
\mathcal V_{\mathfrak C}(h_i)\ge  \min\{\mathcal V_{\mathfrak C}(h^{\lambda_j})\mid j=1,\ldots, t\}
=\mathcal V_{\mathfrak C}(h),
$$
and hence $\mathcal V_{\mathfrak C}(h)=\min\{\mathcal V_{\mathfrak C}(h_{j})\mid j=1,\ldots,t\}$.
\end{proof}

\subsection{Leaves}\label{leaves1}

Let $\gls{Lcdag}\subseteq L^{{\mathfrak C}}$ be the submonoid of tuples 
such that the first non-zero entry is positive. Let $g\in R\setminus\{0\}$ be a regular function on $\hat X$. The first non-zero entry of 
$\mathcal V_{\mathfrak C}(g)$ is a rescaled vanishing multiplicity, so $\mathcal V_{\mathfrak C}:R\setminus\{0\}\rightarrow L^{\mathfrak C}$ is a 
valuation with values in $L^{{\mathfrak C},\dag}$.
The elements in this submonoid satisfy an additional compatibility 
property: for any $\underline a,\underline b,\underline c\in L^{\mathfrak C,\dag}$ we have:

\begin{equation}\label{newinequality}
\text{if }\underline{a}\ge \underline b, \text{ then } \underline a + \underline c\ge \underline b+\underline c\ge \underline b.
\end{equation} 

This property ensures that the subspaces
$$
{R}_{\ge \underline a}^{\mathfrak{C}}:=\{ g\in   R\setminus\{0\}\mid  \mathcal V_{\mathfrak C}(g)\ge \underline a\}\cup\{0\},
$$
respectively
$$
{R}_{> \underline a}^{\mathfrak{C}}:=\{ g\in   R\setminus\{0\}\mid  \mathcal V_{\mathfrak C}(g)> \underline a\}\cup\{0\}
$$
for $\underline a\in  L^{{\mathfrak C},\dag}$ are ideals. It follows that the associated graded vector space
$$
\gls{grcR}=\bigoplus_{\underline a\in  L^{{\mathfrak C},\dag}} {  R}_{\ge \underline a}^{\mathfrak{C}}/{  R}_{> \underline a}^{\mathfrak{C}}
$$
is a $\mathbb K$-algebra. The subquotient ${R}_{\ge \underline a}^{\mathfrak{C}}/{R}_{> \underline a}^{\mathfrak{C}}$ for $\underline a\in  L^{{\mathfrak C},\dag}$
is called a \emph{leaf of the valuation} $\mathcal V_{\mathfrak C}$. 

According to Example~\ref{extremalvaluation}, $ \mathcal V_{\mathfrak C}$ is a full rank valuation, i.e. the lattice generated by $\mathbb V_{\mathfrak C}(X)$
has rank $\dim X+1$. As a consequence of the Abhyankar's inequality one obtains:

\begin{theorem}[{\cite{KM}}]\label{chainleavedimensiontwo}
The valuation $ \mathcal V_{\mathfrak C}$ has at most one-dimensional leaves, i.e. for any $\underline a\in  L^{{\mathfrak C},\dag}$, $\dim_{\mathbb{K}}{R}_{\ge \underline a}^{\mathfrak{C}}/{R}_{> \underline a}^{\mathfrak{C}}\leq 1$.
\end{theorem}

\begin{coro}\label{chainleavedimensionone}
Let $g,h\in R\setminus\{0\}$. Assume that for a maximal chain $\mathfrak C\subseteq A$, $\mathcal V_{\mathfrak C}(g)=\mathcal V_{\mathfrak C}(h)$ holds. Then there exists $\lambda\in\mathbb K^*$ and $h'\in R$ such that $g=\lambda h+h'$. If $h'\neq 0$, then $\mathcal V_{\mathfrak C}(h')>\mathcal V_{\mathfrak C}(h)=\mathcal V_{\mathfrak C}(g)$ holds.
\end{coro}

\begin{proof}
This is just a reformulation of Theorem~\ref{chainleavedimensiontwo}.
\end{proof}

\begin{lemma}\label{powerrelation3}
Let $\mathfrak C\subseteq A$ be a maximal chain and let $g_1,g_2,g_3,g_4\in R\setminus\{0\}$. If
$\mathcal V_{\mathfrak C}(g_1)+\mathcal V_{\mathfrak C}(g_2)=\mathcal V_{\mathfrak C}(g_3)+\mathcal V_{\mathfrak C}(g_4)$,
then there exist $\lambda\in\mathbb K^*$ and $h'\in R$ such that $g_1g_2=\lambda g_3g_4+ h'$. If $h'\neq 0$, then $\mathcal V_{\mathfrak C}(h')>\mathcal V_{\mathfrak C}(g_1)+\mathcal V_{\mathfrak C}(g_2)$ holds. If the functions $g_1,g_2,g_3,g_4$ are homogeneous, so is $h'$. 
\end{lemma}

\begin{proof}
For the first statement it suffices to apply Corollary~\ref{chainleavedimensionone} to $g=g_1g_2$ and $h=g_3g_4$. 
Assume that $g_1,g_2,g_3,g_4$ are homogeneous, and let $h'=h_1+h_2+\ldots+h_t$, $h_k\neq 0$ for $k=1,\ldots,t$, be a decomposition of $h'$ into its 
homogeneous components. If $g_1g_2$ and $g_3g_4$ are not of the same degree, then the equality  $g_1g_2=\lambda g_3g_4+ h'$ is only possible if $\lambda g_3g_4=-h_j$ for some $1\leq j\leq t$. By Lemma \ref{homogeneousparts}, this is impossible because
$$\mathcal V_{\mathfrak C} (h_j)\geq \min\{\mathcal{V}_{\mathfrak{C}}(h_\ell)\mid \ell=1,\ldots,t\}=\mathcal{V}_{\mathfrak{C}}(h')>\mathcal V_{\mathfrak C} (g_1g_2)=\mathcal V_{\mathfrak C} (g_3g_4).$$
Therefore $g_1g_2$ and $g_3g_4$ are necessarily of the same degree,
and so is $h'$.
\end{proof}

Lemma~\ref{powerrelation3} suggests that the structure of $\mathrm{gr}_{\mathfrak C}R$ is very similar to that of the semigroup algebra $\mathbb K[\mathbb V_{\mathfrak C}(X)]$. Indeed, it follows from Definition \ref{quasidef} (c') that $\mathrm{gr}_{\mathfrak C}R$
has no zero divisors. Applying {\cite[Remark 4.13]{BG}} (see also {\cite[Proposition 2.4]{KM}}) gives:

\begin{coro}\label{nearlytoric}
As $\mathbb{K}$-algebra, $\mathrm{gr}_{\mathfrak C}R$ is isomorphic to $\mathbb K[\mathbb V_{\mathfrak C}(X)]$.
\end{coro}

\section{Localization and finite generation}\label{Sec:FiniteGen}

In this section we introduce the core of the valuation $\mathcal{V}_{\mathfrak{C}}$ associated to a maximal chain $\mathfrak{C}$: it is a finitely generated submonoid of the valuation monoid $\mathbb{V}_{\mathfrak{C}}(X)$. Further results on the relation between the core and standard monomials will be given in Section \ref{fanandcone}.

\begin{definition}
The \emph{core} $\gls{PCX}$ of the valuation monoid $\mathbb{V}_\mathfrak{C}(X)$ is defined as its intersection with the positive orthant: 
$$
P_{\mathfrak C}(X):=\mathbb V_{\mathfrak C}(X)\cap \mathbb Q_{\ge 0}^{\mathfrak C}.
$$
\end{definition}

As the intersection of two monoids, $P_{\mathfrak C}(X)$ is a monoid. By Example~\ref{extremalvaluation} we know $\mathcal V_{\mathfrak C}(f_{i})= e_{i}$
and hence $ \mathbb{N}^\mathfrak{C}\subseteq P_{\mathfrak C}(X)$. So the monoid has an additional structure, it is endowed with a natural $\mathbb{N}^\mathfrak{C}$-action:
\begin{equation}\label{NC-action}
\mathbb N^{\mathfrak C}\times P_{\mathfrak C}(X)\rightarrow P_{\mathfrak C}(X),\ \ (\underline n, \underline{k})\mapsto \underline n\circ \underline{k}:=\underline{n}+\underline{k}.
\end{equation}
\begin{lemma}\label{finitegenoverN}
The core $P_{\mathfrak C}(X)$ is a finitely generated $\mathbb N^{\mathfrak C}$-module.
\end{lemma}

\begin{proof}
For the proof we use Dickson's Lemma (\cite{Dick}, Lemma A) which states (the formulation has been adapted to our situation): every monomial ideal in the polynomial ring
$\mathbb K[y_p\vert p\in \mathfrak C]$ is finitely generated. By the standard bijection $\underline{n}=(n_p)_{p\in\mathfrak C}\mapsto
\underline{y}^{\underline{n}}=\prod_{p\in\mathfrak C} y_p^{n_p}$ between $\mathbb N^{\mathfrak C}$
and the monomials in $\mathbb K[y_p\vert p\in \mathfrak C]$ we get a bijection between monomial ideals in $\mathbb K[y_p\vert p\in \mathfrak C]$
and $\mathbb N^{\mathfrak C}$-submodules in $\mathbb N^{\mathfrak C}$, where the latter is acting on itself by addition. And 
Dickson's Lemma can be reformulated as:  every $\mathbb N^{\mathfrak C}$-submodule $M\subseteq \mathbb N^{\mathfrak C}$ is finitely generated
as $\mathbb N^{\mathfrak C}$-module.

Since $e_r,\ldots,e_1,e_0\in L^\mathfrak{C}$, adding elements 
of $\mathbb Z^{\mathfrak C}$: $\underline{n}\circ \underline{\ell}:=\underline{n}+ \underline{\ell}$ defines an action of $\mathbb Z^{\mathfrak C}$ 
on the lattice $L^\mathfrak{C}$. We get an induced action by $\mathbb N^{\mathfrak C}$ on the intersection 
$L^\mathfrak{C}\cap\mathbb{Q}_{\geq 0}^\mathfrak{C}$, which by \eqref{NC-action}
stabilizes the submonoid $P_{\mathfrak C}(X)$.

We decompose $L^\mathfrak{C}$ into a finite number of $\mathbb Z^{\mathfrak C}$-submodules (with respect to the action ``$\circ$'' defined above).
Since $L^\mathfrak{C}$ is a lattice, the intersection
$$
Q:=L^{\mathfrak C}\cap\{\underline{a}=(a_r,\ldots,a_0)\in \mathbb Q_{\ge 0}^{\mathfrak C}\mid 0 \le a_j < 1\}
$$ 
is a finite set. By construction we have $L^\mathfrak{C}=\bigoplus_{\underline{a}\in Q} \mathbb Z^{\mathfrak C}\circ \underline{a}$ and correspondingly $L^\mathfrak{C}\cap\mathbb{Q}_{\geq 0}^\mathfrak{C}=\bigoplus_{\underline{a}\in Q} \mathbb N^{\mathfrak C}\circ \underline{a}$.
It follows: 
$$ 
P_{\mathfrak C}(X)=  \bigoplus_{\underline{a}\in Q} \big(\mathbb N^{\mathfrak C}\circ \underline{a}\cap P_{\mathfrak C}(X)\big)
 \subseteq  L^\mathfrak{C}\cap\mathbb{Q}_{\geq 0}^\mathfrak{C}.
$$
Since $Q$ is a finite set, to show that $P_{\mathfrak C}(X)$ is a finitely generated  $\mathbb N^{\mathfrak{C}}$-module
it suffices to show that the intersection $(\mathbb N^{\mathfrak C}\circ \underline{a})\cap P_{\mathfrak C}(X)$
is a finitely generated $\mathbb N^{\mathfrak{C}}$-module for all $\underline{a}\in Q$. As $\mathbb N^{\mathfrak C}$-modules,
$\mathbb N^{\mathfrak C}\circ \underline{a}$  and  $\mathbb N^{\mathfrak C}$ are isomorphic, and 
$(\mathbb N^{\mathfrak C}\circ \underline{a})\cap P_{\mathfrak C}(X)$ is an intersection of submodules and hence a submodule.
So one can apply (the reformulated version of) Dickson's Lemma and hence $(\mathbb N^{\mathfrak C}\circ \underline{a})\cap P_{\mathfrak C}(X)$
is finitely generated as $\mathbb N^{\mathfrak C}$-module, which in turn implies $P_{\mathfrak C}(X)$ is a finitely generated  $\mathbb N^{\mathfrak{C}}$-module.
\end{proof}

We introduce an open affine subset of $\hat{X}$ to pinch the valuation monoid $\mathbb{V}_{\mathfrak{C}}(X)$.

Let $U_{\mathfrak C}$ be the open affine subset of $\hat X$ defined by
$$
\gls{UC}=\{x\in \hat X\mid \prod_{p\in\mathfrak C}f_p(x)\not=0\}.
$$
Its coordinate ring will be denoted by $\mathbb K[U_{\mathfrak C}]$: it is the
localization of $\mathbb K[\hat X]$ at $\prod_{p\in\mathfrak C}f_p$.

The valuation $\mathcal{V}_\mathfrak{C}:\mathbb{K}(\hat{X})\setminus\{0\}\to L^\mathfrak{C}$ induces by restriction a valuation $\mathcal V_{\mathfrak C}:\mathbb K[U_{\mathfrak C}]\setminus\{0\}\to L^{\mathfrak C}$
which has one-dimensional leaves. Let $\gls{VCUC} :=\{\mathcal V_{\mathfrak C}(g)\mid g\in 
\mathbb K[U_{\mathfrak C}]\setminus\{0\}\}\subseteq L^{\mathfrak C}$ denote the associated valuation monoid. These monoids are contained in each other:
$P_{\mathfrak C}(X)\subseteq \mathbb{V}_{\mathfrak C}(X)\subseteq\mathbb{V}_{\mathfrak C}(U_{\mathfrak C})$.
The core $P_{\mathfrak C}(X)$ can be thought of as a condensed version of the valuation monoid $\mathbb V(U_{\mathfrak C})$:

\begin{lemma}\label{finGenr}
We have 
$$\mathbb{V}_{\mathfrak C}(U_{\mathfrak C})=P_{\mathfrak C}(X)+\mathbb{Z}^\mathfrak{C}=\{\underline{a}+\underline{m}\mid \underline{a}\in P_{\mathfrak C}(X),\ \underline{m}\in\mathbb Z^{\mathfrak C}\}.$$
In particular, as $\mathbb Z^{\mathfrak C}$-module, $\mathbb V(U_{\mathfrak C})$ is finitely generated.
\end{lemma}

\begin{proof}
For $g\in R\setminus\{0\}$, there exists an $m\in \mathbb N$ such that  
$$
\mathcal V_{\mathfrak C}(g (\prod_{p\in\mathfrak C}f_p)^m) =\mathcal V_{\mathfrak C}(g)+\sum_{p\in\mathfrak C} me_p\in P_{\mathfrak C}(X).
$$
Hence every element in $\mathbb K[U_{\mathfrak C}]\setminus\{0\}$ can be written as a quotient $\frac{h}{(\prod_{p\in\mathfrak C}f_p)^n}$, $h\in R\setminus\{0\}$, 
in such a way that $\mathcal V_{\mathfrak C}(h)\in P_{\mathfrak C}(X)$, which proves the claim.
\end{proof}

\begin{coro}\label{powerrelation}
Let $g\in R\setminus\{0\}$.
There exist $m\in\mathbb N$, $\lambda\in \mathbb K^*$ and an element $g'\in R$
with $\mathcal V_{\mathfrak C}(g')>m \mathcal V_{\mathfrak C}(g)$ as long as $g'\neq 0$, such that in $\mathbb K[U_{\mathfrak C}]$ we have
$$
g^m=\lambda f_{p_r}^{ma_r} \cdots f_{p_0}^{ma_0}+g'.
$$
If $g$ is homogeneous, then so are $ f_{p_r}^{ma_r} \cdots f_{p_0}^{ma_0}$ and $g'$, and they are of the same degree.
\end{coro}

\begin{proof}
Fix $m\in \mathbb Z_{>0}$ such that $m\mathcal V_{\mathfrak C}(g)=(ma_r,\ldots,ma_0)\in \mathbb Z^{\mathfrak C}$.
It follows by Example~\ref{extremalvaluation}: $\mathcal V_{\mathfrak C}(g^m)=m\underline{a}= \mathcal V_{\mathfrak C}( f_{p_r}^{ma_r} \cdots f_{p_0}^{ma_0})$.
The remaining part of the proof is the same as Lemma~\ref{powerrelation3}, where we only used the fact that the leaves are one dimensional.
\end{proof}

As a consequence we get the following formula recovering the degree of a homogeneous element from its valuation (one easily verifies it in Example \ref{grasstwofour3}):

\begin{coro}\label{degreerelation}
If $g\in R\setminus\{0\}$ is homogeneous and $\mathcal V_{\mathfrak C}(g)=\underline a$, then
$$
\deg g=a_r\deg f_{p_r} +\ldots+a_1\deg f_{p_1}+a_0\deg f_{p_0}.
$$
\end{coro}

\begin{proof}
By Corollary~\ref{powerrelation}, we can find an $m\in\mathbb N$ such that 
$g^m=\lambda f_{p_r}^{ma_r} \cdots f_{p_0}^{ma_0}+g'$ for some $\lambda\in\mathbb K^*$, 
and all are of the same degree. It follows that 
the degree of $g$ is $\deg g=\frac{1}{m}\deg(f_{p_r}^{ma_r} \cdots f_{p_0}^{ma_0})$, which proves the claimed formula.
\end{proof}

The valuation $\mathcal V_{\mathfrak C}$ induces a filtration on $\mathbb K[U_{\mathfrak C}]$, let 
$\mathrm{gr}_{\mathfrak C}\mathbb K[U_{\mathfrak C}]$ be the associated graded algebra. 
The same arguments as for Corollary~\ref{nearlytoric} imply: the latter is isomorphic to the 
semigroup algebra $\mathbb K[\mathbb V_{\mathfrak C}(U_{\mathfrak C})]$,
and hence by Lemma~\ref{finitegenoverN} and \ref{finGenr}, this algebra is finitely generated and integral.

\begin{coro}\label{assgrad1}
$\mathrm{Spec}(\mathrm{gr}_{\mathfrak C}\mathbb K[U_{\mathfrak C}])\simeq \mathrm{Spec}(\mathbb K[\mathbb V_{\mathfrak C}(U_{\mathfrak C})])$
is a toric variety.
\end{coro}

We finish this section by taking a different point of view: looking at $U_{\mathfrak{C}}$ as being fibered over a torus.

\begin{coro}
The map 
$$
\phi_{\mathfrak C}: U_{\mathfrak C}\rightarrow (\mathbb K^*)^{r+1},\quad
u\mapsto (f_{p_r}(u),\ldots, f_{p_0}(u))
$$ 
is a finite morphism.
\end{coro}

\begin{proof}
According to Example \ref{extremalvaluation}, the valuations $\mathcal{V}_{\mathfrak{C}}(f_{p_r}),\ldots,\mathcal{V}_{\mathfrak{C}}(f_{p_0})$ are linearly independent, hence the functions $f_{p_r},\ldots,f_{p_0}$ are algebraically independent. The induced map between the coordinate rings 
$\phi^*_{\mathfrak C}:\mathbb K[x_r^{\pm 1},\ldots,x_0^{\pm 1}]\rightarrow \mathbb K[U_{\mathfrak C}]$
sending $x_j$ to $f_{p_j}$ is therefore injective and, by Lemma \ref{finGenr} and Corollary \ref{powerrelation}, it makes $\mathbb K[U_{\mathfrak C}]$
into a finite $\mathbb K[x_r^{\pm 1},\ldots,x_0^{\pm 1}]$-module.
\end{proof}

\section{Globalization: a non-negative quasi-valuation}\label{nonegativequasival}

To better understand the role of the finitely generated submonoid $P_{\mathfrak C}(X)$ of $\mathbb V_{\mathfrak C}(X)$, we consider in the following all the valuations $\mathcal V_{\mathfrak C}$ at once, where $\mathfrak{C}$ runs over all maximal chains in $A$. To do so, let $\mathbb Q^{A}$ be the $\mathbb Q$-vector space spanned by the basis elements $\{e_q\mid q\in A\}$. If $\mathfrak C$ is a maximal chain in $A$, we identify $\mathbb Q^{\mathfrak C}$ with the subspace of $\mathbb Q^{A}$ spanned by the basis elements  $\{e_p\mid p\in \mathfrak C\}$. To be able to compare for a given $g\in R\setminus\{0\}$  the various valuations  $\mathcal V_{\mathfrak C}(g)\in \mathbb Q^{\mathfrak C}\subseteq  \mathbb Q^{A}$, we need an order on $\mathbb Q^{A}$.

To define a total order on  $\mathbb Q^{A}$, fix a total order ``$\gls{geqt}$'' on $A$ refining the given partial order, and such that $\ell(p)>\ell(q)$ (see Definiton~\ref{length}) implies $p>^tq$\footnote{This second condition on the length is in fact not necessary. Results in the article hold with this condition removed. For details, see \cite[Section 2.6]{CFL4}.}. Such a total order exists since $A$ is a graded poset. Let 
\begin{equation}\label{totorder}
q_M >^t q_{M-1} >^t \ldots >^t q_0
\end{equation}
be an enumeration of the elements of $A$ depicting the total order.

Writing $v=(a_M,\ldots,a_0)$ for the vector $v=\sum_{j=0}^M a_je_{q_j}\in \mathbb Q^{A}$, we endow $ \mathbb Q^{A}$ with the lexicographic order as total order. This total order is compatible with the addition of vectors. 

We will denote by $\gls{CC}$ the set of all maximal chains in $A$.

\subsection{A non-negative quasi-valuation}
By Lemma \ref{minquasi}, taking the minimum value over a finite number of valuations defines a quasi-valuation. Recall that we think of $\mathbb Q^{\mathfrak C}$ as the subspace of $\mathbb Q^{A}$ spanned by the $e_p$, $p\in \mathfrak C$. So it makes sense to write $\mathcal V_{\mathfrak C}(g)\in \mathbb Q^{A}$ for a regular function $g\in R\setminus\{0\}$. Note that the total order on 
$\mathbb Q^{\mathfrak C}$ defined in Definition~\ref{orderchain} coincides with the total order on $\mathbb Q^{\mathfrak C}$ induced by that defined on $\mathbb Q^A$.

\begin{definition}
\begin{enumerate}
\item We define the quasi-valuation $\mathcal V:R\setminus\{0\}\rightarrow \mathbb Q^{A}$ by
$$
\gls{Vg}:=\min\{\mathcal V_{\mathfrak C}(g) \mid \mathfrak C\in\mathcal{C}\}.
$$
\item For $g\in R\setminus\{0\}$ with $\mathcal V(g)=(a_p)_{p\in A}$, the \textit{support} of $g$ is defined by
$$\gls{suppg}:=\{p\in A\mid a_p\not=0\}.$$ 
\end{enumerate}
\end{definition}

\begin{rem}
Let $g\in R\setminus\{0\}$. Unless $\supp\mathcal V(g)$ is a maximal chain, there might be several
maximal chains $\mathfrak{C}$ such that $\mathcal V(g)=\mathcal V_{\mathfrak C}(g)$.
\end{rem}

As an example let us consider an extremal function.
\begin{lemma}\label{extremalvaluationII}
For any $q\in A$, $\mathcal V(f_q)=e_q$.
\end{lemma}
\begin{proof}
By  Example~\ref{extremalvaluation} we know $\mathcal V_{\mathfrak C}(f_q)=e_q$ as long as $q\in \mathfrak C=(p_r,\ldots,p_0)$.
If $q\not\in \mathfrak C$, then let $p_k\in \mathfrak C$ be the unique element such that $\ell(p_k)=\ell(q)$. Since $p_k$ and $q$ are not comparable with 
respect to the partial order on $A$, $f_q$ vanishes on $X_{p_k}$. But $f_q$ is not the zero function, so there exist elements 
$p_{j}>p_{j-1}\ge p_k$ in $\mathfrak C$
such that $p_{j}$ covers $p_{j-1}$, $f_q\vert_{X_{p_{j}}}\not\equiv 0$, but $f_q$ vanishes on $X_{p_{j-1}}$.
It follows (compare with Example~\ref{fpitupel}):
$$
(f_q)_{\mathfrak C}=(f_q,f_q^N,\ldots,f_q^{N^{r-j}},\ldots)\Longrightarrow \mathcal V_{\mathfrak C}(f_q)
=\frac{\nu_{p_{j},p_{j-1}}(f_q^{N^{r-j}})} {N^{r-j}b_{p_{j},p_{j-1}}} e_{p_j} +\sum_{i<j} a_ie_{p_i}
$$
for some rational numbers $a_i\in \mathbb Q$, $0\le i\le j-1$.
Since $\nu_{p_j,p_{j-1}}(f_q)>0$ and $\ell(p_j)>\ell(q)$, this implies $\mathcal V_{\mathfrak C}(f_q)>e_q$,
and hence $\mathcal V(f_q)=e_q$.
\end{proof}

\begin{rem}
In general, the value $\mathcal V(g)$ for $g\in R\setminus\{0\}$ depends on the choice of the total order $>^t$ in the construction of $\mathcal V$. According to the lemma, for $q \in A$, $\mathcal V(f_q)$ is independent of the choice.
\end{rem}

\begin{example}\label{Ex:Bertini2}
Consider the generic hyperplane stratification introduced in Section \ref{Sec:Generic}. We study the quasi-valuation $\mathcal{V}$ on some particular functions.

First we compute $\mathcal{V}(g_{k})$ for the function $g_{k}$, where $1\leq k\leq s$, defined in Example \ref{Ex:Bertini1}. Let $\mathfrak C=(p_r,\ldots,p_0)$ be a maximal chain in $A$.  
\begin{enumerate}
\item If $p_0\neq q_{0,k}$, the function $g_k$ vanishes in $X_{p_0}$. This implies that the support $\supp\mathcal V_{\mathfrak C}(g_k)$ contains an element $p_k$ for some $k\geq 1$.
\item If $p_0=q_{0,k}$, then the non-vanishing of $g_k$ in $X_{q_{0,k}}$ implies that $\supp\mathcal V_{\mathfrak C}(g)=\{p_0\}$. 
\end{enumerate}

It follows for the quasi-valuation: $\supp \mathcal V(g_k)=\{p_0\}$, and $\mathcal V_{\mathfrak C}(g_k)=\mathcal V(g_k)$
if and only if $q_{0,k}\in\mathfrak C$. In this case we have $g_{\mathfrak{C}}=(g,g^N,\ldots,g^{N^{r}})$. Recall
that the valuation $\nu_{0}$ in the last step is just the degree of a corresponding homogeneous function (see Remark \ref{extrabond}). We obtain: 
$$
\mathcal V(g_k)=\mathcal V_{\mathfrak C}(g_k)=\frac{N^r(s-1)}{N^rb_{p_0,p_{-1}}}e_0=\frac{s-1}{\deg f_{0,k}}e_0.
$$

Then we consider a linear function $h\in V^*$ which does not vanish in any of the points $\varpi_1,\ldots,\varpi_s$, the same arguments as above imply: if $q_{0,k}\in\mathfrak C$, then
$$
\mathcal V(gh)=\mathcal V_{\mathfrak C}(gh)=\frac{s}{\deg f_{0,k}}e_0.
$$
\end{example}

The following non-negativity property of the quasi-valuation will be crucial in determining relations in $R$ in the spirit of Corollary \ref{powerrelation} (see Corollary \ref{powerrelation2}).

\begin{proposition}\label{positivity}
For all $g \in R\setminus\{0\}$, $\mathcal V(g)\in \mathbb Q_{\ge 0}^{A}$.
\end{proposition}

\begin{proof}
For $\mathcal V(g)= \mathfrak m=(m_p)_{p\in A}$ let $\mathfrak C=(p_r,\ldots,p_0)$ be a maximal chain such that $\mathcal V_{\mathfrak C}(g)= \mathfrak m$.
Then $m_p=0$ for all $p\not\in\mathfrak C$. We denote by $\mathfrak m_{\mathfrak C}$ the vector $(m_{p_r},\ldots,m_{p_0})\in \mathbb Q^{r+1}$.
To prove the proposition, it remains to show that $\mathfrak m_{\mathfrak C}\in \mathbb Q_{\ge 0}^{r+1}$.

Let $g_{\mathfrak C}=(g_r,\ldots, g_0)$ be the sequence of rational functions associated to $g$ and $\mathfrak C$.
Since $g$ is regular, we know $\nu_{p_r,p_{r-1}}(g)\ge 0$, and hence $m_{p_r}\ge 0$. We proceed by decreasing induction.

Assume that for some $1\le j \le r$, we know that $g_{j}\in \mathbb K(\hat X_{p_j})$ is integral over $\mathbb K[\hat X_{p_j}]$ and hence $m_{p_{j}}\ge 0$. If $q\in A$ is covered by $p_{j}$, then we can find a maximal chain of the form $\mathfrak C'=(p_r,\ldots,p_{j},q,\ldots)$. By the minimality assumption we know $\mathcal V_{\mathfrak C}(g)\le \mathcal V_{\mathfrak C'}(g)$.
In $\mathfrak C$ and $\mathfrak C'$, the entries with indexes $r,r-1,\ldots,j$ coincide, which implies by the minimality assumption:
$$ 
\frac{\nu_{p_{j},p_{j-1}}(g_{j})}{b_{p_{j},p_{j-1}}}\le \frac{\nu_{p_{j},q}(g_{j})}{b_{p_{j},q}}.
$$ 
By Proposition~\ref{normularinherited}, the rational function $g_{j-1}\in \mathbb K(\hat X_{p_{j-1}})$ is hence integral over 
$\mathbb K[\hat X_{p_{j-1}}]$ and thus $m_{p_{j-1}}\ge 0$, which finishes the proof by induction. 
\end{proof}

\subsection{A characterization via support}
Given a function $g\in R\setminus\{0\}$,
let $\gls{Cg}\subseteq \mathcal C$ be the set of maximal chains along which the quasi-valuation attains its minimum:
$$
\mathcal C(g):=\{ \mathfrak C\in \mathcal C\mid \mathcal V_{\mathfrak C}(g)=\mathcal V (g)\}.
$$
 
\begin{proposition}\label{supportcondition}
The set $\mathcal{C}(g)$ consists of maximal chains in $A$ containing $\mathrm{supp}\mathcal{V}(g)$, i.e. $\mathfrak C\in \mathcal C(g)$ if and only if $\supp\mathcal V(g)\subseteq \mathfrak C$.
\end{proposition}

\begin{proof}
If $\mathfrak C\in \mathcal C(g)$, then $\supp \mathcal V(g)=\supp \mathcal V_{\mathfrak C}(g)\subseteq \mathfrak C$.

To prove the opposite inclusion, let $\mathfrak{C}'=(p_r',\ldots,p_0')$ be a maximal chain in $A$ satisfying $\mathrm{supp}\mathcal{V}(g)\subseteq\mathfrak{C'}$, we show that $\mathcal{V}(g)=\mathcal{V}_{\mathfrak{C}'}(g)$, hence $\mathfrak C'\in \mathcal C(g)$. We proceed by induction on the length of the maximal chain in a poset. If all maximal chains in $A$ have length $0$, then $A$ has only one element and there is nothing to prove.

In the inductive procedure, the functions showing up are not necessarily regular
but those in $\mathbb K(\hat X)$ which are integral over $\mathbb K[\hat X]$, so we start with such a function $g$ and keep in mind  the setup in Section~\ref{normalization}. 

Let $\mathfrak C=(p_r,\ldots,p_0)$ be a maximal chain such that 
$\mathfrak C\in \mathcal C(g)$, i.e. $\mathcal{V}(g)=\mathcal{V}_{\mathfrak{C}}(g)$. We will denote $g_{\mathfrak{C}}=(g_r,\ldots,g_0)$, $g_{\mathfrak{C}'}=(g_r',\ldots,g_0')$, $\mathcal{V}(g)=(a_r,\ldots,a_0)$ and $\mathcal{V}_{\mathfrak{C}'}(g)=(a_r',\ldots,a_0')$. 

Both maximal chains $\mathfrak{C}$ and $\mathfrak{C}'$ start with $p_{\max}=p_r=p_r'$. We consider the following cases:
\begin{enumerate}
\item Assume that $p_{\max}\not\in \supp\mathcal V_{\mathfrak C}(g)$, then let
$0\leq j<r$ be maximal such that $p_j\in \supp\mathcal V_{\mathfrak C}(g)$. It follows that $a_r=\ldots=a_{j+1}=0$
and the sequence of function is $g_{\mathfrak C}=(g,g^N,\ldots,g^{N^{r-j}} ,\ldots)$. In particular, $g\vert_{\hat X_{p_j}} $
is a rational function, integral over $\mathbb K[\hat X_{p_j}]$, which does not vanish on an open and dense subset of $\hat X_{p_j}$.
Since $\mathrm{supp}\mathcal{V}(g)\subseteq\mathfrak{C}'$, $p_j=p_j'$ and hence $\hat X_{p_j}=\hat X_{p'_j}$. This  implies that $g$ does not vanish on an open and dense subset of the subvarieties
$\hat X_{p'_k}$ for $k\ge j$. It follows that $a'_r=\ldots=a'_{j+1}=0$ and $g_{\mathfrak C'}=(g,g^N,\ldots,g^{N^{r-j}} ,\ldots)$.
Replacing $X$ by $X_{p_j}$, $g$ by $g^{N^{r-j}}$, and considering the subposet $A_{p_j}$ whose maximal chains have smaller length, we can proceed by induction.
\item Now we assume that $p_{\max}\in\supp\mathcal V_{\mathfrak C}(g)$. There are two cases to consider:
\begin{enumerate}
\item $p_{r-1}\in\supp\mathcal V_{\mathfrak C}(g)$: In this case, $p_{\max}$ and $p_{r-1}$ are contained in both $\mathfrak C$ and $\mathfrak C'$. Recall that the function $g_{r-1}$ satisfies the desired assumption 
by Proposition~\ref{normularinherited} because $\mathfrak C\in \mathcal{C}(g)$.
Replacing $X$ by $X_{p_{r-1}}$ and $g$ by $g_{r-1}$, we can proceed by induction.
\item $p_{r-1}\notin\supp\mathcal V_{\mathfrak C}(g)$: Let $0\leq j<r-1$ be maximal such that $p_j\in \supp\mathcal V_{\mathfrak C}(g)$. From this assumption we know $a_{r-1}=\ldots=a_{j+1}=0$. We look at the functions $g_{r-1}$ and $g_{r-1}'$:
$$
g_{r-1}=g^N/f_{p_r}^{Na_r}, \quad g'_{r-1}=g^N/f_{p_r}^{Na'_r},
$$
then the sequence of functions is of the form
$$g_{\mathfrak C}=(g,g_{r-1},g_{r-1}^N,\ldots,g_{r-1}^{N^{r-j-1}},\ldots).$$  
In particular, $g_{r-1}^{N^{r-j-1}}$, 
and hence $g_{r-1}\vert_{\hat X_{p_j}} $ is a rational function which is integral over $\mathbb K[\hat X_{p_j}]$, and it does not vanish on an open and dense subset of $\hat X_{p_j}$. 

From $\mathfrak C\in \mathcal{C}(g)$ one knows $a_r\le a'_r$. If the inequality were strict, then $g_{r-1}$
would vanish on the open and dense subset of $\hat X_{p'_{r-1}}$ where $g_{r-1}$ is defined. 
Notice that $\hat X_{p'_{r-1}}\supseteq \hat X_{p'_j}=\hat X_{p_j}$, and 
$g_{r-1}$ does not vanish identically on the latter set, we get a contradiction. 

It follows that $a_r= a'_r$ and hence $g_{r-1}=g'_{r-1}$. Since $g_{r-1}=g'_{r-1}$ does not vanish on an open and dense subset of 
$\hat X_{p'_j}=\hat X_{p_j}$, the inclusions 
$$\hat X_{p'_j}\subseteq \hat X_{p'_{j+1}}\subseteq\ldots\subseteq \hat X_{p'_{r-2}}$$ 
imply $a'_{r-1}=\ldots=a'_{j+1}=0$, and hence 
$g_j=g_{r-1}^{N^{r-j-1}}=g_j'$. Now replacing $X$ by $X_{p_j}$, $g$ by $g_{r-1}^{N^{r-j-1}}$ and consider the subposet $A_{p_j}$, we can proceed again by induction.
\end{enumerate}
\end{enumerate}
\end{proof}

\subsection{Some consequences}
We give some consequences of Proposition \ref{supportcondition} which will be used later. As an immediate corollary of Proposition~\ref{supportcondition}, we have:

\begin{coro}\label{schnittlemma}
If $g,h\in R\setminus\{0\}$, then $\mathcal C(g)\cap \mathcal C(h)\subseteq \mathcal C(gh)$.
\end{coro}

\begin{proof}
If $\mathfrak C\in \mathcal C(g)\cap \mathcal C(h)$, then
$\mathcal V_{\mathfrak C}(gh)=\mathcal V_{\mathfrak C}(g)+\mathcal V_{\mathfrak C}(h)=\mathcal V(g)+\mathcal V(h)$.
This implies $\mathcal V(gh)\ge \mathcal V(g)+\mathcal V(h)=\mathcal V_{\mathfrak C}(gh)$ and hence 
$\mathcal V(gh)=\mathcal V_{\mathfrak C}(gh)$, i.e. $\mathfrak C\in \mathcal C(gh)$.
\end{proof}

Since $\mathcal V$ is a quasi-valuation, the inequality $\mathcal V(gh) \ge \mathcal V(g)+\mathcal V(h)$ holds. We can be more precise:

\begin{proposition}\label{quasivaluationA}
For $g,h\in R\setminus\{0\}$, 
$\mathcal V(gh) = \mathcal V(g)+\mathcal V(h)$ if and only if $\mathcal C(g)\cap \mathcal C(h)\not=\emptyset.$
\end{proposition}

\begin{proof}
For a maximal chain $\mathfrak C$ one has:
\begin{equation}\label{schnittlemmadiagramm}
\begin{matrix}
\mathcal V_{\mathfrak C}(gh)&\ge&\mathcal V(gh)\\
\raisebox{0.6em}{\rotatebox{270}{$=$}}&& \raisebox{0.6em}{\rotatebox{270}{$\ge$}}\\
\mathcal V_{\mathfrak C}(g)+\mathcal V_{\mathfrak C}(h)&\ge&\mathcal V_{}(g)+\mathcal V_{}(h)
\end{matrix}
\end{equation}
If there exists a maximal chain $\mathfrak C\in \mathcal C(g)\cap \mathcal C(h)$, then by Corollary~\ref{schnittlemma},
$$\mathcal V(g)+\mathcal V_{}(h)=\mathcal V_{\mathfrak C}(g)+\mathcal V_{\mathfrak C}(h)=\mathcal V_{\mathfrak C}(gh)=\mathcal V(gh).$$

On the other hand, if $\mathcal V(gh)=\mathcal V(g)+\mathcal V(h)$ and $\mathfrak C$ is a maximal chain such that $\mathcal V_{\mathfrak C}(gh)=\mathcal V(gh)$,
then \eqref{schnittlemmadiagramm} implies $\mathcal V_{\mathfrak C}(g)+\mathcal V_{\mathfrak C}(h)=\mathcal V_{}(g)+\mathcal V_{}(h)$.
But this is only possible when $\mathcal V_{\mathfrak C}(g)=\mathcal V(g)$ and 
$\mathcal V_{\mathfrak C}(h)=\mathcal V(h)$, and hence $\mathfrak C\in \mathcal C(g)\cap \mathcal C(h)$.
\end{proof}

As an immediate consequence of Proposition~\ref{quasivaluationA} we get:

\begin{coro}\label{polysubalgebra}
Let $p_1,\ldots,p_s\in A$. For $a_1,\ldots,a_s\in\mathbb{N}$, $\mathcal V(f_{p_1}^{a_1}\cdots f_{p_s}^{a_s})=\sum_{i=1}^s a_ie_{p_i}$ if and only if there exists
a maximal chain $\mathfrak C$ containing $p_1,\ldots, p_s$.
\end{coro}

\section{Fan monoids associated to quasi-valuations}\label{fanandcone}

As suggested by Proposition \ref{quasivaluationA}, the image of the quasi-valuation $\mathcal{V}$ in $\mathbb{Q}^A$ is no longer a monoid. Nevertheless, by Corollary \ref{polysubalgebra}, it is not too far away from being a monoid. In the next sections we will study the algebraic and geometric structures of this image. 

\subsection{The fan algebra}\label{Sec:FanAlgbera}

We start with fixing notations. Let
$$\gls{Gamma}:=\{\mathcal V(g)\mid g\in R\setminus\{0\}\}\subseteq \mathbb Q^A$$ 
denote the image of the quasi-valuation in $\mathbb Q^A$: as mentioned, $\Gamma$ is in general not a monoid. Let $\gls{K}$ be the set of all chains in $A$. To every (not necessarily maximal) chain $C\in\mathcal{K}$ we associate the cone $\gls{Kc}$ in $\mathbb R^A$ defined as 
$$
\gls{Kc}=\sum_{p\in C}\mathbb R_{\ge 0}e_p.
$$

The collection $\{K_C\mid C\text{ is a chain in }A\}$, together with the origin $\{0\}$, defines a fan $\mathcal{F}_A$ in $\mathbb{R}^A$. By Lemma \ref{Lem:SS}, the fan $\mathcal{F}_A$ is pure (i.e. all its maximal cones share the same dimension) of dimension $\dim X+1$. Its maximal cones are the cones $K_{\mathfrak C}$ associated to the maximal chains $\mathfrak{C}\in\mathcal{C}$.

For $\underline{a} \in \Gamma$ recall that $\supp \underline{a}:=\{p\in A\mid a_p\not=0\}$. From the definition of $\mathcal{V}$, $\supp \underline{a}\subseteq \mathfrak C$ for some maximal chain $\mathfrak{C}$ and hence by Proposition \ref{positivity}, $\underline{a}\in K_{\mathfrak C}$.
Conversely, for a maximal chain $\mathfrak C$ let $\gls{GammaC}\subseteq \Gamma$ be the subset
$$
\gls{GammaC}:=\{\underline{a}\in\Gamma\mid \supp \underline{a}\subseteq \mathfrak C\}= K_{\mathfrak C}\cap \Gamma.
$$ 
Proposition~\ref{supportcondition}, together with Proposition~\ref{quasivaluationA}, implies:
if $g,h\in R\setminus\{0\}$ are such that $\mathcal V(g)$, $\mathcal V(h)\in\Gamma_\mathfrak{C}$, then 
$$
\mathfrak C\in  \mathcal C(g)\cap  \mathcal C(h)\Longrightarrow \mathcal{V}(g)+\mathcal{V}(h)= \mathcal V(gh)\in \Gamma_{\mathfrak C}.
$$
In other words:

\begin{coro}\label{finiteunionmonoid}
The quasi-valuation image $\Gamma$ is a finite union of monoids: 
$$\gls{Gamma}=\bigcup_{\mathfrak C\in\mathcal{C}} \Gamma_{\mathfrak C}.$$
\end{coro}

\begin{definition}
We call $\gls{Gamma}$ the \emph{fan of monoids} associated to the quasi-valuation $\mathcal V$.
\end{definition}

We associate to the fan of monoids $\Gamma$ the \textit{fan algebra} given in terms of generators and relations. Since the various monoids $\Gamma_{\mathfrak C}$ are submonoids of possibly different lattices, to avoid confusion with the group algebra of one lattice, we prefer to write the elements of $\Gamma$ as lower indexes instead of using them as exponents (upper indexes). 

\begin{definition}\label{Defn:FanAlgebra}
The fan algebra $\gls{KGamma}$ associated to the fan of monoids $\Gamma$ is defined as
$$
\mathbb K[\Gamma]:=\mathbb K[x_{\underline{a}}\mid \underline{a}\in \Gamma] / I(\Gamma)
$$
where $I(\Gamma)$ is the ideal generated by the following elements:
$$
\left\{
\begin{array}{rl}
x_{\underline{a}}\cdot x_{\underline{b}}-x_{\underline{a}+\underline{b}},&\textrm{if there exists a chain $C\subseteq A$ such that }\underline{a},\underline{b}\in K_C;\\
x_{\underline{a}}\cdot x_{\underline{b}},&\textrm{if there exists no such a chain.}\\
\end{array}
\right.
$$
\end{definition}

To simplify the notation, we will write $x_{\underline{a}}$ also for its class in $\mathbb{K}[\Gamma]$.
For a maximal chain $\mathfrak C$ denote by $\mathbb K[\Gamma_{\mathfrak C}]$ the subalgebra:
$$
\mathbb K[\Gamma_{\mathfrak C}]:=\bigoplus_{\underline{a}\in\Gamma_{\mathfrak C}} \mathbb Kx_{\underline{a}}\subseteq \mathbb K[\Gamma],
$$
then $\mathbb K[\Gamma_{\mathfrak C}]$ is naturally isomorphic to the usual semigroup algebra associated to the monoid $\Gamma_{\mathfrak C}$.

We endow the algebra $\mathbb{K}[\Gamma]$ with a grading inspired by Corollary~\ref{degreerelation}: for $\underline{a}\in\mathbb Q^A$, the degree of $x_{\underline{a}}$ is defined by
$$
\deg x_{\underline{a}}=\sum_{p\in A} a_p\deg f_p.
$$

\subsection{Weakly positivity versus standardness}

For a fixed maximal chain $\mathfrak C$, we start with an algebraic relation between the monoid $\Gamma_{\mathfrak C}$ and the core $P_{\mathfrak C}$.

\begin{definition}
A regular function $g\in R\setminus \{0\}$ is called
\begin{itemize} 
\item[{(1)}] \emph{weakly positive} along $\mathfrak C$ if 
$\mathcal V_{\mathfrak C}(g)\in P_{\mathfrak C}$, i.e. $\mathcal V_{\mathfrak C}(g)\in\mathbb Q_{\ge 0}^{\mathfrak C}$. 
\item[{(2)}] \emph{standard} along $\mathfrak C$ if $\mathcal V_{\mathfrak C}(g)\in \Gamma_{\mathfrak C}$, 
 i.e. $\mathcal V(g)=\mathcal V_{\mathfrak C}(g)$.
\end{itemize}
\end{definition}
By Proposition~\ref{positivity}, if $g\not=0$ is standard along $\mathfrak C$, then $g$ is also positive along $\mathfrak C$. So we have a natural inclusion of submonoids of the cone $K_{\mathfrak C}$:
$$
\Gamma_{\mathfrak C}\subseteq P_{\mathfrak C} \subseteq K_{\mathfrak C}.
$$

\begin{rem}
A regular function, which is weakly positive along $\mathfrak{C}$, is not necessarily standard along $\mathfrak{C}$ (see Example \ref{CounterExamplePositive}). One should think of the \textit{weak positivity property} as a local property, whereas the \textit{standardness property} involves all maximal chains, and it is in this sense a global property.
\end{rem}

One might ask up to what extent the two submonoids differ.  

\begin{lemma}\label{coreminimaldifference}
\begin{itemize}
\item[{i)}] The monoid $\Gamma_{\mathfrak C}$ is finitely generated.
\item[{ii)}] The ring extension $\mathbb K[\Gamma_{\mathfrak C}]\subseteq \mathbb K[P_{\mathfrak C}]$ is finite and integral. 
\end{itemize}
\end{lemma}

\begin{proof}
By Corollary~\ref{polysubalgebra},
one has $\mathbb N^{\mathfrak C}\subseteq \Gamma_{\mathfrak C}\subseteq P_{\mathfrak C}$, and hence $\Gamma_{\mathfrak C}$
is a natural $\mathbb N^{\mathfrak C}$-submodule of $P_{\mathfrak C}$. Let $Q$ be as in  
 the proof of Lemma~\ref{finitegenoverN}. We have:
$$
\Gamma_{\mathfrak C}=  \bigcup_{\underline{a}\in Q} \big(\mathbb N^{\mathfrak C}\circ \underline{a}\cap \Gamma_{\mathfrak C}\big)
 \subseteq  \mathbb Q^{\mathfrak C}_{\ge 0}.
$$
The same arguments as in the proof of Lemma~\ref{finitegenoverN}
show that $\mathbb N^{\mathfrak C}\circ \underline{a}\cap \Gamma_{\mathfrak C}$ is finitely generated as $\mathbb N^{\mathfrak C}$-module for all $\underline{a}$ in $Q$.

For each $\underline{a}\in Q$ fix a finite generating system $\mathbb B_{\underline{a}}$ for $\mathbb N^{\mathfrak C}\circ \underline{a}\cap \Gamma_{\mathfrak C}$ as $\mathbb N^{\mathfrak C}$-module. Then the union of the $\mathbb B_{\underline{a}}$, $\underline{a}\in Q$, together with 
$\{e_0,\ldots,e_r\}$, is a finite generating system for the monoid $\Gamma_{\mathfrak C}$.

Since $P_{\mathfrak C}$ is a finitely generated module over $\mathbb N^{\mathfrak C}$, it is hence a
finitely generated module over $\Gamma_{\mathfrak C}$, and thus $\mathbb K[P_{\mathfrak C}]$ is a finite
$\mathbb K[\Gamma_{\mathfrak C}]$-module. So one can find a finite number of elements $\underline{a}^1,\ldots,\underline{a}^s$ in $P_{\mathfrak C}$
such that $\mathbb K[P_{\mathfrak C}]= \mathbb K[\Gamma_{\mathfrak C}][x_{\underline{a}^1},\ldots,x_{\underline{a}^s}]$.
To prove \textit{ii)}, it remains to show that these generators are integral over $\mathbb K[\Gamma_{\mathfrak C}]$.

So suppose $\underline{a}\in P_{\mathfrak C}$, and let $g\not=0$ be a regular function such that 
$\mathcal V_{\mathfrak C}(g)=\underline{a}=\sum_{i=0}^r a_ie_{p_i}\in\mathbb{Q}_{\geq 0}^{\mathfrak{C}}$. 
By Corollary~\ref{powerrelation}, one can find an $m\geq 1$ such that 
$$
\mathcal V_{\mathfrak C}(g^m)=\mathcal V_{\mathfrak C}( f_{p_r}^{ma_r} \cdots f_{p_0}^{ma_0})=m\underline{a}\in 
\mathbb N^{\mathfrak C},
$$
and hence $m\underline{a}\in \Gamma_{\mathfrak C}$.
It follows that $x_{\underline{a}}\in \mathbb K[P_{\mathfrak C}]$ satisfies the equation $p(x^{\underline{a}})=0$ for the monic polynomial $p(y)= y^m-x_{m\underline{a}}\in \mathbb K[\Gamma_{\mathfrak C}][y]$, which finishes the proof.
\end{proof}

\begin{rem}
The proof of Lemma~\ref{coreminimaldifference} does not imply that $g^m$ is standard along 
$\mathfrak C$. But one knows from Corollary~\ref{polysubalgebra} that the function $f_{p_r}^{ma_r} \cdots f_{p_0}^{ma_0}$ is standard along $\mathfrak{C}$.
\end{rem}

Since $\mathbb K[\Gamma_{\mathfrak C}]$ is the algebra associated to a finitely generated submonoid 
of the real cone $K_{\mathfrak C}$, it is a finitely generated integral domain, and the associated variety $\mathrm{Spec}(\mathbb K[\Gamma_{\mathfrak C}])$ is an affine toric variety.
Corollary~\ref{finiteunionmonoid} implies hence that $\mathbb K[\Gamma]$ is a reduced, finitely generated algebra,
so $\mathrm{Spec}(\mathbb K[\Gamma])$ is an affine variety. The geometry of the fan of monoids $\Gamma$ is summarized in the following proposition:

\begin{proposition}\label{Prop:DecompsitionIrrComp}
The affine variety $\mathrm{Spec}(\mathbb K[\Gamma])$ is scheme-theoretically the irredundant union of the toric varieties 
$\mathrm{Spec}(\mathbb K[\Gamma_{\mathfrak C}])$ with $\mathfrak C$ running over the set of maximal chains in $A$. 
Each of these toric varieties is irreducible and of dimension $\dim X+1$.
\end{proposition}

\begin{proof}
Let $I_{\mathfrak C}=\textrm{Ann}(x_{\mathfrak C})$ be the annihilator of the element 
$x_{\mathfrak C}=\prod_{p\in \mathfrak C} x_{e_p}\in\mathbb{K}[\Gamma]$. The multiplication rules in Definition \ref{Defn:FanAlgebra} imply that 
$\mathbb K[\Gamma_{\mathfrak C}]\simeq \mathbb K[\Gamma]/I_{\mathfrak C}$. Since
$\mathbb K[\Gamma_{\mathfrak C}]$ has no zero divisors, it follows that  $I_{\mathfrak C}$ is a prime ideal.

For each maximal chain $\mathfrak C$ the canonical map $\mathbb{K}[\Gamma]\to\mathbb{K}[\Gamma_{\mathfrak C}]$ induces an embedding of a toric variety of dimension $\dim X+1$:
$$\mathrm{Spec}(\mathbb K[\Gamma_{\mathfrak C}])\hookrightarrow \mathrm{Spec}(\mathbb K[\Gamma]).$$

We show that the intersection of the prime ideals $I_\mfC $ is equal to zero: $\bigcap_{\mfC}I_\mfC = (0)$.
Indeed, let $h\in\bigcap_{\mfC}I_\mfC$ be a linear combination of monomials. For a monomial $\prod_{p\in A} x_{n_pe_p}$
let the support be the set $\{p\in A\mid n_p>0\}$. The support of a monomial is, by the definition
of the fan algebra, always contained in a maximal chain. Since $h$ has no constant term, if $h$ is non-zero, there would be at least one non-zero monomial in the linear combination. Let $\mathfrak C$ be a maximal chain containing the support of one of these monomials.
In the product $x_{\mathfrak C}h$, all monomials are supported in $\mathfrak{C}$ and they stay linearly independent. It follows $x_{\mathfrak C}h\not=0$, and thus $h\not\in \bigcap_{\mfC}I_\mfC$, in contradiction to the assumption.

It remains to show the minimality of the intersection. Given a maximal chain $\mfC$ in $A$, we have $ x_\mfC \in I_{\mfC'}$ 
for any $\mfC'\neq \mfC$, while $ x_\mfC\not\in I_\mfC$, which finishes the proof.
\end{proof}

We close the section with some comments on Hilbert quasi-polynomials and a remark on connections to structures similar to the fan algebra. Let $R_1$ be an $\mathbb N$-graded integral $\mathbb K$-algebra 
and let $\mathrm{Quot}(R_1)$ be its quotient field. Denote by $Q(R_1)\subseteq \mathrm{Quot}(R_1)$ the subalgebra generated by the elements 
$\frac{h_1}{h_2}$, $h_1,h_2\in R_1$ homogeneous, it is a $\mathbb Z$-graded algebra.

\begin{lemma}\label{fan_monoid:Hilbert-quasi1}
Suppose that $R_1 \subseteq R_2 \subseteq Q(R_1)$ with $R_2$ a $\mathbb{Z}$-graded $\mathbb K$-algebra, which is finitely generated as an $R_1$--module. If the Hilbert quasi-polynomial of $R_2$ has constant leading coefficient, then the leading term of the Hilbert quasi-polynomial of $R_1$ is the same as that of $R_2$.
\end{lemma}

\begin{proof}
Let $h_1,\ldots, h_q\in R_2$ be the generators of $R_2$ as $R_1$-module. Without loss of generality we may assume that
they are homogeneous and of the form $h_i=\frac{h_{i,1}}{h_{i,2}}$, $h_{i,1},h_{i,2}\in R_1$ homogeneous. Let $h$ be the product
of the denominators. The equality $R_2=\sum_{i=1}^q R_1h_i$ implies $hR_2\subseteq R_1$. The Hilbert quasi-polynomial
of the shifted algebra $hR_2$ has the same leading term as the the one of $R_2$. 
Now the inclusions of graded algebras
$hR_2\subseteq R_1\subseteq R_2$ show that the Hilbert quasi-polynomials of all three have the same leading term.
\end{proof}

\begin{lemma}\label{fan_monoid:Hilbert-quasi2}
Let $\gls{LUC}\subseteq \mathbb Q^{\mathfrak C}$ be the lattice generated by $\Gamma_{\mathfrak C}$
and let
$\tilde \Gamma_{\mathfrak C}=\mathcal L^{\mathfrak C}\cap K_{\mathfrak C}$
be the saturation of the monoid. We endow it with the same grading as $\Gamma_{\mathfrak C}$.
The Hilbert quasi-polynomials associated to $\Gamma_{\mathfrak C}$ and $\tilde \Gamma_{\mathfrak C}$ have the same leading term.
\end{lemma}
\begin{proof}
The quasi-polynomial of $\mathbb K[\tilde \Gamma_{\mathfrak C}]$
is by construction an Ehrhart quasi-polynomial with a constant leading coeffcient. Since 
$\mathbb K[ \tilde \Gamma_{\mathfrak C}]$ is the normalization of $\mathbb K[ \Gamma_{\mathfrak C}]$,
the inclusions $\mathbb K[\Gamma_{\mathfrak C}]\subseteq \mathbb K[\tilde \Gamma_{\mathfrak C}]\subseteq \mathbb K[\mathcal L^{\mathfrak C}]$
fulfill the conditions of Lemma~\ref{fan_monoid:Hilbert-quasi1}.
\end{proof}

\begin{rem}
A structure similar to the fan algebra appears in \cite[Definition 1.1 (3)]{GS} as part of their definition of a toric degeneration of Calabi-Yau varieties. 
We will see later (Theorem~\ref{fanAndDegeneratetheorem}) that the quasi-valuation $\mathcal V$ induces
a filtration on the homogeneous coordinate ring $R$ of $X$, such that the associated graded ring $\text{gr}_{\mathcal V}R$ is isomorphic to the
fan algebra defined above. This leads to a flat degeneration (Theorem~\ref{Thm:Degeneration}) of $X$ into $X_0$, a reduced 
union of equidimensional projective toric varieties. If the Seshadri stratification is normal (see Definition~\ref{def:stratification:normal}) and the partially ordered set $A$ in the Seshadri stratification is shellable, then $X_0$ is Cohen-Macaulay (Theorem~\ref{Thm:ProjNormal}), and the flat degeneration fulfills the conditions in \cite[Definition 1.1(3)]{GS}.

Also the \textit{toric bouquets} associated to a quasifan and a lattice in \cite{AH} lead to a structure similar to the fan algebra. They start with a lattice $M$ and
a quasifan (see \cite{AH} for definitions and details) and associate to this pair the \textit{fan ring}. This ring is defined in a similar way as in Definition~\ref{Defn:FanAlgebra},
see \cite[Definition 7.1]{AH}. The differences between the two approaches are: (1). we do not make the assumption that the semigroups are saturated; (2). we do not need, and hence do not assume to have a lattice $M\subseteq \mathbb Q^A$ such that for all maximal chains we have $\mathcal{L}^{\mathfrak{C}}=M\cap \mathbb Q^{\mathfrak C}$. Recall that $\mathcal{L}^\mathfrak{C}$ is the lattice generated by the semigroup $\Gamma_{\mathfrak C}$ in $\mathbb Q^{\mathfrak C}$. One can find a lattice $M\subseteq \mathbb{Q}^{A}$ such that the $\mathcal{L}^{\mathfrak{C}}$ are of finite index in the intersection $M\cap \mathbb Q^{\mathfrak C}$, but it is an open question to find conditions on the Seshadri stratification ensuring the equality
$\mathcal{L}^{\mathfrak{C}}=M\cap \mathbb Q^{\mathfrak C}$
for all maximal chains.
\end{rem}

\section{Leaves, the associated graded ring and the fan algebra}\label{leavesandgr}

The quasi-valuation $\mathcal{V}:R\setminus\{0\}\to\mathbb{Q}^A$ induces a filtration on the homogeneous coordinate ring $R=\mathbb K[\hat X]$. In this section we further investigate the associated graded algebra $\mathrm{gr}_{\mathcal{V}}R$: we will prove properties parallel to those for valuations in Sections \ref{value_valuations} and \ref{Sec:FiniteGen} and those for fan algebras in Section \ref{fanandcone}.

\begin{definition}
For $\underline{a}\in \Gamma\subseteq \mathbb Q^{A}_{\ge 0}$ we set
$$
\begin{array}{l}
R_{\ge \underline{a}}:=\{g\in R\setminus\{0\}\mid \mathcal V(g)\ge \underline{a}\}\cup\{0\}\ \text{and}\
R_{> \underline{a}}:=\{g\in R\setminus\{0\}\mid \mathcal V(g)> \underline{a}\}\cup\{0\}.
\end{array}
$$
Since the quasi-valuation has only non-negative entries, these subspaces are ideals. Denote the associated graded algebra by: 
$$
\gls{grvR}=\bigoplus_{\underline{a}\in \Gamma} 
R_{\ge \underline{a}}/R_{> \underline{a}}.
$$
Each $R_{\ge \underline{a}}/R_{> \underline{a}}$ for $a\in\Gamma$ will be called a \emph{leaf} of the quasi-valuation.
\end{definition}

We start with the one-dimensional leaves property of this quasi-valuation.

\begin{lemma}\label{NuLeaves}
The leaves $R_{\ge \underline{a}}/R_{> \underline{a}}$,  $\underline{a}\in \Gamma$, are one dimensional.
\end{lemma}

\begin{proof}
Let $f,g$ be non-zero regular functions on $\hat X$ such that $ \mathcal V(g)= \mathcal V(f)=\underline{a}$. Let
$\mathfrak C\in \mathcal C(g)$ and $\mathfrak C'\in \mathcal C(f)$ be such that 
$$
\mathcal V(g) = \mathcal V_{\mathfrak C }(g) =\underline{a} = \mathcal V_{\mathfrak C'}(f)= \mathcal V(f).
$$
These equalities imply $\supp\mathcal V_{\mathfrak C}(g)=\supp\mathcal V_{\mathfrak C'}(f)$ and hence 
$\supp\mathcal V_{\mathfrak C'}(f)\subseteq \mathfrak C$. By Proposition~\ref{supportcondition}, this implies
$\mathfrak C\in \mathcal C(f)$ and hence we can assume $\mathfrak C=\mathfrak C'$. 
By Corollary~\ref{chainleavedimensionone} there exist a scalar $\lambda\in \mathbb K^*$ and $h\in R$ such that $g=\lambda f +h$ with $\mathcal V_{\mathfrak C}(h)>\underline{a}$ when $h\neq 0$.

It remains to show that $\mathcal V(h)>\underline{a}$ if $h\not=0$. Now $\mathcal V$
is a quasi-valuation and hence $\mathcal V(h)\ge \min\{\mathcal V(g),\mathcal V(f)\}=\underline{a}$.
If we have equality, then there exists a maximal chain $\mathfrak C_2$ such that
$\mathcal V_{\mathfrak C_2}(h)=\mathcal V(h)=\underline{a}=\mathcal V_{\mathfrak C}(g)$. A similar argument as above implies
$\supp\mathcal V_{\mathfrak C}(g)=\supp\mathcal V_{\mathfrak C_2}(h)$ and hence
$\supp\mathcal V_{\mathfrak C_2}(h)\subseteq \mathfrak C$. As a consequence one has
$\mathfrak C\in \mathcal C(h)$ and hence
$\mathcal V(h)=\mathcal V_{\mathfrak C}(h)=\underline{a}$, which is a contradiction.
\end{proof}

As a consequence we can prove a version of Corollary~\ref{powerrelation} for the quasi-valuation. The proof relies on the following lemma, which is a generalization of Lemma~\ref{homogeneousparts} to the quasi-valuation $\mathcal V$:

\begin{lemma}\label{homogeneousparts2}
\begin{itemize}
\item[{i)}] For any $g\in R\backslash\{0\}$ and $\lambda\in \mathbb K^*$, $\mathcal{V}(g^\lambda)=\mathcal V(g)$.
\item[{ii)}] Let $g=g_{1}+\ldots +g_{t}\in R=\bigoplus_{i\ge 0} R(i)$ be a decomposition of $g\not=0$ into its homogeneous parts.
Then
$$
\mathcal V(g)=\min\{\mathcal V(g_{j})\mid j=1,\ldots,t\}.
$$
\end{itemize}
\end{lemma}
\begin{proof}
Let $\mathfrak C\in\mathcal C(g)$, then $\mathcal V(g)=\mathcal V_{\mathfrak C}(g)=\mathcal V_{\mathfrak C}(g^\lambda)
\ge \mathcal V(g^\lambda)$. And if $\mathfrak C\in\mathcal C(g^\lambda)$, then we get vice versa:
$\mathcal V(g^\lambda)=\mathcal V_{\mathfrak C}(g^\lambda)=\mathcal V_{\mathfrak C}(g)\ge \mathcal V(g)$,
which proves part i).

To prove ii), recall that $\mathcal V(g)\ge \min\{\mathcal V(g_j)\mid j=1,\ldots,r\}$ by
property (a) in Definition~\ref{quasidef}.
But we can also find pairwise distinct, non-zero scalars $\lambda_1,\ldots,\lambda_t$ such that the linear spans of the following functions coincide:
$
\langle g_1,\ldots,g_t\rangle_{\mathbb K} =\langle g^{\lambda_1},\ldots,g^{\lambda_t}\rangle_{\mathbb K},
$
and hence $\mathcal V(g_j)\ge 
\min\{\mathcal V(g^{\lambda_j})\mid j=1,\ldots,r\}
=\mathcal V(g)$
by part i) of the lemma.
\end{proof}

\begin{coro}\label{powerrelation2}
Let $g\in R\setminus\{0\}$ and suppose $\mathcal V(g)=\sum_{p\in A}a_p e_p$. If $m$ is such that $m a_p\in \mathbb N$ for all $p\in A$, then there exist $\lambda\in\mathbb{K}^*$ and $g'\in R$ such that 
$$
g^m=\lambda \prod_{p\in A} f^{ma_p}_p +g'
$$
with $\mathcal V(g')>\mathcal V(g^m)$ when $g'\neq 0$. If $g$ is homogeneous and $g'\not=0$, 
then $ f_{p_r}^{ma_r} \cdots f_{p_0}^{ma_0}$ and $g'$ are homogeneous of the same degree as $g$.
\end{coro}

Before investigating $\mathrm{gr}_{\mathcal V}R$, we fix a maximal chain $\mathfrak C$ and look at the subspace $\mathrm{gr}_{\mathcal V,\mathfrak C}R$ consisting of leaves supported in $\mathfrak{C}$, i.e.
$$
\gls{grvcR}=\bigoplus_{\underline{a}\in \Gamma_{\mathfrak C}} R_{\ge \underline{a}}/R_{> \underline{a}}
\subseteq \mathrm{gr}_{\mathcal V}R.
$$
This is actually a subalgebra: for $\underline{a},\underline{b}\in \Gamma_{\mathfrak C}$  let  
$g,h\in R$ be representatives of $\bar g\in R_{\ge \underline{a}}/R_{> \underline{a}}\setminus\{0\}$ and 
$\bar h\in R_{\ge \underline{b}}/R_{> \underline{b}}\setminus\{0\}$. Since 
$\mathcal V(g),\mathcal V(h)\in \Gamma_{\mathfrak C}$, we know by Proposition~\ref{supportcondition}: 
$\mathcal V(g)=\mathcal V_{\mathfrak C}(g)$, $\mathcal V(h)=\mathcal V_{\mathfrak C}(h)$, and hence by Lemma~\ref{quasivaluationA},
$\mathcal V(gh)=\mathcal V_{\mathfrak C}(gh)$. This implies $\mathcal V(gh)=\mathcal V(g)+\mathcal V(h)=\underline{a}+\underline{b}$, and therefore
$$
(R_{\ge \underline{a}}/R_{> \underline{a}})\cdot(R_{\ge \underline{b}}/R_{> \underline{b}})
\subseteq 
R_{\ge \underline{a}+\underline{b}}/R_{> \underline{a}+\underline{b}}\subseteq \mathrm{gr}_{\mathcal V,\mathfrak C}R.
$$
It follows that $\mathrm{gr}_{\mathcal V,\mathfrak C}R$ is a subalgebra of $\mathrm{gr}_{\mathcal V}R$ graded by the monoid $\Gamma_\mathfrak{C}$, without zero divisors and
every graded component is one dimensional. This allows us to apply again \cite[Remark 4.13]{BG}:

\begin{lemma}\label{proposition_component}
There exists an isomorphism of algebras $\mathrm{gr}_{\mathcal V,\mathfrak C}R\simeq\mathbb K[\Gamma_{\mathfrak C}]$.
\end{lemma}

By Corollary~\ref{coreminimaldifference} one has the following immediate consequence:

\begin{coro}\label{all_finite-generated}
\begin{enumerate}
    \item[{i)}] The $\mathbb{K}$-algebra $\mathrm{gr}_{\mathcal V}R$ is finitely generated and reduced.
    \item[{ii)}] For any maximal chain $\mathfrak{C}$, the $\mathbb{K}$-algebra $\mathrm{gr}_{\mathcal V,\mathfrak{C}}R$ is a finitely generated integral domain.
\end{enumerate}
\end{coro}

It follows that $\mathrm{Spec}(\mathrm{gr}_{\mathcal V,\mathfrak C}R)$ is an irreducible affine variety,
Lemma~\ref{proposition_component} implies that it is a toric variety. 

Similarly, $\mathrm{Spec}(\mathrm{gr}_{\mathcal V}R)$ is an affine variety. We give a decomposition of $\mathrm{Spec}(\mathrm{gr}_{\mathcal V}R)$ into irreducible components as in Proposition \ref{Prop:DecompsitionIrrComp}. The following proposition can be deduced as a corollary of Theorem \ref{fanAndDegeneratetheorem}. We sketch a direct proof which is similar to Proposition \ref{Prop:DecompsitionIrrComp}.

For a maximal chain $\mfC=(p_r,\ldots,p_0)$ we set $x_\mfC =   f_{p_0} \cdots f_{p_r}\in R$ and $\gls{Ic} := \mathrm{Ann}(\overline{x}_{\mathfrak{C}})\subseteq \mathrm{gr}_{\mathcal V}R$.

\begin{proposition}\label{proposition_quotient_decomposition}
Given a maximal chain $\mfC$, the quotient $\mathrm{gr}_{\mathcal V}R / I_\mfC$ is isomorphic to $\mathrm{gr}_{\mathcal V,\mathfrak C}R$.
In particular, $I_\mfC$ is a homogeneous prime ideal of $\mathrm{gr}_{\mathcal V}R$. Moreover 
$\bigcap_{\mfC}I_\mfC = (0)$ is the minimal prime decomposition of the zero ideal in $\mathrm{gr}_{\mathcal V}R$.
\end{proposition}

\begin{proof}
First note that $\mathcal{V}(x_\mathfrak{C}) = e_{p_r} + e_{p_{r-1}} + \cdots + e_{p_0}$, hence $\supp\mathcal{V}(x_\mathfrak{C}) = \mathfrak{C}$ and $\mathcal{C}(x_\mathfrak{C}) = \{\mathfrak{C}\}$ by Proposition~\ref{supportcondition}. Now, for $g\in R$ we have $\overline{g}\cdot\overline{x_\mathfrak{C}} = 0$ in $\mathrm{gr}_{\mathcal V}R$ if and only if $\mathcal{V}(g\cdot x_\mathfrak{C}) > \mathcal{V}(g) + \mathcal{V}(x_\mathfrak{C})$ and, by Proposition \ref{quasivaluationA}, this is equivalently to 
$\supp\mathcal{V}(g)\not\subseteq\mathfrak{C}$, i.e. $\mathcal{V}(g)\not\in\Gamma_\mathfrak{C}$.

This characterization of the elements of $I_\mfC = \mathrm{Ann}(x_\mathfrak{C})$ proves that
\[
I_\mfC = \bigoplus_{\underline{a}\in\Gamma\setminus\Gamma_\mfC} R_{\geq \underline{a}} / R_{> \underline{a}}.
\]
It follows: $\mathrm{gr}_{\mathcal V} R / I_\mfC \simeq \mathrm{gr}_{\mathcal V,\mfC}R$
and, since  $\mathrm{gr}_{\mathcal V,\mfC}R$ is an integral domain, $I_\mfC$ is a prime ideal. It also follows that $\bigcap_{\mfC}I_\mfC = (0)$.

We want to show that $I_\mathfrak{C}$ is a minimal prime. So, suppose that $I$ is an ideal properly contained in $I_\mathfrak{C}$. The quotient $\mathrm{gr}_{\mathcal V}R / I$ contains an element $\overline{g}$ with $\mathcal{V}(g)\not\in\Gamma_\mathfrak{C}$, in particular we have $\overline{g}\cdot\overline{x_\mathfrak{C}} = 0$ in $\mathrm{gr}_{\mathcal V}R / I$ again by Proposition \ref{quasivaluationA}. Hence this quotient is not a domain and $I$ is not prime.

Finally, note that, for a maximal chain $\mathfrak{D}$, $\overline{x_\mathfrak{D}}$ is a non-zero element in the intersection $\bigcap_{\mfC\neq\mathfrak{D}}I_\mfC$. This shows that the intersection $\bigcap_{\mfC}I_\mfC$ is non-redundant.
\end{proof}

\begin{coro}\label{irredUnion}
The variety $\mathrm{Spec}(\mathrm{gr}_{\mathcal V}R)$ is scheme-theoretically the irredundant union of the irreducible varieties 
$\mathrm{Spec}(\mathrm{gr}_{\mathcal V,\mathfrak C}R)$ with $\mfC$ running over the set of maximal chains of $A$; 
each of these varieties is irreducible and of dimension $\dim X+1$.
\end{coro}
\begin{proof} 
This follows by Lemma~\ref{proposition_component} and Proposition~\ref{proposition_quotient_decomposition}.
\end{proof}

\section{The fan algebra and the degenerate algebra}\label{Sec:FanDeg}

The goal of this section is to give a description of the associated graded algebra $\mathrm{gr}_{\mathcal{V}}R$ (Theorem \ref{fanAndDegeneratetheorem}). It is a wide generalization of the conjecture on special LS-algebra stated in \cite[Remark 1]{Chi}.

We fix throughout this section a vector space basis $\mathbb B$ of the degenerate algebra $\text{gr}_{\mathcal V}R$: 
$$
\mathbb B=\{\bar g_{\underline{a}}\mid \underline{a}\in\Gamma, 0\neq\bar{g}_{\underline{a}}\in R_{\ge \underline{a}}/R_{> \underline{a}} \}.
$$
By Proposition~\ref{quasivaluationA} we know for $\underline{a},\underline{b}\in\Gamma$: 
$\bar g_{\underline{a}}\cdot \bar g_{\underline{b}}\not=0$ in $\text{gr}_{\mathcal V}R$ if and only if
one can find a chain $C$ in $A$ (not necessarily maximal) such that $\supp \underline{a},\supp\underline{b}\subseteq C$. If this holds,
then one finds some non-zero $c_{\underline{a},\underline{b}}\in\mathbb K^*$ such that 
$\bar g_{\underline{a}}\cdot \bar g_{\underline{b}}=c_{\underline{a},\underline{b}} \bar g_{\underline{a}+\underline{b}}$. 
This list of non-zero coefficients $c_{\underline{a},\underline{b}}$ provides
a complete description of $\text{gr}_{\mathcal V}R$ in terms of generators and relations.

\begin{theorem}\label{fanAndDegeneratetheorem} 
The  degenerate algebra $\mathrm{gr}_{\mathcal V}R$ is isomorphic to the fan algebra $\mathbb K[\Gamma]$.
\end{theorem}
\begin{proof}
First note that if we fix constants $c_{\underline{a}}\in\mathbb{K}^*$ for all $\underline{a}\in\Gamma$, the linear map 
$$\mathbb K[\Gamma]\rightarrow \mathrm{gr}_{\mathcal V}R,\ \  x_{\underline{a}}\mapsto c_{\underline{a}}\bar g_{\underline{a}}$$
is a vector space isomorphism between the fan algebra $\mathbb{K}[\Gamma]$ and the degenerate algebra $\mathrm{gr}_{\mathcal{V}}R$.

Suppose we have already had for all $\underline{a}\in\Gamma$ non-zero elements $c_{\underline{a}}\in \mathbb K^*$ with the 
following property: whenever $\supp \underline{a},\supp\underline{b}\subseteq C$ for some chain $C$ in $A$, then 
\begin{equation}\label{triplenumbersone}
c_{\underline{a}}\cdot c_{\underline{b}}=c_{\underline{a}+\underline{b}}\cdot c_{\underline{a},\underline{b}}.
\end{equation}
We rescale the basis $\mathbb B$ and get a new basis: $\mathbb B'=\{\bar h_{\underline{a}}=\frac{1}{c_{\underline{a}}}\bar g_{\underline{a}}\mid 
\underline{a}\in\Gamma\}$. The corresponding rescaled vector space isomorphism is defined on the new basis as follows:
$$
\chi:\mathbb K[\Gamma]\rightarrow \text{gr}_{\mathcal V}R, \quad x_{\underline{a}}\mapsto \bar h_{\underline{a}}\quad\textrm{\ for all $\underline{a}\in\Gamma$}.
$$
This is in fact an algebra isomorphism: if $\underline{a},\underline{b}\in\Gamma$ are such that there is no chain in $A$ containing both $\supp \underline{a}$ and $\supp\underline{b}$, then $\bar h_{\underline{a}}\cdot \bar h_{\underline{b}}=0$, and otherwise we get
$$
\bar h_{\underline{a}}\cdot \bar h_{\underline{b}}=\frac{1}{c_{\underline{a}}c_{\underline{b}}}\bar g_{\underline{a}}\bar g_{\underline{b}}=
\frac{c_{\underline{a},\underline{b}}}{c_{\underline{a}}c_{\underline{b}}}\bar g_{\underline{a}+\underline{b}}=\bar h_{\underline{a}+\underline{b}}.
$$
It remains to prove the existence of the $c_{\underline{a}}\in\mathbb K^*$ for $\underline{a}\in\Gamma$, satisfying the condition
in \eqref{triplenumbersone}. This will be done in the next subsections.
\end{proof}
\subsection{Existence of the rescaling  coefficients}\label{existence}
To prove the existence of the rescaling coefficients $c_{\underline{a}}$, $\underline{a}\in\Gamma$, we construct an affine variety $Z$
having the following properties: 
\begin{enumerate}\label{varietycondition}
\item\label{varietycondition1} There is a natural interpretation of the $\bar g_{\underline{a}}$, $\underline{a}\in\Gamma$, as functions on $Z$.
\item\label{varietycondition2} The relations $\bar g_{\underline{a}}\bar g_{\underline{b}}=c_{\underline{a},\underline{b}}\bar g_{\underline{a}+\underline{b}}$
for $\underline{a},\underline{b}\in \Gamma$ hold also in $\mathbb K[Z]$ whenever there exists a chain $C$ in $A$ such that 
$\supp \underline{a},\supp\underline{b}\subseteq C$.
\item\label{varietycondition3}  There exists a point $z\in Z$ such that $\bar g_{\underline{a}}(z)\not=0$ for all $\underline{a}\in\Gamma$.
\end{enumerate}
\vskip 3pt
\noindent
\textbf{Proof of the existence of the rescaling coefficients.} Suppose $Z$ is an affine variety having the above properties.
Let $z\in Z$ be a point as in \ref{varietycondition}.\ref{varietycondition3}. By  \ref{varietycondition}.\ref{varietycondition1} it is possible to define 
$c_{\underline{a}}=\bar g_{\underline{a}}(z)$ for all  $\underline{a}\in\Gamma$, and these are all elements in $\mathbb K^*$.
Then property \ref{varietycondition}.\ref{varietycondition2} implies:  if there exists a chain $C$ in $A$ such that $\supp \underline{a},\supp\underline{b}\subseteq C$, then
$$
c_{\underline{a}}\cdot c_{\underline{b}}=g_{\underline{a}}(z)\cdot g_{\underline{b}}(z)=c_{\underline{a},\underline{b}}g_{\underline{a}+\underline{b}}(z)
=c_{\underline{a},\underline{b}}\cdot c_{\underline{a}+\underline{b}}
$$
So this collection of non-zero rescaling coefficients has the desired property \eqref{triplenumbersone}.
\hfill$\bullet$
\subsection{The varieties $Z_{\mathcal M}$}
It remains to construct a variety $Z$ with the properties \ref{varietycondition}.\ref{varietycondition1} -- \ref{varietycondition}.\ref{varietycondition3}. This will be done using an inductive procedure.
Let $C\subseteq A$ be a chain, not necessarily maximal. We associate to $C$
the submonoid $\Gamma_C=\{\underline{a}\in\Gamma\mid \supp(\underline{a})\subseteq C\}$.
Denote by $M_C=\langle \Gamma_C\rangle_{\mathbb Z}\subseteq \mathbb Q^A$ the lattice generated by 
$\Gamma_C$ in $\mathbb Q^A$, and let $N_C$ be the dual lattice of $M_C$. Set
$$
Q_C=\bigoplus_{\underline{a}\in\Gamma_C} R_{\ge \underline{a}}/R_{>\underline{a}}\subseteq \mathrm{gr}_{\mathcal V}R.
$$ 
In Section~\ref{leavesandgr} we have discussed the case where $C=\mathfrak C$ is a maximal chain. The same arguments can be applied to show that $Q_C$ is a finitely generated integral domain. Denote  $Y_C=\mathrm{Spec\,}Q_C$ the associated affine variety. Again,
as in Section~\ref{leavesandgr}, $Y_C$ is a toric variety for the torus $T_C:=T_{N_C}$, i.e. the torus associated to the 
lattices $N_C$. The $\mathbb K$-algebra $Q_C$ is positively graded, so the affine variety $Y_C$ has a unique vertex which we denote by $0$.

The set of all chains in $A$ is partially ordered with respect to the inclusion relation. 

We fix an inclusion of chains
$C_1\subseteq C_2$. It induces an inclusion of monoids $i_{C_1,C_2}:\Gamma_{C_1}\hookrightarrow  \Gamma_{C_2}$,
which in turn induces a morphism of algebras: 
$i_{C_1,C_2}:Q_{C_1}\hookrightarrow  Q_{C_2}$. By \cite{CLS}, Proposition~1.3.14, we get an induced toric morphism
between the associated varieties $\psi_{C_1,C_2}:Y_{C_2}\rightarrow  Y_{C_1}$ and a group homomorphism
$\phi_{C_1,C_2}:T_{C_2}\rightarrow  T_{C_1}$. All these maps are compatible: for a sequence of inclusions of chains
$C_1\hookrightarrow C_2\hookrightarrow C_3$, one has $i_{C_2,C_3}\circ i_{C_1,C_2}=i_{C_1,C_3}$ and 
$\psi_{C_1,C_2} \circ \psi_{C_2,C_3} =\psi_{C_1,C_3}$ and so on.

Moreover, we also have a closed immersion $\Psi_{C_1,C_2}:Y_{C_1}\hookrightarrow  Y_{C_2}$ of varieties: the subspace
$$
I_{C_1,C_2}=\bigoplus_{\substack{\underline{a}\in\Gamma_{C_2}\setminus\Gamma_{C_1}}} R_{\ge \underline{a}}/R_{>\underline{a}}\subseteq Q_{C_2}
$$ 
is an ideal and we have a natural isomorphism of algebras $Q_{C_2}/I_{C_1,C_2}\rightarrow Q_{C_1}$. This isomorphism induces
the desired closed immersion $\Psi_{C_1,C_2}:Y_{C_1}=\mathrm{Spec\,}Q_{C_1}\hookrightarrow  Y_{C_2}=\mathrm{Spec\,}Q_{C_2}$. 
The composition of algebra morphisms
$$
Q_{C_1}\hookrightarrow  Q_{C_2}\rightarrow Q_{C_1}\simeq Q_{C_2}/I_{C_1,C_2}
$$ 
is the identity map on $Q_{C_1}$. The composition of morphisms of affine varieties:
\[\begin{tikzcd}        
    Y_{C_1} \arrow[hook]{rr}{\Psi_{C_1,C_2}} & & Y_{C_2}  \arrow{rr}{\psi_{C_1,C_2}} & & Y_{C_1}
\end{tikzcd}\]
is hence the identity map on $Y_{C_1}$.

\begin{rem}
We can extend the above definitions to the case when $C=\emptyset$ is the empty set. In this case $\Gamma_\emptyset=\{0\}\in\Gamma$, $M_\emptyset=\{0\}\in\mathbb{Q}^A$, $Q_\emptyset=\mathbb{K}$ and $Y_{\emptyset}$ is the vertex $0$ in all toric varieties $Y_C$. Such  definitions are compatible with the inclusion $\emptyset\subseteq C$ for any chain $C$; this allows us to define $i_{\emptyset,C}$, $\psi_{\emptyset,C}$, $\Psi_{\emptyset,C}$, etc.
\end{rem}

\begin{definition}
A subset $\mathcal{M}$ of the set of all chains in $A$ is called \emph{saturated} if for any chain $C\in\mathcal{M}$ and a subset $C'\subseteq C$, we have $C'\in\mathcal{M}$. The set of all saturated subsets in $A$ will be denoted by $\gls{Ks}$.
\end{definition}

Notice that the empty set is contained in any saturated subset. The set $\mathcal K^s$ is partially ordered with respect to inclusion.

We now associate to  $\mathcal M\in \mathcal K^s$ an affine variety $Z_{\mathcal M}$ as follows: 

\begin{definition}
For $\mathcal M\in \mathcal K^s$, we define $Z_{\mathcal M}$ as the following closed subset of the product of varieties $\prod_{C\in\mathcal M} Y_C$:
$$
\gls{ZM}:=\left\{(y_C)_{C}\in \prod_{C\in\mathcal M} Y_C\mid \forall\, C_1\subseteq C_2\in\mathcal{M},\ \psi_{C_1,C_2}(y_{C_2})=y_{C_1}\right\}.
$$
We endow $Z_{\mathcal M}$ with the induced reduced structure as an affine variety.
\end{definition}

For $C\in \mathcal M$, the canonical projection $p_C:\prod_{C'\in\mathcal M} Y_{C'}\rightarrow Y_C$ restricts to a morphism $p_C:Z_{\mathcal M}\rightarrow Y_C$.

\begin{lemma}\label{Lem:PCSurj}
The morphism $p_C:Z_{\mathcal M}\rightarrow Y_C$ is surjective.
\end{lemma}
\begin{proof}
To prove the lemma, we associate to a chain $C\in \mathcal M$ a subvariety $Z_C\subseteq Z_{\mathcal M}$ having the property: 
$p_C\vert_{Z_C}:Z_C\rightarrow Y_C$ is an isomorphism.
As a first step we define $Z_C$ as a subvariety of $\prod_{C'\in \mathcal M} Y_{C'}$ and set:
$$
Z_C=\left\{\left.(y_{C'})_{C'}\in \prod_{C'\in \mathcal M} Y_{C'}\right\vert 
\begin{array}{ll}
y_C\in Y_C&\ \\ 
y_{C'}=\Psi_{C\cap C',C'}\circ\psi_{C\cap C',C}(y_{C})&\forall C'\in\mathcal{M}
\end{array}
 \right\}.
$$
The subset $Z$ is closed, we endow $Z_C$ with the induced reduced structure. 

We show $Z_C\subseteq Z_{\mathcal M}$.
Given a point $(y_{C'})_{C'}\in Z_C$ and two elements $C_1\subseteq C_2$ in $\mathcal{M}$, we have to verify that $\psi_{C_1,C_2}(y_{C_2})=
y_{C_1}$. Indeed, set $C_i'=C\cap C_i$ for $i=1,2$. The inclusion $C_1\subseteq C_2$ induces an inclusion $C_1'\subseteq C_2'$. The two inclusions
induce two algebra homomorphisms: $i_{C_1,C_2}:Q_{C_1}\hookrightarrow  Q_{C_2}$ followed by the quotient 
$q_{C_2,C_2'}:Q_{C_2}\rightarrow Q_{C_2'}\simeq Q_{C_2}/I_{C_2',C_2}$, and the quotient map:
$q_{C_1,C_1'}:Q_{C_1}\rightarrow Q_{C_1'}\simeq Q_{C_1}/I_{C_1',C_1}$ followed by the 
monomorphism $i_{C'_1,C'_2}:Q_{C'_1}\hookrightarrow  Q_{C'_2}$. By construction, these two algebra homomorphisms are the same,
which in terms of morphisms of varieties implies:
$$
\psi_{C_1,C_2}\circ\Psi_{C_2',C_2}=\Psi_{C_1',C_1}\circ\psi_{C_1',C_2'}.
$$
And hence we conclude:
$$
\begin{array}{rll}
\psi_{C_1,C_2}(y_{C_2})=\psi_{C_1,C_2}\circ(\Psi_{C_2',C_2}\circ\psi_{C_2',C}(y_C))
&=&\Psi_{C_1',C_1}\circ\psi_{C_1',C_2'}\circ\psi_{C_2',C}(y_C)\\
&=&\Psi_{C_1',C_1}\circ\psi_{C_1',C}(y_C)\\
&=&y_{C_1},
\end{array}
$$
which shows: $Z_C\subseteq Z_{\mathcal M}$.

The restriction of $p_C$ to $Z_C$ induces a bijection between $Z_C$ and $Y_C$.
The above description of $Z_C$ can be used to define an inverse map. It follows that $Z_C$ is isomorphic to $Y_C$.
\end{proof}

The varieties $Z_{\mathcal M}$ satisfy the usual universal property of fibre products in the category of affine varieties:

\begin{lemma}\label{universprop}
Let $V$ be an affine variety. Given morphisms $\Phi_C:V\rightarrow Y_C$ for all $C\in \mathcal M$ such that for all $C_1\subseteq C_2$ in $\mathcal M$, $\Phi_{C_1}= \psi_{C_1,C_2}\circ \Phi_{C_2}$ holds. Then there exist a unique morphism $\eta:V\rightarrow Z_{\mathcal M}$ such that $p_C\circ\eta=\Phi_C$. 
\end{lemma}
\begin{proof}
The image of the obvious map $\prod \Phi_C:V\rightarrow \prod_{C\in\mathcal M}Y_C$ 
is contained in  $Z_{\mathcal M}$ and this map has the desired property.
\end{proof}

When writing down a saturated subset, we usually omit the empty set.

\begin{example}\label{lengthzero}
If $\mathcal M=\{C\}$ is just a chain of length zero, then $Z_{\mathcal M}=Y_C$ is a one-dimensional (not necessarily normal) toric
variety. If $\mathcal M=\{C_1,\ldots,C_s\}$ is a collection of chains, all of length zero, then the fibre product becomes a direct product
and $Z_{\mathcal M}= Y_{C_1}\times\cdots\times Y_{C_s}$ is an $s$-dimensional toric variety.
\end{example}

\begin{lemma}\label{uniquemaxexample}
If $\mathcal M$ has a unique maximal element $C$, then $Z_{\mathcal M} \simeq Y_C$. 
\end{lemma}

\begin{proof}
By Lemma~\ref{universprop}, the morphisms $\psi_{C',C}$ for $C'\subseteq C$ induce a morphism $\phi:Y_C\rightarrow \prod_{C'\in \mathcal M}Y_{C'}$ with image in $Z_{\mathcal M}$. This image is in fact the entire $Z_{\mathcal M}$ by construction. The compositions $p_C\circ\phi$ and $\phi\circ p_C$, where $p_C:Z_{\mathcal{M}}\to Y_C$ is the surjective map in Lemma \ref{Lem:PCSurj}, are identity maps.
\end{proof}

Let $\mathcal M',\mathcal M\in \mathcal K^s$ be such that $\mathcal M'\subseteq \mathcal M$. The projection
$\prod_{C''\in \mathcal M}Y_{C''}\rightarrow \prod_{C''\in \mathcal M'}Y_{C''}$ induces a morphism between the fibered products:

\begin{lemma}\label{Lem:InduceMorphism}
Any inclusion $\mathcal M'\subseteq\mathcal M$ induces a morphism $\psi_{\mathcal M',\mathcal M}:Z_{\mathcal M}\rightarrow Z_{\mathcal M'}$.
These morphisms are compatible with the inclusion relations, i.e. for  $\mathcal M''\subseteq\mathcal M' \subseteq\mathcal M$ one has
$\psi_{\mathcal M'',\mathcal M}=\psi_{\mathcal M'',\mathcal M'}\circ \psi_{\mathcal M',\mathcal M}$.
\end{lemma}

These morphisms can be used to describe an inductive procedure to construct $Z_{\mathcal M}$. Let $C\in \mathcal M$ be a chain which is 
maximal in $\mathcal M$ with respect to the inclusion relation. Let $\mathcal M_C=\{C'\in\mathcal M\mid C'\subseteq C\}$, $\mathcal M'=\mathcal M\setminus\{C\}$ and $\mathcal M_C'=\mathcal M_C\setminus\{C\}$. These sets are all saturated. Recall that the fibre products are considered in the category of varieties. 

\begin{lemma}\label{fibreproduct}
There exists an isomorphism of varieties $Z_{\mathcal M}\simeq Y_C\times_{Z_{\mathcal M_C'}}Z_{\mathcal M'}$.
\end{lemma}

\begin{proof}
If $C$ is the only chain in $\mathcal M$ which is maximal with respect to the inclusion relation, then $\mathcal M=\mathcal M_C$
and $\mathcal M'=\mathcal M'_C$, and hence $Y_C\times_{Z_{\mathcal M'}}Z_{\mathcal M'}\simeq Y_C\simeq Z_{\mathcal M}$
by Lemma~\ref{uniquemaxexample}. So without loss of generality we may assume in the following
that $C$ is not the only chain in $\mathcal M$ which is maximal with respect to the inclusion relation.

By the universal property of the fibre product, the commutative square on the right hand side, induced by the commutative square on the left hand side (Lemma \ref{Lem:InduceMorphism}), gives rise to a morphism $Z_{\mathcal M}\rightarrow Y_C\times_{Z_{\mathcal M_C'}}Z_{\mathcal M'}$.
\[
\xymatrix{
\mathcal{M}_C'\ar[r] \ar[d] & \mathcal{M}'\ar[d]\\
\mathcal{M}_C\ar[r] & \mathcal{M}
}
\hskip 20pt
\xymatrix{
Z_{\mathcal{M}}\ar[r] \ar[d] & Z_{\mathcal{M}_C}\cong Y_C \ar[d] \\
Z_{\mathcal{M}'}\ar[r] & Z_{\mathcal{M}_C'}
}
\]
By construction, the fibre product is a subvariety of
$$
Y_C\times \prod_{C'\in \mathcal M_C'} Y_{C'}\times  \prod_{C'\in \mathcal M_C'}  Y_{C'} \times \prod_{C''\in \mathcal M'\setminus \mathcal M_C'} Y_{C''}.
$$
where, in addition to the maps defining $Z_{\mathcal M'}\subseteq  \prod_{C'\in \mathcal M_C'}  Y_{C'} \times \prod_{C''\in \mathcal M'\setminus \mathcal M_C'} Y_{C''}$
(here we take the second copy of $\prod_{C'\in \mathcal M_C'}  Y_{C'}$) and $Z_{\mathcal M_C'}\subseteq \prod_{C'\in \mathcal M_C'}  Y_{C'}$ (here we take the first copy of $\prod_{C'\in \mathcal M_C'}  Y_{C'}$) and so on, we have the identity map between the two copies of the product $\prod_{C'\in \mathcal M_C'}  Y_{C'}$. 
So we may omit one copy of the product $\prod_{C'\in \mathcal M_C'}  Y_{C'}$, what is left is the variety $Z_{\mathcal M}$.
\end{proof}

For $\mathcal M\in\mathcal K^s$ let $\mathcal M_0\subseteq \mathcal M$ be the subset of all chains of length zero.

\begin{proposition}\label{finitemap}
The morphism $\psi_{\mathcal M_0,\mathcal M}: Z_{\mathcal M}\rightarrow Z_{\mathcal M_0}$ is finite.
\end{proposition}

Before going to the proof, we discuss some consequences of the proposition.

\begin{coro}\label{dimensionforZM}
We have: $\dim Z_{\mathcal M}=\#\mathcal M_0$.
\end{coro}
If we have saturated sets such that $\mathcal M_0\subseteq \mathcal M'\subseteq \mathcal M$, then, by Proposition~\ref{finitemap},
the morphisms $\psi_{\mathcal M_0,\mathcal M}:Z_{\mathcal M}\rightarrow Z_{\mathcal M_0}$ and 
$\psi_{\mathcal M_0,\mathcal M'}:Z_{\mathcal M'}\rightarrow Z_{\mathcal M_0}$ are finite. Since $\psi_{\mathcal M_0,\mathcal M}$
is the composition of $\psi_{\mathcal M',\mathcal M}$ and $\psi_{\mathcal M_0,\mathcal M'}$, this implies:

\begin{coro}
Let $\mathcal M'\subseteq \mathcal M$ be such that $\mathcal M_0\subseteq \mathcal M'$. Then 
$\psi_{\mathcal M_0,\mathcal M}: Z_{\mathcal M}\rightarrow Z_{\mathcal M'}$ is a finite morphism.
\end{coro}

More generally, for $\mathcal M'\subseteq \mathcal M$ we set $\mathcal M''=\mathcal M'\cup\mathcal M_0$.
The morphism $\psi_{\mathcal M',\mathcal M}$ is the composition of the finite morphism $\psi_{\mathcal M'',\mathcal M}$ and 
the projection $\psi_{\mathcal M',\mathcal M''}$, and hence:

\begin{coro}
If $\mathcal M'\subseteq \mathcal M$, then 
$\psi_{\mathcal M',\mathcal M}: Z_{\mathcal M}\rightarrow Z_{\mathcal M'}$ is a surjective morphism.
\end{coro}

\begin{proof}[Proof of Proposition~\ref{finitemap}] 
The proof is by induction on the number of elements in $\mathcal M$. If $\# \mathcal M=1$, then $\mathcal M=\mathcal M_0$,
and there is nothing to prove. Suppose now $\#\mathcal M >1$ and the claim holds for all saturated sets $\mathcal M'$ having strictly fewer elements than $\mathcal M$. 
Let  $C\in \mathcal M$ be a chain which is maximal in $\mathcal M$ with respect to the inclusion relation. The set $\mathcal M'=\mathcal M\setminus\{C\}$
is a saturated set. There are three possible cases:

(a). $C$ is a chain of length zero. In this case the fibre product in Lemma~\ref{fibreproduct} becomes a cartesian product, i.e.
$Z_{\mathcal M}\simeq Y_C\times Z_{\mathcal M'}$ and $Z_{\mathcal M_0}\simeq Y_C\times Z_{\mathcal M'_0}$. In this case
the morphism $Z_{\mathcal M}\rightarrow Z_{\mathcal M_0}$ is finite by induction. 

(b). $C$ is the only chain in $\mathcal M$ which is maximal with respect to the inclusion relation (hence of length $\geq 1$). In this case we know $\mathcal M_0=\mathcal M'_0$ and by Lemma \ref{uniquemaxexample}, $Z_{\mathcal M}\simeq Y_C$ is an irreducible toric variety.

Chains of length zero in $\mathcal{M}$ are just elements of $A$. For $p\in A$, we denote by $C_p$ the chain in $\mathcal{M}$ consisting of $p$, i.e. $C_p=\{p\}$.

For every chain $C_p\in\mathcal{M}_0$, the algebra $Q_{C_p}=\bigoplus_{\underline{a}\in\Gamma_{C_p}} R_{\ge \underline{a}}/R_{>\underline{a}}\simeq \mathbb K[Y_{C_p}]$ is a subalgebra of $\mathbb K[Y_{C}]$. As a matter of fact, from the construction, the composition of the morphisms $Y_{C}=Z_{\mathcal M}\rightarrow Z_{\mathcal M_0}$ and $Z_{\mathcal M_0}\rightarrow Y_{C_p}$ identifies $\mathbb K[Y_{C_p}]$ with this subalgebra. It follows that the morphism $\psi_{\mathcal M_0,\mathcal M}^*$ identifies
$\mathbb K[Z_{\mathcal M_0}]=\bigotimes_{C_p\in \mathcal M_0}\mathbb K[Y_{C_p}]$ with the subalgebra of $\mathbb K[Y_{C}]$ generated by the $Q_{C_p}$, $C_p\in \mathcal M_0$. In particular, the morphism $\psi_{\mathcal M_0,\mathcal M}^*$ is injective and $\psi_{\mathcal M_0,\mathcal M}$ is hence a dominant morphism.

By adding some elements to a (finite) generating system of $\mathbb K[Y_C]$
if necessary, we may without loss of generality assume it to have the form 
$$
\{\bar f_p\mid C_p\in \mathcal M_0\}\cup\{\bar g_{\underline{a}^1},\ldots,\bar g_{\underline{a}^t}\}.
$$
By Corollary~\ref{powerrelation2}, each of the generators
$\bar g_{\underline{a}^j}$ is integral over the polynomial subring generated by the $\{\bar f_p\mid C_p\in \mathcal M_0\}$. It follows that
$\mathbb K[Z_{\mathcal M}]=\mathbb K[Y_C]$ is a finite module over the subalgebra generated by the $\bar f_p$,  $C_p\in \mathcal M_0$,
and hence $\mathbb K[Z_{\mathcal M}]$ is a finite module over $\mathbb K[Z_{\mathcal M_0}]$.

(c). $C$ is a chain of length $\ge 1$ and $C$ is not the only chain in $\mathcal M$ which is maximal with respect to the inclusion relation. 
We still have $\mathcal M_0=\mathcal M'_0$.
Set $\mathcal M_C =\{ C'\in \mathcal M\mid C'\subseteq C\}$ and $\mathcal M'_C=\mathcal M_C\setminus\{C\}$. These sets
are saturated, and $Z_{\mathcal M}\simeq Y_C\times_{Z_{\mathcal M'_C}}Z_{\mathcal M'}$
by Lemma~\ref{fibreproduct}.
In this case one can use base change arguments. By induction, $Y_C=Z_{\mathcal M_C}$ is finite over $Z_{\mathcal M_{C,0}}$ and so is 
$Z_{\mathcal M'_C}$ over $Z_{\mathcal M'_{C,0}}=Z_{\mathcal M_{C,0}}$. 
Since $\psi_{\mathcal M_{C,0},\mathcal M_C}=\psi_{\mathcal M_{C,0},\mathcal M'_{C}}\circ
\psi_{\mathcal M'_C,\mathcal M_{C}}$, $Z_{\mathcal M_C}$ is finite over $Z_{\mathcal M'_{C}}$.
By base change, this implies $Z_{\mathcal M}\simeq Y_C\times_{Z_{\mathcal M'_C}}Z_{\mathcal M'}$ is finite over 
$Z_{\mathcal M'}$. By induction, $Z_{\mathcal M'}$ is finite over $Z_{\mathcal M_0}=Z_{\mathcal M'_0}$. Since
$\psi_{\mathcal M_{0},\mathcal M}$ is the composition of $\psi_{\mathcal M',\mathcal M}$ and $\psi_{\mathcal M_{0},\mathcal M'}$, $Z_{\mathcal M}$ is finite over $Z_{\mathcal M_0}$. 
\end{proof}

\subsection{The variety $Z_{\mathcal{K}^s}$}
The candidate for the variety $Z$ mentioned in subsection~\ref{existence} is the variety $Z_{\mathcal K^s}$.

\begin{proposition}
The affine variety $Z_{\mathcal K^s}$ satisfies the three properties \ref{varietycondition}.\ref{varietycondition1}--\ref{varietycondition}.\ref{varietycondition3}. 
\end{proposition}

\begin{proof}
Given $\underline{a}\in \Gamma$, let $C,C_1,C_2$ be chains such that $C=\supp \underline{a}$  
and $C\subseteq C_1,C_2$. Denote by $\mathcal M_C$ (resp. $\mathcal M_{C_1}$, $\mathcal M_{C_2}$)
the smallest saturated sets in $\mathcal{K}^s$ containing $C$ (resp. $C_1$, $C_2$). 
By Lemma~\ref{uniquemaxexample}, we have $Z_{\mathcal M_C}\simeq Y_C$ and  $Z_{\mathcal M_{C_i}}\simeq Y_{C_i}$
for $i=1,2$.
The morphisms $\psi_{\mathcal M_C,\mathcal M_{C_i}}$
are just the same as the morphisms $\psi_{C,{C_i}}$, $i=1,2$, and hence $\psi^*_{\mathcal M_C,\mathcal M_{C_i}}$
is just the inclusion $i_{C,C_i}: Q_C\hookrightarrow Q_{C_i}$, $i=1,2$.
 
By construction, $\bar g_{\underline{a}}\in \mathbb K[Y_C]= \mathbb K[Z_{\mathcal M_C}]=Q_C$, and it is also
an element in $\mathbb K[Y_{C_1}]$ and $\mathbb K[Y_{C_2}]$. By the compatibility of the dominant morphisms:
$$
\psi_{M_C,\mathcal K^s}=\psi_{M_C,\mathcal M_{C_1}}\circ\psi_{\mathcal M_{C_1},\mathcal K^s}=
\psi_{M_C,\mathcal M_{C_2}}\circ\psi_{\mathcal M_{C_2},\mathcal K^s}
$$
we see that $\psi^*_{\mathcal M_C,\mathcal K^s}(\bar g_{\underline{a}})
=\psi^*_{\mathcal M_{C_1},\mathcal K^s} \circ\psi^*_{\mathcal M_C,\mathcal M_{C_1}}(\bar g_{\underline{a}})
=\psi^*_{\mathcal M_{C_2},\mathcal K^s} \circ\psi^*_{\mathcal M_C,\mathcal M_{C_2}}(\bar g_{\underline{a}})$.

(1). The variety $Z_{\mathcal K^s}$ satisfies the property \ref{varietycondition}.\ref{varietycondition1}: given $\underline{a}\in\Gamma$, we can view $\bar g_{\underline{a}}$
as a function on $Z_{\mathcal K^s}$: just take any chain $C$ in $A$ such that $\supp \underline{a}\subseteq C$, then
$\psi^*_{\mathcal M_{C},\mathcal K^s}(\bar g_{\underline{a}})\in\mathbb K[Z_{\mathcal K^s}]$ is well defined and 
independent of the choice of $C$.

(2).  The variety $Z_{\mathcal K^s}$ satisfies the property \ref{varietycondition}.\ref{varietycondition2}: if there exists a chain $C$ in $A$ such that 
$\supp \underline{a},\supp\underline{b}\subseteq C$, then we have 
$$
\psi^*_{\mathcal M_{C},\mathcal K^s}(\bar g_{\underline{a}})\cdot \psi^*_{\mathcal M_{C},\mathcal K^s}(\bar g_{\underline{b}})
=\psi^*_{\mathcal M_{C},\mathcal K^s}(\bar g_{\underline{a}}\cdot\bar g_{\underline{b}})
=c_{\underline{a},\underline{b}}\psi^*_{\mathcal M_{C},\mathcal K^s}(\bar g_{\underline{a}+\underline{b}}).
$$
So by abuse of notation we write for $\underline{a}\in \Gamma$ just $\bar g_{\underline{a}}\in \mathbb K[Z_{\mathcal K^s}]$
for the function $\psi^*_{\mathcal M_C,\mathcal K^s}(\bar g_{\underline{a}})$, where $C$ is a chain in $A$ containing $\supp\underline{a}$. 

(3).  The variety $Z_{\mathcal K^s}$ satisfies the property \ref{varietycondition}.\ref{varietycondition3}: 
By Proposition~\ref{finitemap}, the map $\psi_{\mathcal K^s_0,\mathcal K^s}:Z_{\mathcal K^s}\rightarrow Z_{\mathcal K^s_0}$ is finite. We adopt the notation $C_p$ from the proof of Proposition \ref{finitemap}.

For an element $q\in A$ consider the projection $p_{C_q}: Z_{\mathcal K^s_0}=\prod_{C_p\in \mathcal K^s_0} Y_{C_p}\rightarrow Y_{C_q}$.
The coordinate ring $Y_{C_q}$ contains the class of the extremal function $\bar f_q$, by abuse of notation we also write
$\bar f_q\in \mathbb K[Z_{\mathcal K^s_0}]$ for its image via $p^*_{C_q}$. Due to the tensor product structure of $\mathbb K[Z_{\mathcal K^s_0}]$
(compare to Example~\ref{lengthzero}), the classes of the extremal functions $\bar f_p$, $p\in A$, generate in $\mathbb K[Z_{\mathcal K^s_0}]$ 
a polynomial algebra of dimension $\# A$. 

By the finiteness
of $\psi_{\mathcal K^s_0,\mathcal K^s}$, the classes of the extremal functions  $\bar f_p\in \mathbb K[Z_{\mathcal K^s}]$ (or rather their images
via $\psi^*_{\mathcal K^s_0,\mathcal K^s}$) generate a polynomial algebra in $\mathbb K[Z_{\mathcal K^s}]$, which by Corollary~\ref{dimensionforZM} 
is of the same dimension
as $\dim Z_{\mathcal K^s}$. It follows: there exists an open and dense subset $U\subseteq Z_{\mathcal K^s}$ such that 
$\bar f_p(z)\not=0$ for all $z\in U$.

Let $\underline{a}$ be an element in $\Gamma$. By property \ref{varietycondition}.\ref{varietycondition2}, 
the multiplication relations also hold in $\mathbb K[Z_{\mathcal K^s}]$.
An appropriate power of $\bar g_{\underline{a}}$ (viewed as a function on $Z$) is hence in $\mathbb K[Z]$ a non-zero scalar multiple of a product 
the $\bar f_p$'s with $p\in\supp \underline{a}$. It follows: $\bar g_{\underline{a}}(z)\not=0$ for all $z\in U$ and all $\underline{a}\in \Gamma$.
\end{proof}
The proof that $Z_{\mathcal K^s}$ has the desired properties finishes thus the proof of Theorem~\ref{fanAndDegeneratetheorem}.

\section{Flat degenerations}\label{flatdegen}

In this section we construct a flat degeneration from $X$ to the reduced projective scheme $\mathrm{Proj}(\mathrm{gr}_{\mathcal{V}}R)$ using a Rees algebra construction. We start by extending the valuation image by a total degree.

For $\underline{a}\in \Gamma$, there exists a homogeneous
function $g\in R$ (Lemma~\ref{homogeneousparts2}) such that $\mathcal V(g)=\underline{a}$ with degree given by the formula in Corollary~\ref{degreerelation}. This motivates the following definition of an extra $\mathbb{N}$-grading on $\Gamma$.
Let $\mathcal J\subseteq \mathbb N\times \Gamma$ be the subset
$$
\mathcal J=\left\{(m,\underline{a})\in \mathbb N\times \Gamma \left\vert \ \underline{a}\in\Gamma,
m=\sum_{p\in\supp\underline{a}} a_p\deg f_p\right.\right\}.
$$
Denote by $\succ$ the lexicographic order on the set $\mathbb N\times \Gamma$. 
For a pair $(m,\underline{a})\in \mathbb N\times \Gamma$ let $\mathcal I_{\succeq (m,\underline{a})}$ be the following homogeneous ideal in $R$: 
$$
\mathcal I_{\succeq (m,\underline{a})}=\left\langle g\in R\mid g\textrm{\rm\ homogeneous and\ }\left\{
\begin{array}{l}
\textrm{\rm either\ } \deg g=m \textrm{\rm\ and\ } \mathcal V(g)\ge\underline{a}\\
\textrm{\rm or\ }\deg g>m\\
\end{array}\right.
\right\rangle.
$$
It follows $\mathcal I_{\succeq (m,\underline{a})}\mathcal I_{\succeq (k,\underline{b})}\subseteq  
\mathcal I_{\succeq (m+k,\underline{a}+\underline{b})}$, so the ideals define a filtration on the ring $R$.

We define $\mathcal I_{\succ (m,\underline{a})}$ similarly. 
By Lemma~\ref{NuLeaves}, the quotient
\begin{equation}\label{umnummerierung}
\mathcal I_{\succeq (m,\underline{a})}/\mathcal I_{\succ (m,\underline{a})}=\begin{cases} \{0\}, & \text{if}\ (m,\underline{a})\not\in \mathcal J; \\ R_{\ge \underline{a}}/R_{> \underline{a}}, & \text{if}\ (m,\underline{a})\in \mathcal J.
\end{cases}
\end{equation}
In particular, if $(m,\underline{a})\in \mathcal J$, the quotient space $\mathcal I_{\succeq (m,\underline{a})}/\mathcal I_{\succ (m,\underline{a})}$ is one-dimensional.
 
Let $\pi:\mathcal J\rightarrow \mathbb N$ be an enumeration of the countable many elements respecting the total order,
i.e. $\pi((m,\underline{a}))<\pi((k,\underline{b}))$ if and only if 
$(m,\underline{a})\prec(k,\underline{b})$ with respect to the lexicographic order. We write sometimes just $\mathcal I_j$ if $\pi((m,\underline{a}))=j$ 
instead of $\mathcal I_{\succeq(m,\underline{a})}$. In this way we get a decreasing filtration:
$$
R=\mathcal I_0 \supseteq\mathcal I_1\supseteq\mathcal I_2\supseteq\ldots.
$$
Let 
$$\mathcal A=\cdots\oplus R t^2\oplus R t\oplus R \oplus\mathcal I_1 t^{-1}\oplus\mathcal I_2 t^{-2}\oplus\cdots \subseteq R[t,t^{-1}]$$ 
be the associated Rees algebra.
The natural inclusion $\mathbb K[t]\hookrightarrow \mathcal A$ induces a morphism $\phi:\mathrm{Spec}(\mathcal A)\rightarrow \mathbb A^1$.

The ring $R=\mathbb{K}[\hat{X}]$ is clearly graded since $\hat{X}$ is the cone over the projective variety $X$. We extend this grading to $R[t,t^{-1}]$ by declaring that the degree of $t$ is $0$. Being generated by homogeneous elements, the ideals $\mathcal{I}_j$ are homogeneous, hence the subalgebra $\mathcal{A}\subseteq R[t,t^{-1}]$ is graded. So we have natural $\mathbb{G}_m$-actions on $\mathrm{Spec}(\mathcal A)$ and on $\mathbb{A}^1$; note that the latter is trivial since $\mathbb{K}[t]$ is in degree $0$.

\begin{theorem}\label{Thm:SemiToricDegen}
The morphism $\phi$ is flat and $\mathbb{G}_m$--equivariant. The general fibre for $t\neq 0$ is isomorphic to $\hat X$, the special fibre for $t=0$ is isomorphic to $\mathrm{Spec}(\mathrm{gr}_{\mathcal V}R)$.
\end{theorem}

\begin{proof}
The ring $\mathcal A$ is a torsion-free $\mathbb K[t]$ module by the inclusion $\mathcal A\subseteq R[t,t^{-1}]$, and hence it is a flat module.
For $b\not=0$ it is easily seen that $\mathcal A/(t-b)\simeq R$, and for $t=0$ we get 
$$
\mathcal A/(t)\simeq R/\mathcal I_1\oplus \mathcal I_1/\mathcal I_2 \oplus \mathcal I_2/\mathcal I_3\oplus\cdots\simeq 
\bigoplus_{(\sum_{p\in A}a_p\deg f_p,\underline{a})\in \mathcal J} R_{\ge \underline{a}}/R_{> \underline{a}}\simeq \mathrm{gr}_{\mathcal V}R.
$$
The morphism $\phi$ is $\mathbb{G}_m$--equivariant simply because its fibers are stable by the $\mathbb{G}_m$--action on $\mathrm{Spec}(\mathcal A)$ since $\mathbb{K}[t]$ is in degree $0$.
\end{proof}

An immediate consequence of the $\mathbb{G}_m$--equivariance of $\phi$ is the existence of an induced morphism $\psi:\mathrm{Proj} (\mathcal A)\rightarrow\mathbb A^1$; in particular we get the following result.

\begin{theorem}\label{Thm:Degeneration}
The morphism $\psi$ is flat. The general fibre for $t\neq0$ is isomorphic to $X$, and the special fibre for $t=0$ is isomorphic to $\mathrm{Proj}(\mathrm{gr}_{\mathcal V}R)$.
\end{theorem}

Combining Theorem \ref{Thm:Degeneration} with the existence of a generic hyperplane stratification in Proposition \ref{Prop:Generic}, we have:

\begin{coro}\label{Cor:SemiToric}
Every embedded projective variety $X \subseteq \mathbb{P}(V)$, smooth in codimension one, admits a flat degeneration into $X_0$, a reduced union of projective toric varieties. Moreover, $X_0$ is equidimensional, the number of irreducible components in $X_0$ coincides with the degree of $X$.
\end{coro}

The varieties appearing as irreducible components of $X_0$ are in general not linear because the valuations $\mathcal{V}_{\mathfrak{C}}$ may take values in different lattices. See next section, especially Section \ref{Sec:DegFormulaGHS} for a convex geometric aspect of this corollary.

As mentioned in Remark \ref{Rmk:Hibi}, the above corollary resembles a geometric counterpart of the result by Hibi \cite{Hibi} that finitely generated positively graded rings admit Hodge algebra structures.

\section{The Newton-Okounkov simplicial complex}\label{Nocomplex}

In the theory of Newton-Okounkov bodies \cite{O1, KK,LM}, one associates a convex body to a valuation on a positively graded algebra by taking the convex closure of the degree-normalized valuation images. The (normalized) volume of this convex body computes the leading coefficient of the Hilbert polynomial of the graded algebra.

In our setup, the image $\Gamma:=\{\mathcal V(g)\mid g\in R\setminus\{0\}\}\subseteq \mathbb Q^A$ of the quasi-valuation $\mathcal{V}$ is not necessarily a monoid. We need variations of the approaches in \emph{loc.cit.} to define an analogue of the Newton-Okounkov body for a quasi-valuation, so that its volume computes the degree of the embedded projective variety $X$ in $\mathbb{P}(V)$.

From now on we will work in real vector spaces with the usual Euclidean topology. The spirit of the construction is much the same as in \cite[Section 3]{Dehy}.

\subsection{The order complex}
Recall that in Section \ref{Sec:FanAlgbera}, to every chain $C$ in $A$ we have associated a cone $K_C$ in $\mathbb R^A$ defined by
$$
K_C:=\sum_{p\in C}\mathbb R_{\ge 0}e_p.
$$
The collection of these cones, together with the origin $\{0\}$, is the fan $\mathcal{F}_A$ associated to the poset $A$. Its maximal cones are 
the cones $K_{\mathfrak C}$ associated to the maximal chains in $A$. Each of them comes endowed with a submonoid, the monoid 
$\Gamma_{\mathfrak C}\subseteq K_{\mathfrak C}$ (see Section \ref{Sec:FanAlgbera}).

The order complex $\Delta(A)$ associated to the poset $(A,\leq)$ is the simplicial complex having $A$ as vertices and all chains $C\subseteq A$ as faces, i.e. $\Delta(A)=\{C\subseteq A\mid C\text{ is a chain}\}$.
A geometric realization of $\Delta(A)$ can be constructed by intersecting the cones $K_C$ with appropriate hyperplanes:
for a chain $C\subseteq A$ denote by $\Delta_{C}\subseteq  \mathbb R^{A}$
the simplex: 
\begin{equation}\label{chainsimplex}
\Delta_{C}:=\textrm{convex hull}\ \left\{\frac{1}{\deg f_p} e_p\mid p\in C\right\}.
\end{equation}
The union of the simplexes 
$$|\Delta(A)|=\bigcup_{C\subseteq A\textrm{\ chain}}\Delta_{C}\subseteq\mathbb{R}^A$$ 
is the desired \emph{geometric realization} of $\Delta(A)$. The maximal simplexes are those $\Delta_{\mathfrak{C}}$ 
arising from maximal chains $\mathfrak{C}$ in $A$.

\subsection{The Newton-Okounkov simplicial complex}\label{simplicialcomplex1}
Let $\mathfrak C=\{p_r>\ldots>p_0\}$ be a maximal chain. To avoid an inundation with indexes, we abbreviate throughout this section the natural basis $e_{p_j}$ of $\mathbb R^{\mathfrak C}$ by $e_j$, for $j=0,\ldots,r$. Similarly, the functions $f_{p_j}$ (resp. the bonds $b_{p_j,p_{j-1}}$) will be simplified to $f_j$ (resp. $b_j$) for $j=0,\ldots,r$.

Let $\deg: \mathbb R^A\rightarrow \mathbb R$ be the degree function defined by: for $\underline{a}\in\mathbb{R}^A$, set
$$\deg \underline{a}:=\sum_{p\in A} a_p\deg f_p.$$ 
If $g\in R\setminus\{0\}$ is homogeneous, then
Corollary~\ref{degreerelation} implies that $\deg \mathcal V(g)$ is the degree of $g$.

Recall that $\mathcal{C}$ is the set of all maximal chains in $A$.
\begin{definition}\label{Defn:NOSC}
The \emph{Newton-Okounkov simplicial complex} $\Delta_{\mathcal V}$ associated to the quasi-valuation $\mathcal{V}$ is defined as
$$
\gls{DV}:=\overline{\left\{\frac{\underline{a}}{\deg\underline{a}}\,|\,\underline{a}\in\Gamma\setminus \{0\}\right\}}\subseteq \mathbb R^A.
$$
\end{definition}

\begin{rem}
\begin{enumerate}
\item It is straightforward to show that 
$$\Delta_{\mathcal{V}}=\bigcup_{\mathfrak C\in\mathcal{C}}
\overline{ \bigcup_{m\ge 1}
\left\{
\frac{1}{m}\,\underline{a}\mid\underline{a}\in \Gamma_{\mathfrak C},\deg \underline{a}=m \right\} 
}.
$$
\item In the definition of the Newton-Okounkov bodies in \cite{O1, KK,LM} as closed convex hulls of points, the convex hull operation is not necessary (see for example \cite{O2}).
\end{enumerate}
\end{rem}

Recall that a simplicial complex $\Delta$ of dimension $r$ is called \emph{homogeneous} when for any face $F$ of $\Delta$, there exists an $r$-simplex containing $F$ as a face.

\begin{proposition}\label{simplicialdecomp}
We have $\Delta_{\mathcal V}=|\Delta(A)|$. In particular,
$\Delta_{\mathcal V}$ is a homogeneous simplicial complex of dimension $r$.
\end{proposition}

\begin{proof}
We show $\Delta_{\mathcal V}\subseteq |\Delta(A)|$. Let $g\in R\setminus\{0\}$ be a homogeneous element with $\mathcal{V}(g)=\underline{a}\in\Gamma_{\mathfrak{C}}$. By Corollary~\ref{powerrelation2}, there exists $m\in\mathbb{N}$ such that for any $p\in\mathfrak{C}$, $ma_p\in\mathbb{N}$ and $\mathcal V(g^m)=\mathcal V(\prod_{p\in \mathfrak{C}}f_p^{ma_p})$. Noticing that 
$$\frac{\mathcal V(g)}{\deg g} =\frac{\mathcal V(g^m)}{m\deg g},$$
it suffices to assume that $g=\prod_{p\in C}f_p^{n_p}$ for $n_p\in\mathbb{N}$ and a chain $C\subseteq A$. Let $m=\deg g=\sum_{p\in C}n_p\deg f_p$ be the degree of $g$. If $\mathfrak C$ is a maximal chain containing $C$, then 
$$
\frac{\mathcal V(g)}{\deg g}=\frac{\mathcal V_{\mathfrak C}(g)}{m}=\sum_{p\in C} \frac{n_p\deg f_p}{m}\left(\frac{1}{\deg f_p}e_p\right)\in \Delta_C\subseteq |\Delta(A)|.
$$
Conversely, we show $|\Delta(A)|\cap\mathbb{Q}^A\subseteq\Delta_{\mathcal{V}}$. For a chain $C\subseteq A$, a point in $\Delta_C\cap\mathbb{Q}^A\subseteq |\Delta(A)|\cap\mathbb{Q}^A$ is a convex linear combination 
$$\sum_{p\in C} a_p\frac{1}{\deg f_p}e_p,\ \text{with}\ a_p\in\mathbb{Q}_{\geq 0}\ \text{and}\ \sum_{p\in C} a_p=1.$$
Choose a non-zero $m\in\mathbb N$ to be such that for all $p\in C$: $m\frac{a_p}{\deg f_p}\in\mathbb N$.
Let $g=\prod_{p\in C}f_p^{m\frac{a_p}{\deg f_p}}$. It is now easy to verify that 
$$\frac{\mathcal V(g)}{\deg g}=\frac{\mathcal V(g)}{m}=\sum_{p\in C} a_p\frac{1}{\deg f_p}e_p.$$
This terminates the proof since $\Delta_{\mathcal{V}}$ is closed.
\end{proof}

\subsection{Rational and integral structures on simplexes}\label{sec:volumeandcovolume}

Let $H_R(t)$ be the Hilbert polynomial for the homogeneous cooordinate ring $R=\mathbb K[\hat{X}]$. Our aim is to translate
the calculation of its leading coefficient into a problem of calculating the leading coefficient of certain Ehrhart quasi-polynomials. Let $\Gamma_n=\{\underline{a}\in\Gamma\mid \deg \underline{a}=n\}$ be the elements in $\Gamma$ of degree $n$.
Since the leaves of the quasi-valuation are one dimensional, one has 
$H_R(n)=\# \Gamma_n=\# \Delta_{\mathcal V}(n)$ for $n$ large, where $\Delta_{\mathcal V}(n)$ is defined as the intersection
$n\Delta_{\mathcal V}\cap \Gamma$.

We have already pointed out that the lattice $L^{\mathfrak{C}}$ defined in \eqref{latticedef} is in general too large compared to the monoid
$\Gamma_{\mathfrak C}$. In the following let $\gls{LUC}\subseteq L^{\mathfrak{C}}$ be the sublattice generated by 
$\Gamma_{\mathfrak C}$. 
The affine span of the simplex $\Delta_{\mathfrak C}$ is the affine subspace $U_1=\frac{1}{\deg f_0}e_0+U_0$ of $\mathbb{R}^{\mathfrak{C}}$,
where $U_0$ is the linear subspace 
$$
U_0:=\mathrm{span}_\mathbb{R}\left\{\frac{1}{\deg f_j}e_j-\frac{1}{\deg f_0}e_0\mid j=1,\ldots,r\right\}.
$$
This linear subspace $U_0\subseteq \mathbb{R}^{\mathfrak{C}}$ can also be characterized as the kernel of the 
degree function $\deg:\mathbb{R}^{\mathfrak{C}}\to\mathbb{R}$, $\underline{a}\mapsto \sum_{p\in \mathfrak{C}} a_p\deg f_p$.
Hence for $m\in\mathbb{Z}$, the affine subspace $U_m$ of $\mathbb R^{\mathfrak C}$ of elements of degree $m$ is given by:
$$
U_m=\frac{m}{\deg f_0}e_0+U_0=\{\underline{a}\in \mathbb R^{\mathfrak C}\mid \deg \underline{a}=m\}.
$$

\subsubsection{The projection}\label{sec:projection}

Let $\mathfrak C=(p_r,\ldots,p_0)$ be a maximal chain in $A$.
The degree function $\deg:\mathbb{R}^{\mathfrak{C}}\to\mathbb{R}$ takes integral values on $\Gamma_{\mathfrak C}$, and hence the function 
$\deg$ takes also integral values on $\mathcal{L}^{\mathfrak{C}}$.  Set $\mathcal L^{\mathfrak C}_0=U_0\cap \mathcal L^{\mathfrak C}$. For later purpose
it is important to know that $U_m\cap \mathcal L^{\mathfrak C}\not=\emptyset$ for all $m\in\mathbb Z$, this is a consequence
of the following lemma: 

\begin{lemma}\label{degreeonelattice}
We have: $\frac{1}{\deg f_0}e_0\in  U_{1}\cap \mathcal{L}^{\mathfrak C}$.
\end{lemma}
\begin{proof}
We fix an enumeration of the length $0$ elements in $A$, say: $q_1=p_0,q_2,\ldots,q_s$, and identify $q_j$ with the unique point in $X_{q_j}$. For $i=2,\ldots,s$ let $h_i\in  V^*$ be a linear function such that $h_i(q_j)\neq 0$ for 
$j=1,\ldots,s$, $j\neq i$, and $h_i(q_i)=0$. The function
\begin{equation}\label{pointfunction}
g=\prod_{j=2,\ldots,s} h_j
\end{equation}
vanishes in $q_2,\ldots,q_s$, but not in $p_0=q_1$. 

Let $\mathfrak C'=(p'_r,\ldots,p'_0)$ be another maximal chain and $h\in V^*$ a linear function which does not vanish in any of the points $q_1,\ldots,q_s$. The same calculation as in Example \ref{Ex:Bertini2} shows that 
$$
\mathcal V(g)=\frac{s-1}{\deg f_{p_0}}e_0 \ \ \text{and}\ \ \mathcal V(gh)=\frac{s}{\deg f_{p_0}}e_0.
$$
These elements are both in $\Gamma_{\mathfrak C}$, hence $\frac{1}{\deg f_{p_0}}e_0$ is an element in the lattice $\mathcal L^{\mathfrak C}$ generated by $\Gamma_{\mathfrak C}$.
\end{proof}

Set $\ell_1=\frac{1}{\deg f_{p_0}}e_0$.
The projection $\mathrm{pr}_{\ell_1}:\mathbb{R}^{\mathfrak{C}}\rightarrow \mathbb{R}^{\mathfrak{C}\setminus \{p_0\}}$ $(\simeq\mathbb{R}^{\mathfrak{C}}/\mathbb R\ell_1$)
induces an isomorphism of vector spaces:
$$
\mathrm{pr}_{\ell_1}\vert_{U_0}:U_0\stackrel{\sim}{\longrightarrow} \mathbb{R}^{\mathfrak{C}\setminus \{p_0\}}.
$$ 
We identify for convenience in the following $U_0$ sometimes with $\mathbb{R}^{\mathfrak{C}\setminus \{p_0\}}$, and we identify $\mathcal L^{\mathfrak C}_0\subseteq U_0$
with its image in $\mathbb{R}^{\mathfrak{C}\setminus \{p_0\}}$. Note that $\ell_1$ is primitive, so the image of
$\mathcal L^{\mathfrak C}_0$ is isomorphic to $\mathcal L^{\mathfrak C}/\mathbb Z \ell_1$.

For all $m\ge 1$, the restriction $\mathrm{pr}_{\ell_1}\vert_{U_m}:U_m\rightarrow U_0$ induces a bijection between the affine subspace $U_m$ and the 
vector space $U_0$. An element $\ell_m\in \mathcal L^{\mathfrak C}\cap U_m$ can always be 
written as $\ell_m=m\ell_1 +\ell_0$ for some $\ell_0\in \mathcal L^{\mathfrak C}_0$. The bijection $\mathrm{pr}_{\ell_1}\vert_{U_m}$ induces 
hence in this case a bijection between $\mathcal L^{\mathfrak C}\cap U_m$  and the lattice $\mathcal L^{\mathfrak C}_0\subseteq U_0$.

The map $\mathrm{pr}_{\ell_1}$ is $\mathbb R$-linear, so the restriction $\mathrm{pr}_{\ell_1}\vert_{U_1}$ preserves the notion of the 
affine convex linear hull of elements in $U_1$.
Set $\Delta_{\mathfrak C}^0:=\mathrm{pr}_{\ell_1}\vert_{U_1}(\Delta_{\mathfrak C})\subseteq U_0$, then
$$
\Delta_{\mathfrak C}^0=\text{convex hull}\left\{\frac{1}{\deg f_p} \mathrm{pr}_{\ell_1}(e_p)\mid p\in\mathfrak C\right\}=
\text{convex hull}\left\{0,\frac{1}{\deg f_p} e_{p_j}\mid j=1,\ldots,r\right\}.
$$
The same arguments hold for all $m\ge 1$: 
$$
\mathrm{pr}_{\ell_1}\vert_{U_m}(m\Delta_{\mathfrak C})=m\Delta_{\mathfrak C}^0.
$$

From this construction we have:
\begin{lemma}\label{equalityofpoints}
For all $n\ge 1$, $\# \left(n\Delta_{\mathfrak C}\cap \mathcal L^{\mathfrak C}\right)=\# \left(n\Delta_{\mathfrak C}^0\cap \mathcal L^{\mathfrak C}_{0}\right)$.
\end{lemma}

\subsubsection{Rational and integral structure}\label{integral_structure}
By fixing an ordered basis of $\mathcal L^{\mathfrak C}_{0}$, one gets an isomorphism $\Psi: U_0\rightarrow \mathbb R^{r}$ which identifies
the lattice $\mathcal{L}^{\mathfrak C}_{0}\subseteq U_0$ with the lattice $\mathbb Z^{r}\subseteq \mathbb R^{r}$. The simplex
$\Delta_{\mathfrak C}^0\subseteq U_0$ can be identified in this way with the simplex $D_{\mathfrak C}\subseteq \mathbb R^{r}$, defined as:
$$
D_{\mathfrak C}=\textrm{\rm convex hull}\ \left\{0,\frac{1}{\deg f_p} \Psi(e_p)\mid p\in\mathfrak C\setminus\{p_0\}\right\}\subseteq \mathbb R^{r}.
$$
The simplex has rational vertices. Now Lemma~\ref{equalityofpoints} implies for all $n\in\mathbb N$:
\begin{equation}\label{pointsareequal}
\# \left(n\Delta_{\mathfrak C}\cap \mathcal L^{\mathfrak C}\right)=\# \left(nD_{\mathfrak C}\cap \mathbb Z^{r}\right).
\end{equation}

The correspondence in \eqref{pointsareequal} generalizes \cite[Definition 3.1]{Dehy}; we call it a \textit{rational structure} on $\Delta_{\mathfrak C}$.
In case the vertices of $D_{\mathfrak C}$ are integral points, it will be called an \textit{integral structure} 
on $\Delta_{\mathfrak C}$. 

The map $\Psi$ and the simplex $D_{\mathfrak C}$ depend on the choice of the 
ordered basis of $\mathcal L^{\mathfrak C}_{0}$. But different choices of ordered bases lead to unimodular equivalent
simplexes.

\subsection{The degree formula}
Let $\Delta_{\mathcal V}$ be the Newton-Okounkov simplicial complex associated to $X$ and the quasi-valuation $\mathcal V$.
For each maximal chain $\mathfrak C$ let  $\Delta_{\mathfrak C}\subseteq  \Delta_{\mathcal V}$ be the associated simplex
given by the simplicial decomposition in Proposition~\ref{simplicialdecomp}. We fix for all maximal chains $\mathfrak C$ in $A$
a rational structure on $\Delta_{\mathfrak C}$ (see subsection~\ref{integral_structure}), denote by $D_{\mathfrak C}\subseteq \mathbb R^{r}$
the corresponding simplex with rational vertices. Let $ \textrm{\rm vol}(D_{\mathfrak C})$ be the Euclidean volume of the simplex. This is also the normalized volume of $\Delta_{\mathfrak{C}}^0$ with respect to the lattice $\mathcal{L}_0^{\mathfrak{C}}$.

\begin{theorem}\label{volumetheorem1}
The degree of the embedded variety $X\hookrightarrow\mathbb P(V)$ is equal to 
$$
\deg X=r!\sum_{\mathfrak C\in\mathcal{C}} \textrm{\rm vol}(D_{\mathfrak C}).
$$
\end{theorem}

\begin{proof}
Let $H_R(t)$ be the Hilbert polynomial of the homogeneous coordinate ring $R=\mathbb K[\hat{X}]$.
The degree of $X$ can be computed as $r!\cdot c_{ {r}}$, where $c_{ {r}}$ is the leading coefficient of $H_R(t)$.
By Theorem~\ref{fanAndDegeneratetheorem},  the Hilbert polynomial $H_R(t)$ and the Hilbert polynomial $H_\Gamma(t)$
associated to the fan algebra $\mathbb K[\Gamma]$ coincide. 

Let $H_{\mathfrak C}(t)$ be the Hilbert quasi-polynomial of the algebra $\mathbb K[\Gamma_{\mathfrak C}]$.
The leading coefficient of $H_\Gamma(t)$ is the sum of the leading coefficients of the $H_{\mathfrak C}(t)$,
where the sum runs over all maximal chains $\mathfrak C$ in $A$. By Lemma \ref{fan_monoid:Hilbert-quasi2}, the leading coefficient of $H_{\mathfrak C}(t)$
is equal to the leading coefficient of the Hilbert quasi-polynomial $\tilde{H}_{\mathfrak C}(t)$ of the algebra 
$\mathbb K[\tilde\Gamma_{\mathfrak C}]$.

In the sections~\ref{simplicialcomplex1}, \ref{sec:projection} and \ref{integral_structure} we construct a simplex $D_{\mathfrak C}\subseteq \mathbb Q^r$ such that
$\# (\tilde\Gamma_{\mathfrak C})_{m}=\# \left(m\Delta_{\mathfrak C}\cap \mathcal L^{\mathfrak C}\right)$
is equal to $\# \left(mD_{\mathfrak C}\cap \mathbb Z^{r}\right)$, which implies that the leading coefficient of $\tilde{H}_{\mathfrak C}(t)$
is the same as the leading coefficient of the Ehrhart quasi-polynomial 
$\mathrm{Ehr}_{D_{\mathfrak C}}(t)$ associated $D_{\mathfrak C}$. In this case the leading coefficient  is 
the Euclidean volume of $D_{\mathfrak{C}}$ divided by the co-volume  of the lattice $\mathbb Z^{r}$, which finishes the proof of the theorem.
\end{proof}

The estimates above can be made more precise if the monoid $\Gamma_{\mathfrak{C}}$
is saturated, that is to say, $\mathcal L^{\mathfrak C}\cap  K_{\mathfrak C} =\Gamma_{\mathfrak C}$, for all maximal chains $\mathfrak{C}$. Algebraically it is equivalent to say that the algebra $\mathbb{K}[\Gamma_{\mathfrak{C}}]$ is normal.
Geometrically this condition is equivalent to the normality of all the irreducible components of the degenerate variety
$\mathrm{Spec}(\mathrm{gr}_{\mathcal V}R)$. We give a name to such a situation:

\begin{definition}\label{def:stratification:normal}
A Seshadri stratification is called \emph{normal} if $\Gamma_\mathfrak{C}$ is saturated for every maximal chain $\mathfrak{C}$.
\end{definition}

Further results on normal Seshadri stratifications can be found in \cite{CFL3}. If this condition is fulfilled, then the number of points $\# (n\Delta_{\mathcal V}\cap \Gamma)$ can be determined as an alternating
sum of the set of lattice points in all possible intersections of the simplexes.

\begin{proposition}\label{volumetheorem2}
If the Seshadri stratification is normal, then the Hilbert function $H_R(t)$ of the graded ring $R=\mathbb K[X]$ is an alternating sum of Ehrhart quasi-polynomials. More precisely, 
let $\mathfrak C_1,\ldots,\mathfrak C_t$ be an enumeration
of the maximal chains in $A$, then
$$
H_R(t)=\sum_{1\le i_1<\ldots<i_\ell\le r} (-1)^{\ell -1} \mathrm{Ehr}_{\Delta_{\mathfrak C_{i_1}}\cap\ldots\cap \Delta_{\mathfrak C_{i_\ell}}}(t).
$$
\end{proposition}

\subsection{The degree formula and the generic hyperplane stratification}\label{Sec:DegFormulaGHS}
We fix a generic hyperplane stratification as in the proof of Proposition \ref{Prop:Generic} and choose the functions $f_{0,k}$ as in Example \ref{Ex:Bertini1}.

To see the connection with the degree formula above, note that by Lemma~\ref{Lem:ImageValC}, for every maximal
chain ${\mathfrak C}=(q_r,\ldots,q_1,q_{0,j})$ the valuation ${\mathcal V}_{\mathfrak C}$ takes values in the lattice:
$$
L^{\mathfrak C}=\{(a_r,\ldots,a_0)\in \mathbb Q^{\mathfrak C}\,\vert\, (s-1)a_0\in\mathbb{Z},\, a_j\in\mathbb{Z},\ 1\leq j\leq r\}.
$$
Since $\mathcal{L}^{\mathfrak C}\subseteq L^{\mathfrak C}$, Example~\ref{extremalvaluation} and 
Lemma~\ref{degreeonelattice} imply that one has indeed equality $\mathcal{L}^{\mathfrak C}= L^{\mathfrak C}$.

It follows that the projection $\mathrm{pr}_{\ell_1}:\mathbb{R}^{\mathfrak{C}}\rightarrow \mathbb{R}^{\mathfrak{C}\setminus \{q_{0,j}\}}$ (Section~\ref{sec:projection}) induces an identification
of $\mathcal L_0^{\mathfrak C}\simeq \mathcal L^{\mathfrak C}/\mathbb Z \ell_1$ with $\mathbb Z^r$, and 
$\Delta_{\mathfrak C}^0:=\mathrm{pr}_{\ell_1}\vert_{U_1}(\Delta_{\mathfrak C})$ is just the standard simplex
with vertices $0, e_{q_1},\ldots,e_{q_r}$. Therefore in this case we can take $\Psi$ as the identity map,
and $D_{\mathfrak C}=\Delta_{\mathfrak C}^0$ is a standard simplex of dimension $r$ in $\mathbb{R}^{\mathfrak{C}\setminus \{q_{0,j}\}}$.

As a summary: For every maximal chain in $A$ we get as simplex a standard simplex, its Euclidean volume is $\frac{1}{r!}$. There are $s$ maximal chains in $A$. The degree formula in Theorem \ref{volumetheorem1} reproduces hence the predicted number.

\section{Projective normality}\label{Sec:ProjNormal}

In this section, we give a criterion on the projective normality of the projective variety $X\subseteq\mathbb{P}(V)$ using the normality of the toric varieties $\mathrm{Spec}(\mathbb{K}[\Gamma_{\mathfrak{C}}])$ and the topology of the poset $A$.

Let $\mathrm{SR}(A)$ be the \emph{Stanley-Reisner algebra} of the poset $A$. By definition
$$\mathrm{SR}(A):=\mathbb{K}[t_p\mid p\in A]/(t_pt_q\mid p\text{ and } q \text{ are incomparable}).$$
The Stanley-Reisner algebra has a linear basis consisting of $t_{\underline{a}}:=t_{p_r}^{a_r}\cdots t_{p_0}^{a_0}$ where $\mathfrak{C}=(p_r,\ldots,p_1,p_0)$ runs over all maximal chains in $A$ and $\underline{a}=(a_r,\ldots,a_0)\in\mathbb{N}^{\mathfrak{C}}$.

The poset $A$ is called \emph{Cohen-Macaulay} over $\mathbb{K}$, if the algebra $\mathrm{SR}(A)$ is Cohen-Macaulay. As an example, a shellable poset $A$ is Cohen-Macaulay over any field $\mathbb{K}$ \cite{Bj}.

The goal of this section is to prove the following result:

\begin{theorem}\label{Thm:ProjNormal}
If the Seshadri stratification is normal and the poset $A$ is Cohen-Macaulay over $\mathbb{K}$, then
\begin{itemize}
\item[{i)}] the ring $R$ is normal, hence $X\subseteq\mathbb{P}(V)$ is projectively normal;
\item[{ii)}] the special fibre $X_0$ is Cohen-Macaulay.
\end{itemize}
\end{theorem}

\begin{proof}
The proof follows the ideas in \cite{Chi}. We divide it into several steps. 
\vskip 5pt
\noindent \emph{Step 1}. We construct an embedding of $\mathbb{K}$-algebras $\phi: \mathbb{K}[\Gamma]\to\mathrm{SR}(A)$. This endows $\mathrm{SR}(A)$ with a $\mathbb{K}[\Gamma]$-module structure.
    
Let $M$ be the product of all bonds in $\mathcal{G}_A$. For a fixed maximal chain $\mathfrak{C}$, by Lemma \ref{Lem:ImageValC}, the monoid $\Gamma_{\mathfrak{C}}$ is contained in the lattice $L^{\mathfrak{C}}$, hence $M\Gamma_{\mathfrak{C}}\subseteq\mathbb{N}^{\mathfrak{C}}$. The embedding of monoids $\Gamma_{\mathfrak{C}}\hookrightarrow \mathbb{N}^{\mathfrak{C}}$, $\underline{a}\mapsto M\underline{a}$
induces injective $\mathbb{K}$-algebra morphism
$\phi_{\mathfrak{C}}:\mathbb{K}[\Gamma_{\mathfrak{C}}] \hookrightarrow\mathbb{K}[\mathbb{N}^{\mathfrak{C}}]$. In view of Corollary \ref{finiteunionmonoid}, we define a morphism of $\mathbb{K}$-algebra $\phi:\mathbb{K}[x_{\underline{a}}\mid\underline{a}\in\Gamma]\to\mathrm{SR}(A)$
in such a way that for a maximal chain $\mathfrak{C}$ and $\underline{a}\in\Gamma_{\mathfrak{C}}$, 
$\phi(x_{\underline{a}}):=\phi_{\mathfrak{C}}(x_{\underline{a}})=t_{M\underline{a}}$.
This map is clearly independent of the choice of the maximal chain so it is well-defined. Since the $\phi_{\mathfrak{C}}$ are injective, $\phi$ passes through $I(\Gamma)$, yields an injective $\mathbb{K}$-algebra morphism from $\mathbb{K}[\Gamma]$ to $\mathrm{SR}(A)$, which is also denoted by $\phi$.
\vskip 5pt
\noindent \emph{Step 2.} We define a $\mathbb{K}[\Gamma]$-module morphism $\psi:\mathrm{SR}(A)\to\mathbb{K}[\Gamma]$.

We start from considering a single maximal chain $\mathfrak{C}$. Let $\psi_{\mathfrak{C}}:\mathbb{K}[\mathbb{N}^{\mathfrak{C}}]\to\mathbb{K}[\Gamma_{\mathfrak{C}}]$ be the $\mathbb{K}$-linear map sending $t_{\underline{a}}\in\mathbb{K}[\mathbb{N}^{\mathfrak{C}}]$ to $x_{\underline{a}/M}$ if $\underline{a}\in M\Gamma_{\mathfrak{C}}$ or to $0$ otherwise. We show that $\psi_{\mathfrak{C}}$ is a morphism of $\mathbb{K}[\Gamma_{\mathfrak{C}}]$-modules, where the $\mathbb{K}[\Gamma_{\mathfrak{C}}]$-module structure on $\mathbb{K}[\mathbb{N}^{\mathfrak{C}}]$ is defined by: for $\underline{a}\in\Gamma_{\mathfrak{C}}$ and $\underline{b}\in\mathbb{N}^{\mathfrak{C}}$, $x_{\underline{a}}\cdot t_{\underline{b}}:=t_{M\underline{a}+\underline{b}}$. It suffices to show that for $\underline{a}\in \Gamma_{\mathfrak{C}}$ and $\underline{b}\in \mathbb{N}^{\mathfrak{C}}$, 
$\psi_{\mathfrak{C}}(t_{M\underline{a}}t_{\underline{b}})=x_{\underline{a}}\psi_{\mathfrak{C}}(t_{\underline{b}})$. 
If $\underline{b}\in M\Gamma_{\mathfrak{C}}$ then there is nothing to show. Assume that $\underline{b}\notin M\Gamma_{\mathfrak{C}}$, then $\psi_{\mathfrak{C}}(t_{\underline{b}})=0$. We need to show that $M\underline{a}+\underline{b}\notin M\Gamma_{\mathfrak{C}}$. Assume the contrary, since $M\Gamma_{\mathfrak{C}}$ is saturated, $\underline{b}=(M\underline{a}+\underline{b})-M\underline{a}\in M\Gamma_{\mathfrak{C}}$ gives a contradiction. 

For the general case, take a monomial $t_{\underline{b}}$ in $\mathrm{SR}(A)$ supported on a maximal chain $\mathfrak{C}$ and define $\psi(t_{\underline{b}}):=\psi_{\mathfrak{C}}(t_{\underline{b}})$. It is independent of the choice of the maximal chain $\mathfrak{C}$, hence the map $\psi$ is well-defined.

To show that $\psi$ is a morphism of $\mathbb{K}[\Gamma]$-modules, take $\underline{a}\in\Gamma$, there exists a maximal chain $\mathfrak{C}'$ such that $\mathrm{supp}(\underline{a})\subseteq\mathfrak{C}'$. Take a monomial $t_{\underline{b}}$ in $\mathrm{SR}(A)$ supported on a maximal chain $\mathfrak{C}$. If the support of $t_{\underline{b}}$ is not contained in $\mathfrak{C}'$, $\phi(x_{\underline{a}})\cdot t_{\underline{b}}=0$. Otherwise we assume that they both supported in $\mathfrak{C}$, and $\psi(\phi(x_{\underline{a}})\cdot t_{\underline{b}})=\psi_{\mathfrak{C}}(t_{M\underline{a}+\underline{b}})=\phi(x_{\underline{a}})\cdot \psi(t_{\underline{b}})$.

One easily verifies that as linear maps $\psi\circ\phi=\mathrm{id}_{\mathbb{K}[\Gamma]}$.
\vskip 5pt
\noindent \emph{Step 3.} The $\mathbb{K}$-algebra $\mathbb{K}[\Gamma]$ is Cohen-Macaulay over $\mathbb{K}$.

We start with the Stanley-Reisner algebra, which is Cohen-Macaulay by assumption. Consider the following elements in $\mathrm{SR}(A)$: for $i=0,1,\ldots,r$ with $r=\dim X$,
$$\ell_i:=\sum_{p\in A,\ \ell(p)=i}t_p.$$
Since $\mathrm{SR}(A)$ is Cohen-Macaulay, $\ell_0,\ell_1,\ldots,\ell_r$ form a regular sequence in $\mathrm{SR}(A)$.

We prove the Cohen-Macaulay-ness of $\mathbb{K}[\Gamma]$ by constructing a regular sequence in $\mathbb{K}[\Gamma]$ of length $r+1$. Then the depth of $\mathbb{K}[\Gamma]$ is greater or equal than $r+1$, while the other inequality always holds.

Since any two elements in $A$ having the same length are incomparable: we consider the elements 
$$\ell_i^M=\sum_{p\in A,\ \ell(p)=i}t_p^M,\ \ i=0,1,\ldots,r.$$
These elements form a regular sequence in $\mathrm{SR}(A)$, and they are contained in the image of $\phi$.

We choose the unique element $u_i\in\mathbb{K}[\Gamma]$ such that $\phi(u_i)=\ell_i^M$, and show that the image of $u_i$ is not a zero divisor in $\mathbb{K}[\Gamma]/(u_0,\ldots,u_{i-1})$. Assume the contrary, there exist $h_0,\ldots,h_{i-1}\in\mathbb{K}[\Gamma]$ such that $u_i h=h_0u_0+\ldots+h_{i-1}u_{i-1}$. Applying $\phi$ gives $\ell_i^M \phi(h)=\phi(h_0)\ell_0^M+\ldots+\phi(h_{i-1})\ell_{i-1}^M$. If $\phi(h)=s_0\ell_0^M+\ldots+s_{i-1}\ell_{i-1}^M$ for some $s_0,\ldots,s_{i-1}\in \mathrm{SR}(A)$, applying $\mathbb{K}[\Gamma]$-module morphism $\psi$ gives
$h=\psi(s_0)u_0+\ldots+\psi(s_{i-1})u_{i-1}$, contradicting to the assumption that $h\notin (u_0,\ldots,u_{i-1})$. Therefore $\phi(h)\notin (\ell_0^M,\ldots,\ell_{i-1}^M)$, contradicts to the fact that $\ell_0^M,\ell_1^M,\ldots,\ell_r^M$ is a regular sequence.
\vskip 5pt
\noindent \emph{Step 4.} The ring $R$ is normal, hence $X\subseteq\mathbb{P}(V)$ is projectively normal.

By Theorem \ref{fanAndDegeneratetheorem}, $\mathrm{gr}_{\mathcal{V}}R$ is isomorphic to the fan algebra $\mathbb{K}[\Gamma]$, it is therefore Cohen-Macaulay. Since being Cohen-Macaulay is an open property, Theorem \ref{Thm:SemiToricDegen} implies that $R$ is Cohen-Macaulay. According to the axiom (S1) of a Seshadri stratification, the ring $R$ is smooth in codimension one. By Serre's criterion, $R$ is normal.
\end{proof}

\section{Standard Monomial Theory}\label{section_standard_monomial_theory}

In this section we first review the results from the previous sections in terms of a weak form of a standard monomial theory for the ring $R = \mathbb{K}[\hat{X}]$. 
Imposing additional conditions on the Seshadri stratification gives stronger versions of this theory. We will discuss two of these enhancements: the normality and the balancing of the stratification. In case these requirements are fulfilled, the ring  $R$ admits a structure closely related to the LS-algebras \cite{Chi}.

\subsection{A basis associated to the leaves}

For $\underline{a} \in \Gamma$ we choose a regular function $x_{\underline{a}}\in R$ with quasi-valuation $\mathcal V(x_{\underline{a}}) = \underline{a}$. By Lemma \ref{NuLeaves}, $R = \bigoplus_{\underline{a}\in\Gamma}\mathbb{K}\, x_{\underline{a}}$ as a $\mathbb{K}$--vector space. In particular, for each pair $\underline{a}, \underline{b} \in \Gamma$, there exists a relation, called \emph{straigthening relation},
\begin{equation}\label{straightening:One}
x_{\underline{a}}\cdot x_{\underline{b}} = \sum_{\underline{a} + \underline{b} \leq^t \underline{c}} u_{\underline{c}}^{\underline{a}, \underline{b}}x_{\underline{c}}
\end{equation}
expressing the product $x_{\underline{a}}\cdot x_{\underline{b}}$ in terms of the $\mathbb{K}$--basis elements. In what follows, each time we refer to the coefficients $u_{\underline{c}}^{\underline{a}, \underline{b}}$ of a straightening relation, we will simply write $u_{\underline{c}}$ by omitting the dependence on the leaves $\underline{a}$ and $\underline{b}$ if this does not create any ambiguity.

The restriction
on the indexes of the possibly non-zero terms in the expression above comes from the fact that $\mathcal V$ is a quasi-valuation (Lemma \ref{simpleproperties}).
This implies that $u_{\underline{c}}\neq 0$ only for those $\underline{c}\in\Gamma$ such that $\underline{a} + \underline{b} \leq^t \underline{c}$. We recall that $\underline{a} + \underline{b} \in \Gamma$ if and only if $\supp\underline{a} \cup \supp\underline{b}$ is a chain in $A$ (Corollary \ref{finiteunionmonoid}). In such a case, the term $x_{\underline{a} + \underline{b}}$ does appear in the straightening relation for $x_{\underline{a}}\cdot x_{\underline{b}}$, i.e. $u_{\underline{a} + \underline{b}} \neq 0$. This term is clearly the leading term in the above straightening relation.

\begin{rem}
In analogy to the theory of LS-algebras developed in \cite{Chi} and \cite{Chi2}, one can 
call the data consisting of the generators $x_{\underline{a}}$ with $\underline{a}\in\Gamma$ and the straightening relations an 
\emph{algebra with leaf basis} over $\Gamma$ for the ring $R$.
Results in previous sections can be proved in the general algebraic context of an algebra with leaf basis; this will appear in a forthcoming article.
\end{rem}

\subsection{Normality and standard monomials}\label{SMT:normal:leaf}

Imposing the normality to the Seshadri stratification allows us to make a natural choice to the monomials, and to define the condition of being standard in a standard monomial theory. Let $\mathfrak{C}$ be a maximal chain in $A$.

\begin{definition}\label{Defn:Indec}
An element $\underline{a}\in\Gamma_\mathfrak{C}$ is called \emph{decomposable} if it is $0$, or if there exist $\underline{a}_1,\,\underline{a}_2\in\Gamma_\mathfrak{C}\setminus\{0\}$ with $\min\supp\underline{a}_1\geq\max\supp\underline{a}_2$ such that $\underline{a} = \underline{a}_1 + \underline{a}_2$. We say that $\underline{a}$ is \emph{indecomposable} if it is not decomposable.
\end{definition}

\begin{proposition}\label{prop:indecomposable}
Each $\underline{a}\in\Gamma_\mathfrak{C}$ has a decomposition $\underline{a} = \underline{a}_1 + \underline{a}_2 + \ldots + \underline{a}_n$ with $\underline{a}_1,\underline{a}_2, \ldots,\underline{a}_n\in\Gamma_\mathfrak{C}$ indecomposable such that $\min\supp\underline{a}_j \geq \max\supp\underline{a}_{j+1}$ for each $j=1, 2,\ldots,n-1$.
\end{proposition}

\begin{proof}
If $\underline{a}$ is $0$ or indecomposable, then there is nothing to prove. So we suppose that $\underline{a}$ is decomposable and proceed by induction on the degree of $\underline{a}$.

Let $\underline{a} = \underline{a}_1 + \underline{a}_2$ with $\underline{a}_1,\underline{a}_2\neq 0$ and $\min\supp\underline{a}_1 \geq \max\supp\underline{a}_2$. By induction, both $\underline{a}_1$ and $\underline{a}_2$ have decompositions into indecomposable elements. Putting them together gives a decomposition of $\underline{a}$.
\end{proof}

The decomposition as in the above proposition may not be unique. The uniqueness holds in an important special case. Recall that $\Gamma_\mathfrak{C}$ is \emph{saturated} if $\mathcal{L}^\mathfrak{C}\cap K_{\mathfrak{C}} = \Gamma_\mathfrak{C}$ (see the paragraphs before Proposition \ref{volumetheorem2}).

\begin{proposition}\label{Prop:UniqueDec}
If $\Gamma_\mathfrak{C}$ is saturated, the decomposition of any element in $\Gamma_\mathfrak{C}$ is unique. 
\end{proposition}

\begin{proof}
Let $\underline{a} = \underline{b}_1 + \ldots + \underline{b}_n = \underline{c}_1 + \ldots + \underline{c}_m$ be two different decompositions of $\underline{a}$ into indecomposable elements. Let $k$ be minimal such that $\underline{b}_k\neq\underline{c}_k$. We can assume that either $\min\supp\underline{b}_k > \min\supp\underline{c}_k$ or that $p = \min\supp\underline{b}_k = \min\supp\underline{c}_k$ but the entries of $\underline{b}_k$ in $p$ is strictly less than that of $\underline{c}_k$. In both cases $\underline{d} = \underline{c}_k - \underline{b}_k\neq 0$ is an element of 
$\mathcal{L}^\mathfrak{C}\cap K_{\mathfrak{C}}=\Gamma_\mathfrak{C}$ and hence $\underline{c}_k = \underline{b}_k + \underline{d}$ with $\min\supp\underline{b}_k \geq \max\supp\underline{d}$. This is impossible since $\underline{c}_k$ is indecomposable.
\end{proof}

Assume that the Seshadri stratification is normal (Definition \ref{def:stratification:normal}). Let $\mathbb G\subseteq \Gamma$ be the set of indecomposable elements in $\Gamma$. By Proposition \ref{prop:indecomposable}, $\mathbb{G}$ is a generating set of $\Gamma$. For each $\underline{a}\in\mathbb G$, we fix a regular function $x_{\underline{a}}\in R$ satisfying $\mathcal V(x_{\underline{a}}) = \underline{a}$ and denote $\mathbb{G}_R:=\{x_{\underline{a}}\mid \underline{a}\in\mathbb G\}$.

\begin{definition}\label{def:standard:monomial}
A monomial $x_{\underline{a}_1}\cdots x_{\underline{a}_n}$ with $\underline{a}_1,\ldots,\underline{a}_n\in\mathbb{G}$ is called \emph{standard} if for each $j$ we have $\min\supp\underline{a}_j \geq \max\supp\underline{a}_{j+1}$.
\end{definition}

When writing down a standard monomial $x_{\underline{a}_1}\cdots x_{\underline{a}_n}$, it is understood that for each $j$, $\min\supp\underline{a}_j \geq \max\supp\underline{a}_{j+1}$ holds.

By Proposition \ref{Prop:UniqueDec}, any element $\underline{a}\in \Gamma$ has a unique decomposition $\underline{a} = \underline{a}_1 + \ldots + \underline{a}_n$ into indecomposable elements. By Proposition \ref{prop:indecomposable}, there exists a maximal chain $\mathfrak C$ such that 
$\supp\underline{a}_j\subseteq \mathfrak C$ for all $j=1,\ldots,n$. By Proposition~\ref{quasivaluationA}, this implies that the quasi-valuation 
is additive on standard monomials. Summarizing we have: 

\begin{proposition}\label{proposition_standard_monomial_basis}
\begin{enumerate}
\item[i)] The set $\mathbb G_R$ is a generating set for $R$.
\item[ii)] The set of standard monomials in $\mathbb G_R$ is a vector space basis for $R$. 
\item[iii)]
If $\underline{a} = \underline{a}_1 + \underline{a}_2 + \ldots + \underline{a}_n$ is the decomposition of $\underline{a}\in\Gamma$ into indecomposables,
then the standard monomial $x_{\underline{a}} := x_{\underline{a}_1}\cdots x_{\underline{a}_n}$ is such that $\mathcal V(x_{\underline{a}}) = \underline{a}$.
\item[iv)]
If a monomial $x_{\underline{a}_1}\cdots x_{\underline{a}_n}$ is not standard, then there exists a straightening relation expressing it as a linear combination of standard monomials
\[
x_{\underline{a}_1}\cdots x_{\underline{a}_n} = \sum_h u_h x_{\underline{a}_{h,1}}\cdots x_{\underline{a}_{h,n_h}},
\]
where, as in \eqref{straightening:One}, $u_h\not=0$ only if 
$\underline{a}_1+\ldots+\underline{a}_n\le^t \underline{a}_{h,1}+\ldots+\underline{a}_{h,n_h}$.
\item[v)] If in \emph{iv)} there exists a chain $\mathfrak{C}$ such that $\supp\underline{a}_i\subseteq\mathfrak{C}$ for all $i=1,\ldots,n$, and $\underline{a}'_1 + \cdots + \underline{a}'_m $ is the decomposition of $\underline{a}_1 + \cdots +\underline{a}_n \in \Gamma$ then the standard monomial $x_{\underline{a}'_1}\cdots x_{\underline{a}'_m}$ appears in the right side of the straightening relation in \emph{iv)} with a non-zero coefficient.
\end{enumerate}
\end{proposition}

We must note here that it is \emph{not} true in general that, given $p\in A$, the element $e_p\in\Gamma$ is indecomponsable. But there exists a positive integer $m_p$ such that $\frac{1}{m_p} e_p$ is the unique indecomponsable element with support $\{p\}$. 
Indeed, assume that $p\in\mathfrak{C}$ for a maximal chain $\mathfrak{C}$, then $e_p\in\Gamma_{\mathfrak{C}}$ and there exist indecomponsable elements $u_1e_p,\ldots,u_ne_p$, with $u_1,\ldots,u_n\in\mathbb{Q}_{>0}$, such that
\[
e_p = \sum_{i=1}^n u_i e_p
\]
is the decomposition of $e_p$. Now let $\frac{a}{m} e_p$, with $a,m\in\mathbb{Q}_{>0}$, be an indecomponsable element. Then we have
\[
\sum_{i=1}^n a u_i e_p = ae_p = m\cdot \frac{a}{m}e_p
\]
and, since the decomposition of $ae_p$ is unique, we have $au_i = \frac{a}{m}$ for each $i$. So $u_i = \frac{1}{m}$, $a=1$ and $\frac{1}{m}e_p$ is the unique indecomponsable element with support $\{p\}$. 

In the following, whenever the stratification is normal, we fix a function $h_p\in R$ such that $\mathcal{V}(h_p) = \frac{1}{m_p}e_p$ and we set $x_{e_p/m_p} = h_p\in\mathbb{G}_R$. The function $h_p$ does not vanish identically on $X_p$ since $\mathcal{V}_\mathfrak{C}(h_p) = \frac{1}{m_p}e_p$, where $\mathfrak{C}$ is any maximal chain in $A$ containing $p$. In some special situations, as will be shown in the proof of Theorem \ref{prop:SMT:for:subvarieties} below, the functions $h_p$, $p\in A$, have vanishing properties similar to those of the functions $f_p$.

\subsection{Balanced stratification}\label{subsec:balanced}

In the process of associating a quasi-valuation to a Seshadri stratification, the only choice we made is a total order $\leq^t$ on $A$ refining the given partial order (see Section~\ref{nonegativequasival} for all possible choices for such a refinement). To emphasize this dependence, we write ${\mathcal V}_{\leq^t}$ for the quasi-valuation and $\Gamma_{\leq^t}$ for the fan of monoids in this subsection.

\begin{definition}\label{definition_balanced_statification}
A Seshadri stratification of the variety $X$ is called \emph{balanced}\footnote{The notion of a balanced Seshadri stratification got generalized in \cite[Section 2.9]{CFL4}, where the length-preserving condition on the refinements of the partial order is removed.} if the following two properties hold:
\begin{enumerate}
\item the set $\Gamma_{\leq^t}$ of leaves for the quasi-valuation ${\mathcal V}_{\leq^t}$ is independent of the choice of the total order $\leq^t$: this allows us to simply write $\Gamma$ instead of $\Gamma_{\leq^t}$;
\item for each $\underline{a}\in\Gamma$ there exists a regular function $x_{\underline{a}}\in R$ such that ${\mathcal V}_{\leq^t}(x_{\underline{a}}) = \underline{a}$ for each possible total order $\leq^t$.\footnote{This second condition is a consequence of the first one, see \cite[Section 2.9]{CFL4} for a proof.}
\end{enumerate}
\end{definition}

\begin{rem}\label{Rmk:balance}
If the Seshadri stratification is normal and balanced, the second condition above can be weakened: it suffices to require this condition to hold for $\underline{a}\in\mathbb{G}$, the indecomposable elements in $\Gamma$.
\end{rem}

For a balanced stratification the order requirement in the straightening relations \eqref{straightening:One} is much stronger. It carries more the spirit of the classical Pl\"ucker relations.

\begin{definition}\label{definition_balanced_order}
Let $\underline{a},\,\underline{b}\in\mathbb{Q}^A$. We write $\underline{a}\trianglelefteq\underline{b}$ if $\underline{a}\leq^t\underline{b}$ for each total order $\leq^t$ on $A$ extending the given partial order $\leq$ and such that $p<^t q$ if $\ell(p) < \ell(q)$.
\end{definition}

Note that the partial order $\trianglelefteq$ is not necessarily a total order on $\mathbb{Q}^A$. Indeed, if $p$ and $q$ are not comparable in $A$ and $\ell(p) = \ell(q)$ then $e_p\not\trianglelefteq\,e_q$ and $e_q\not\trianglelefteq\,e_p$.

\begin{proposition}\label{proposition:straightening:relation}
If the Seshadri stratification is balanced, then in the straightening relation \eqref{straightening:One} with the choice of $x_{\underline{a}}$ as in Definition 15.8, we have: $u_{\underline{c}} \neq 0$ only if $\underline{a} + \underline{b} \trianglelefteq \underline{c}$.
\end{proposition}

\subsection{Compatibility with the strata}\label{subvarieties}

We have seen in Remark~\ref{induction} that each stratum $X_p$ in the Seshadri
stratification of $X$ is naturally endowed with a Seshadri stratification.
Let $\mathcal V$ be the quasi-valuation on $R=\mathbb K[\hat{X}]$ and let $\ge^t$ be the total
order on $A$ chosen in the construction of $\mathcal V$. We denote by $\ge^t$ the induced
total order in $A_p$, and let $\mathcal V_p$ be the associated quasi-valuation on
$R_p=\mathbb K[\hat{X}_p]$. So we get an associated fan of monoids etc.
It is natural to ask under which conditions these objects are compatible
with the corresponding ones for $X$. The best result is obtained in case the Seshadri 
stratification of $X$ is balanced and normal. Indeed, in this case one gets automatically a standard monomial theory
for each subvariety $X_p$.

We assume for the rest of this subsection: the Seshadri stratification of $X$ is balanced and normal.

Let $\mathbb G\subseteq \Gamma$ be the subset of indecomposable
elements and $\mathbb G_R:=\{x_{\underline{a}} \mid \underline{a}\in\mathbb G\}\subseteq R$ be a set 
of regular functions chosen as in Definition \ref{definition_balanced_statification} (2) (see Remark \ref{Rmk:balance}).

\begin{definition}\label{compatibility_subvariety}
Let $p\in A$. A standard monomial $x_{\underline{a}_1}\cdots x_{\underline{a}_n}$
on $X$ is called \emph{standard on $X_p$} if  $\max\supp{\underline{a}_1}\leq p$.
\end{definition}

By fixing an element $p\in A$ one gets various natural objects associated to $A_p$: $\mathbb G_p:=\{\underline{a}\in \mathbb G\mid
\supp \underline{a}\subseteq A_p\}$, $\mathbb G_{R_p}:=\{x_{\underline{a}}\vert_{X_p} \mid \underline{a}\in\mathbb G_p\}$, and
$\Gamma_p:=\{\underline{a}\in \Gamma\mid \supp \underline{a} \subseteq A_p\}$.
We consider the vector space $\mathbb Q^{A_p}$ as a subspace of $\mathbb Q^{A}$. Since every
maximal chain in $A_p$ is a chain in $A$, by abuse of notation we write $\mathcal V_p(f)\in \mathbb Q^{A}$ 
for a non-zero function $f\in \mathbb K[\hat{X}_p]$.

\begin{theorem}\label{prop:SMT:for:subvarieties}
If the Seshadri stratification of $X$ is balanced and normal, then the following holds:
\begin{enumerate}
    \item[i)]\label{balancednormal} For all $p\in A$, the induced Seshadri stratification on $X_p$ is balanced and normal.
    \item[ii)]\label{fan} The fan of monoids associated to $\mathcal V_p$ is equal to $\Gamma_p$, $\mathbb G_p$ is its generating set of indecomposables
    and $\mathbb G_{R_p}$ is a generating set for $R_p=\mathbb K[\hat{X}_p]$.
    \item[iii)]\label{valuation-compatibility} If $x_{\underline{a}}$ is a standard monomial, standard on $X_p$, then $\mathcal V_p(x_{\underline{a}}\vert_{X_p})=\mathcal V(x_{\underline{a}})={\underline{a}}$.
    \item[iv)]\label{restriction} The restrictions of the standard monomials $x_{\underline{a}}\vert_{X_p}$, standard on $X_p$, form a basis of $\mathbb K[\hat{X}_p]$.
    \item[v)]\label{vanishing} A standard monomial $x_{\underline{a}}$ on $X$ vanishes on the subvariety $X_p$ if and only if $x_{\underline{a}}$ is not standard on $X_p$.
    \item[vi)]\label{ideal} The vanishing ideal $\mathcal I(X_p)\subseteq R=\mathbb K[\hat{X}]$ is generated by the elements in 
    $\mathbb G_{R}\setminus \mathbb G_{R_p}$, and the ideal has as vector space basis the set of all standard monomials on $X$ which are not standard on $X_p$.
    \item[vii)]\label{intersection} For all pairs of elements $p,q\in A$, the scheme theoretic intersection $X_p\cap X_q$ is reduced. 
    It is the union of those subvarieties $X_r$ such that $r\le p$ and $r\le q$,
    endowed with the induced reduced structure. 
\end{enumerate}
\end{theorem}

\begin{proof}
We start with the most important property v). Let $x_{\underline{a}}:=x_{\underline{a}_1}\cdots x_{\underline{a}_n}$ be a standard monomial and $q:=\max\supp\underline{a}_1$. 

The monomial $x_{\underline{a}}$ is standard on $X_p$ if and only if $q\leq p$. In this case, by Proposition \ref{proposition_standard_monomial_basis} iii), we have $\mathcal{V}(x_{\underline{a}})=\underline{a}$. Let $\mathfrak{C}$ be a maximal chain containing $\supp\underline{a}$ and $p$ (the existence is guaranteed by $q\leq p$). By Proposition~\ref{supportcondition}, $\mathcal V(x_{\underline{a}}) = \mathcal V_{\mathfrak{C}}(x_{\underline{a}})$ and hence $x_{\underline{a}}$ does not vanish on $X_p$.

For the other implication, assuming $q\not\leq p$, we show that $x_{\underline{a}}$ vanishes on $X_p$. The proof is by descending induction on the length of $q$. There are two cases:
\begin{enumerate}
\item[(a).] $q$ and $p$ are comparable, i.e. $q>p$. Notice that the base step $q=p_{\max}$ of the induction is included in this case.
Let $M$ be a positive integer such that $M\underline{a}\in\mathbb{N}^A$.
\begin{enumerate}
    \item[(a.1).] Assume that $\supp\underline{a} = \{q\}$. Then $\underline{a} = n\cdot\frac{1}{m_q}e_q$ by the discussion after Proposition \ref{proposition_standard_monomial_basis}; so we can choose $M=m_q$. The monomial $x_{\underline{a}}^M = x_{M\underline{a}}$ is standard and it is equal to $h_q^{nM}$. Note that $f_q^n$ and $x_{M\underline{a}}$ have the same value $ne_q$ for each quasi-valuation $\mathcal{V}_{\leq^t}$; hence we have
    \[
    f_q^n = c\cdot h_q^{nM} + g
    \]
    where $c\in\mathbb{K}^*$ and $g$ is a linear combination of standard monomials $x_{\underline{b}}$ of the same degree of $f_q^n$ with $ne_q \trianglelefteq \underline{b}$, $ne_q\neq\underline{b}$. So, if we set $q'=\max\supp\underline{b}$ we have $\ell(q') > \ell(q)$. Hence $q' \not\leq p$ and, by induction, $g$ vanishes on $X_p$. Since $f_q$ vanishes on $X_p$ as well, it follows that $h_q$, and hence $x_{\underline{a}}$ vanish on $X_p$ too.
    \item[(a.2).] We assume that $\supp\underline{a}\neq\{q\}$. Since the Seshadri stratification is normal, there exist $0\neq\underline{a}'\in\Gamma$ and $n\geq 1$ such that $M\underline{a}=ne_q+\underline{a}'$ and $\max\supp\underline{a}'<q$. In this case $x_{M\underline{a}}=h_q^{nm_q}x_{\underline{a}'}$. Notice that $x_{\underline{a}}^M$ is \emph{not} a standard monomial: in its straightening relation we choose a standard monomial $x_{\underline{b}}$ with non-zero coefficient. By Proposition~\ref{proposition:straightening:relation}, $M\underline{a}\trianglelefteq\underline{b}$. We denote $q':=\max\supp\underline{b}$, then either $q=q'$ or $\ell(q')>\ell(q)$. 

In the first case, we can write $\underline{b}=me_q+\underline{b'}$ with a rational number $m\geq n$. Since $\underline{b}-ne_q\in\Gamma$, $x_{\underline{b}}=h_q^{n m_q}x_{\underline{b}-ne_q}$. Since we have already proved in (a.1) that $h_q$ vanishes on $X_p$, so does $x_{\underline{b}}$. In the second case, $x_{\underline{b}}$ vanishes on $X_p$ by induction.
\end{enumerate}

\item[(b).] $q$ and $p$ are not comparable. In this case $x_{\underline{a}}h_p$ is not standard: in its straightening relation we take a standard monomial $x_{\underline{b}}$ with non-zero coeffcient. By Proposition~\ref{proposition:straightening:relation}, $\underline{a} + e_p/m_p\trianglelefteq\underline{b}$. We denote $q':=\max\supp\underline{b}$, then $\ell(q')>\ell(q)$ and $\ell(q')>\ell(p)$. The last inequality implies $q'\not\leq p$, and the induction can be applied: $x_{\underline{b}}$ vanishes on $X_p$. This shows that $x_{\underline{a}}h_p$ vanishes on $X_p$. Since $h_p$ does not vanish on $X_p$ and $X_p$ is irreducible, $x_{\underline{a}}$ vanishes on $X_p$.
\end{enumerate}
This completes the proof of v).

We prove iii). Let $x_{\underline{a}}$ be a standard monomial, standard on $X_p$. By v) we know that the restriction $x_{\underline{a}}\vert_{X_p}$ does not identically vanish, so it makes sense to consider the quasi-valuation $\mathcal V_{p}(x_{\underline{a}}\vert_{X_p})$.
Let $\mathfrak{C}'\subseteq A_p$ be a maximal chain in $A_p$ and let $\mathfrak{C}$ be an extension of the chain
to a maximal chain in $A$. Keeping in mind the identification of $\mathbb Q^{A_p}$ as a subspace of $\mathbb Q^{A}$, the 
renormalizing coefficients of the valuation in Definition~\ref{chainvaluation} are chosen such that 
$\mathcal V_{p,\mathfrak{C'}}(x_{\underline{a}}\vert_{X_p})=\mathcal V_\mathfrak{C}(x_{\underline{a}})$. 
Conversely, since $\supp \underline{a} \subseteq A_p$, one can always find a maximal chain 
$\mathfrak{C'}$ in $A_p$ containing the support and extend this chain to a maximal chain $\mathfrak{C}$ in $A$.
Since $\mathfrak{C}$ contains $\supp\underline{a}$ we get: $\mathcal V(x_{\underline{a}})=\mathcal V_{\mathfrak{C}}(x_{\underline{a}})$,
and the latter is by the above equal to $\mathcal V_{p,\mathfrak{C'}}(x_{\underline{a}}\vert_{X_p})$. It follows:
$\mathcal V_{p}(x_{\underline{a}}\vert_{X_p}) =\mathcal V(x_{\underline{a}})$.

Parts i), ii), iv) follow from part iii). Part \textit{vi)} is an immediate consequence of iv) and v).

It remains to prove vii). Let $Y\subseteq X$ be the union of the subvarieties $X_r$ such that $r\le p$ and $r\le q$,
endowed with the induced reduced structure. We say a standard monomial $x_{\underline{a}}$ on $X$ is standard on $Y$
if it is standard on at least one of its irreducible components. It follows that the restriction $x_{\underline{a}}\vert_{Y}$
of a standard monomial vanishes identically on $Y$ if and only if it is not standard on $Y$. So the restrictions 
$x_{\underline{a}}\vert_{Y}$ of the standard monomials, standard on $Y$, span the homogeneous coordinate ring 
$\mathbb K[Y]$ as a vector space. Indeed, it is a basis: given a linear dependence relation between standard monomials $x_{\underline{a}}\vert_{Y}$, 
standard on $Y$, fix a $t\in A$ such that at least one summand with a non-zero coefficient is standard on $X_t$. All summands which are not standard
on $X_t$ vanish on $X_t$, so after restricting the linear dependence relation to $X_t$ one gets a non-trivial linear dependence relation
between standard monomials, standard on $X_t$, which is not possible by \textit{iv)}. It follows that the homogeneous coordinate ring 
of $Y$ has as vector space basis the standard monomials, standard on $Y$, and the vanishing ideal $\mathcal I(Y)\subseteq R$ has
as vector space basis the standard monomials which are not standard on $Y$. Using the decomposition of an element $\underline{a}\in\Gamma$
into indecomposables, we see that $\mathcal I(Y)$ is generated by those $x_{\underline{a}}\in\mathbb G_R$ such that 
$x_{\underline{a}}$ is not standard on $Y$. But this implies $x_{\underline{a}}$ is either not standard on $X_p$ and hence 
$x_{\underline{a}}\in \mathcal I(X_p)$, or $x_{\underline{a}}$ is not standard on $X_q$ and hence 
$x_{\underline{a}}\in \mathcal I(X_q)$. It follows by \textit{vi)}: $\mathcal I(X_q)+\mathcal I(X_p) = \mathcal I(Y)$, which finishes the proof.
\end{proof}

\subsection{Algorithmic aspects of standard monomials}

We restate some of the above results in the language of Khovanskii basis and discuss an implementation of the subduction algorithm (see for example \cite{KM}) to write a mononial in $R=\mathbb{K}[\hat{X}]$ into a linear combination of standard monomials. We keep notations as in previous subsections.

\begin{definition}
\begin{enumerate}
\item[(1)] A subset $\mathbb{B}\subseteq R$ is called a \emph{Khovanskii basis for the quasi-valuation} $\mathcal{V}_{\leq^t}$, if the image of $\mathbb{B}$ in $\mathrm{gr}_{\mathcal{V}_{\leq^t}}R$ generates the algebra $\mathrm{gr}_{\mathcal{V}_{\leq^t}}R$.
\item[(2)] A subset $\mathbb{B}\subseteq R$ is called a \emph{Khovanskii basis for the Seshadri stratification}, if it is a Khovanskii basis for all possible $\mathcal{V}_{\leq^t}$, where $\leq^t$ is a linear extensions of $\leq$ satisfying: if $\ell(p) < \ell(q)$ then $p<^t q$.
\end{enumerate}
\end{definition}

Combining Lemma \ref{coreminimaldifference} and Theorem \ref{fanAndDegeneratetheorem} gives

\begin{coro}
\begin{enumerate}
\item[i)] For any total order $\leq^t$, there exists a finite Khovanskii basis $\mathbb{B}_{\leq^t}$ for the quasi-valuation $\mathcal{V}_{\leq^t}$.
\item[ii)] If the Seshadri stratification is balanced, there exists a finite Khovanskii basis for the Seshadri stratification.
\item[iii)] If the Seshadri stratification is normal and balanced, the set $\mathbb{G}_R$ is a Khovanskii basis for the Seshadri stratification.
\end{enumerate}
\end{coro}

Assume hereafter that the Seshadri stratification is normal. For $g\in R$, denote by $\overline{g}$ its image in $\mathrm{gr}_{\mathcal{V}}R$.

\begin{algo}[Subduction algorithm] $ $\\
\textit{Input:} A non-zero homogeneous element $f\in R$.

\noindent
\textit{Output:} $f=\sum c_{\underline{a}_1,\ldots,\underline{a}_n}x_{\underline{a}_1}\cdots x_{\underline{a}_n}$ where $c_{\underline{a}_1,\ldots,\underline{a}_n}\in\mathbb{K}^*$ and $x_{\underline{a}_1}\cdots x_{\underline{a}_n}$ is a standard monomial.

\noindent
\textit{Algorithm:}
\begin{enumerate}
\item[(1).] Compute $\underline{a}:=\mathcal{V}(f)$ and choose a maximal chain $\mathfrak{C}$ such that $\underline{a}\in\Gamma_{\mathfrak{C}}$.
\item[(2).] Decompose $\underline{a}$ into a sum of indecomposable elements $\underline{a}=\underline{a}_1+\ldots+\underline{a}_n$ such that $\min\supp\underline{a}_i\geq\max\supp\underline{a}_{i+1}$.
\item[(3).] Compute $\overline{f}$ and $\overline{x}_{\underline{a}_1}\cdots\overline{x}_{\underline{a}_n}$ in $\mathrm{gr}_{\mathcal{V}}R$ to find $\lambda\in\mathbb{K}^*$ such that $\overline{f}=\lambda\overline{x}_{\underline{a}_1}\cdots\overline{x}_{\underline{a}_n}$.
\item[(4).] Print $\lambda x_{\underline{a}_1}\cdots x_{\underline{a}_n}$ and set $f_1:=f-\lambda x_{\underline{a}_1}\cdots x_{\underline{a}_n}$. If $f_1\neq 0$ then return to Step (1) with $f$ replaced by $f_1$.
\item[(5).] Done. 
\end{enumerate}
\end{algo}

\begin{proposition}
For any valid input, the algorithm terminates and prints out the output as in the description.
\end{proposition}

\begin{proof}
We first prove the termination of the algorithm. Assume the contrary. In the first step $\mathcal{V}(f)<\mathcal{V}(f_1)$ holds by construction. Iterating this argument gives an infinite long sequence
$$\mathcal{V}(f)<\mathcal{V}(f_1)<\mathcal{V}(f_2)<\ldots,$$
hence $f$, $f_1$, $f_2$, $\ldots$ are linearly independent. In the first step $f$ and $x_{\underline{a}_1}\cdots x_{\underline{a}_n}$ have the same degree, so does $f_1$. Repeating this argument shows that $f$, $f_1$, $f_2$, $\ldots$ are in the same homogeneous component of $R$, contradicting to the linear independency.

Once the termination is established, the correctness holds by Lemma \ref{NuLeaves} and Proposition \ref{proposition_standard_monomial_basis}.
\end{proof}

The subduction algorithm allows us to lift relations in $\mathrm{gr}_{\mathcal{V}} R$ to $R$.

Let $S$ denote the polynomial ring in variables $y_{\underline{a}}$ for $\underline{a}\in\mathbb{G}$. There is a surjective algebra morphism $\widetilde{\varphi}:S\to R$ sending the variable $y_{\underline{a}}$ to the regular function $x_{\underline{a}}\in R$. Let $I$ denote the kernel of $\widetilde{\varphi}$.

We apply the subduction algorithm to lift relations from $\mathrm{gr}_{\mathcal{V}}R$ to $R$. Let $\varphi:\mathbb{K}[t_{\underline{a}}\mid \underline{a}\in\mathbb{G}]\to \mathrm{gr}_{\mathcal{V}}R$ be the algebra morphism defined by $t_{\underline{a}}\mapsto \overline{x}_{\underline{a}}$ and $I_{\mathcal{V}}$ be its kernel. By Proposition \ref{proposition_standard_monomial_basis}, $\varphi$ is surjective. For $r(t_{\underline{a}})\in I_{\mathcal{V}}$ we set $g:=r(y_{\underline{a}})\in S$. The setup is summarized in the following diagram:
\[\xymatrix{
S:=\mathbb{K}[y_{\underline{a}}] \ar[d] \ar[r]^(.66){\widetilde{\varphi}} & R\ar[d]\\
\mathbb{K}[t_{\underline{a}}] \ar[r]^{\varphi} & \mathrm{gr}_{\mathcal{V}}R.
}
\]

We apply the subduction algorithm to $\widetilde{\varphi}(g)$ and denote the output by $h$. Being a linear combination of standard monomials, $h$ can be looked as an element in $S$. We set $\widetilde{r}:=g-h\in \ker\widetilde{\varphi}\subseteq S$.

From the subduction algorithm, all these relations $\widetilde{r}$ for $r\in I_{\mathcal{V}}$ are sufficient to rewrite a product of elements in the generating set $\mathbb{G}_R$ as a linear combination of standard monomials. This proves

\begin{coro}
The ideal $I$ is generated by $\{\widetilde{r}\mid r\in I_{\mathcal{V}}\}$.
\end{coro}

As a matter of fact, lifting a particular generating set of $I_{\mathcal{V}}$ gives a Gr\"obner basis of $I$. Details on this lifting, applications to the Koszul and Gorenstein properties, as well as the relations to Lakshmibai-Seshadri algebras, can be found in \cite{CFL3}.

\section{Examples}\label{backexamples}

\subsection{Seshadri stratifications of Hodge type}

We consider a special case of a Seshadri stratification where all bonds appearing in the extended Hasse diagram $\mathcal{G}_{\hat{A}}$ are $1$; such a Seshadri stratification will be called of \emph{Hodge type}. Their properties are summarized in the following proposition.

\begin{proposition}\label{Prop:Hodge}
The following statements hold for a Seshadri stratification of Hodge type.
\begin{enumerate}
\item[i)] For any maximal chain $\mathfrak{C}$, the monoid $\Gamma_{\mathfrak{C}}$ coincides with $\mathbb{N}^{\mathfrak{C}}$, hence the Seshadri stratification is normal.
\item[ii)] The degenerate algebra $\mathrm{gr}_{\mathcal{V}}R$ is isomorphic to the Stanley-Reisner algebra $\mathrm{SR}(A)$.
\item[iii)] There exists a flat degeneration of $X$ into a union of weighted projective spaces, one for each maximal chain in $\mathcal{C}$.
\item[iv)] If the poset $A$ is Cohen-Macaulay over $\mathbb{K}$, $X\subseteq\mathbb{P}(V)$ is projectively normal.
\item[v)] The degree of the embedded variety $X\subseteq\mathbb{P}(V)$ is
$$\sum_{\mathfrak{C}\in\mathcal{C}}\frac{1}{\prod_{q\in\mathfrak{C}}\mathrm{deg}f_q}.$$
\item[vi)] The Hilbert polynomial of $R$ is given by the formula in Proposition \ref{volumetheorem2}.
\item[vii)] The Seshadri stratification is balanced.
\item[viii)] If $deg(f_p)=1$ for all $p\in A$, then the subvarieties $X_p$ are defined in $X$ by linear equations.
\item[ix)] For any $p,q\in A$, $X_p\cap X_q$ is a reduced union of the subvarieties contained in both of them.
\end{enumerate}
\end{proposition}

\begin{proof}
By Remark \ref{Rmk:Bond1}, for any maximal chain $\mathfrak{C}$ in $A$, the lattice $L^{\mathfrak{C}}\cong\mathbb{Z}^{\mathfrak{C}}$. Together with Example \ref{extremalvaluation}, it implies $\Gamma_{\mathfrak{C}}=\mathbb{N}^{\mathfrak{C}}$. Therefore $\mathbb{K}[\Gamma_{\mathfrak{C}}]$ are polynomial algebras, and the associated (projective) toric varieties are weighted projective spaces. The first six statements follow from Proposition \ref{Prop:DecompsitionIrrComp}, Theorem \ref{fanAndDegeneratetheorem}, Theorem \ref{volumetheorem1}, Theorem \ref{Thm:ProjNormal} and Proposition \ref{volumetheorem2}.

It remains to show vii), then viii) and ix) follow from Theorem \ref{prop:SMT:for:subvarieties}. By i), all monoids $\Gamma_{\mathfrak{C}}=\mathbb{N}^{\mathfrak{C}}$ do not depend on the possible choices of the total order $\leq^t$. The set $\mathbb{G}$ of indecomposable elements is given by $\{e_p\mid p\in A\}$. In view of Remark \ref{Rmk:balance}, we choose $f_p$ as the regular function, its image under the quasi-valuation $\mathcal{V}_{\leq^t}$ is $e_p$, which is independent of the possible choices of $\leq^t$.
\end{proof}

According to Remark \ref{induction}, for any $p\in A$, the induced Seshadri stratification on $X_p$ is of Hodge type, so results in the above proposition hold for $X_p$.

We apply it to low degree projective varieties.

\begin{coro}\label{Cor:Degree2}
If $X\subseteq\mathbb{P}(V)$ is a projective variety of degree $2$ which is smooth in codimension one, then $X$ is projectively normal.
\end{coro}

\begin{proof}
Under the degree $2$ assumption, the generic hyperplane stratification for $X$ in Proposition \ref{Prop:Generic} with choices of $f_{0,k}$ as in Example \ref{Ex:Bertini1} is of Hodge type. The statement follows from Proposition \ref{Prop:Hodge} iv).
\end{proof}

\begin{example}\label{Example:Grassmann}
We consider the Seshadri stratification of the Grassmann variety $\mathrm{Gr}_d\mathbb{K}^n$ with the Pl\"ucker embedding in Example \ref{GmodBexample} where all the bonds are $1$ and the extremal functions are of degree $1$. In this case the poset $A$ is a distributive lattice, hence for any $p\in A$, the subposet $A_p$ is shellable and hence Cohen-Macaulay over any field $\mathbb{K}$ (\cite{Bj}). Applying Proposition \ref{Prop:Hodge}, we obtain:
\begin{enumerate} 
\item[i)] a degeneration of Schubert varieties $X(\underline{i})\subseteq\mathrm{Gr}_d\mathbb{K}^n$ for $\underline{i}\in I_{d,n}$ into a union of projective spaces using quasi-valuations, recovering the main results in \cite{FL};
\item[ii)] the projective normality of the Schubert varieties $X(\underline{i})$ for $\underline{i}\in I_{d,n}$ in the Pl\"ucker embedding;
\item[iii)] the degree of the embedded Schubert varieties $X(\underline{i})$ as the cardinality of $\mathcal{C}_{\underline{i}}$ for $\underline{i}\in I_{d,n}$;
\item[iv)] the Schubert varieties are defined by linear equations in the Grassmann variety;
\item[v)] the intersection of two Schubert varieties $X(\underline{i})\cap X(\underline{j})$ is a reduced union of Schubert varieties.
\end{enumerate}
The projective normality of the Schubert varieties in Grassmann varieties are proved by Hochster \cite{Hoch}, Laksov \cite{Lak}, and Musili \cite{Mus} (see also the work of Igusa \cite{Igu} for the Grassmann varieties themselves). Our approach is a geometrization of the one in the framework of Hodge algebra by De Concini, Eisenbud and Procesi in \cite{DEP}.
\end{example}

\begin{rem}\label{Rmk:Positroid}
On Grassmann varieties there exist many Seshadri stratifications. For example, all positroid varieties \cite{KLS} in $X$, together with some well-chosen extremal functions, form a Seshadri stratification on $\mathrm{Gr}_2\mathbb{C}^4$. Details will be given in a forthcoming work.
\end{rem}

\begin{example}
For any finite lattice $\mathcal{L}$ with meet operation $\wedge$ and join operation $\vee$, Hibi \cite{Hibi} introduced a graded $\mathbb{K}$-algebra 
$$\mathcal{R}_{\mathbb{K}}(\mathcal{L}):=\mathbb{K}[x_\ell\mid\ell\in\mathcal{L}]/(x_{\ell}x_{\ell'}-x_{\ell\wedge\ell'}x_{\ell\vee\ell'}\mid \ell,\ell'\in\mathcal{L}\text{ non comparable}),$$ 
which is an integral domain if and only if $\mathcal{L}$ is a distributive lattice. In this case, one obtains a projective toric variety $Y_{\mathcal{L}}\subseteq\mathbb{P}(\mathbb{K}^{|\mathcal{L}|})$, called Hibi toric variety. For $p\in\mathcal{L}$, $\mathcal{L}_p:=\{\ell\in\mathcal{L}\mid \ell\leq p\}$ is again a distributive lattice. We leave to the reader to verify that the collection of projective subvarieties $Y_{\mathcal{L}_p}$ and $f_p:=x_p$ for $p\in\mathcal{L}$ defines a Seshadri stratification of Hodge type on $Y_{\mathcal{L}}$. Proposition \ref{Prop:Hodge} recovers the well-known degeneration of the Hibi toric variety into the variety associated to the Stanley-Reisner algebra $\mathrm{SR}(\mathcal{L})$ and its projective normality.
\end{example}

\subsection{Compactification of a maximal torus}\label{subsection_torus_compactification}

The action of the torus $T$ of $\mathrm{PSL}_3(\mathbb{C})$ on $\mathfrak{sl}_3(\mathbb{C})$ defines an embedding of $T$ in $\mathbb{P}(\textrm{End}(\mathfrak{sl}_3(\mathbb{C})))$. Since this is a diagonal action and the weights of $\mathfrak{sl}_3(\mathbb{C})$ are the union $\Phi_0$ of $0$ and the root system $\Phi$, we get an embedding
\[
T\,\ni\,t\,\longmapsto\,[t^\gamma\,|\,\gamma\in\Phi_0]\in\mathbb{P}^6(\mathbb{C}).
\]
Let $X\subseteq\mathbb{P}^6(\mathbb{C})$ be the torus compactification given by the closure of the image of this embedding. We want to show that the $T$--orbit closures in $X$ are the strata for a Seshadri stratification.

If we denote by $[x_\gamma\,|\,\gamma\in\Phi_0]$ the homogeneous coordinates in $\mathbb{P}^6(\mathbb{C})$ then the equations defining $X$ in $\mathbb{P}^6(\mathbb{C})$ are
\[
x_\gamma x_\delta = x_\epsilon x_\eta,\quad\gamma,\delta,\epsilon,\eta\in\Phi_0\text{ such that }\gamma+\delta = \epsilon+\eta.
\]
It is well known that the $T$-orbit closures in $X$ are in bijection with the faces of the polyhedral decomposition of $\mathbb{R}^2$ given by the Weyl chambers. Further, it is easy to check that:
\begin{itemize}
    \item[(i)] we have six $0$--dimensional (closed) orbits $\{p_\gamma\}$, $\gamma\in\Phi$, where $p_\gamma\in X$ has coordinates $x_\gamma(p_\gamma) = 1$ and $x_\delta(p_\gamma) = 0$ for $\delta\neq\gamma$,
    \item[(ii)] the six $1$--dimensional orbit closures are the projective lines $\ell_{\gamma,\delta}$ passing through $p_\gamma,\,p_\delta$, with $\gamma,\delta\in\Phi$, $\gamma\neq\delta$ and $\gamma+\delta\not\in\Phi_0$.
    \item[(iii)] $X$ is the unique $2$--dimensional orbit closure.
\end{itemize}
Clearly these varieties are all smooth in codimension one.

Let us denote by $\alpha,\beta,\theta$ the three positive roots of $\Phi$ such that $\alpha + \beta = \theta$. The inclusion relations among these orbits, i.e. the same inclusion relations among the faces of the polyhedral decomposition of $\mathbb{R}^2$, are as follows:
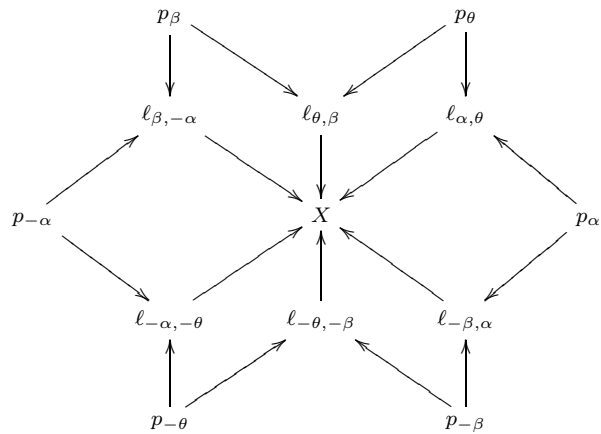
\[{\tiny
\xymatrix{
 & p_\beta\ar[d]\ar[dr] & & p_\theta\ar[d]\ar[dl] & \\
 & \ell_{\beta, -\alpha}\ar[dr] & \ell_{\theta,\beta}\ar[d] & \ell_{\alpha,\theta}\ar[dl] & \\
p_{-\alpha}\ar[ur]\ar[dr] & & X & & p_{\alpha}\ar[ul]\ar[dl] \\
 & \ell_{-\alpha, -\theta}\ar[ur] & \ell_{-\theta,-\beta}\ar[u] & \ell_{-\beta,\alpha}\ar[ul] & \\
 & p_{-\theta}\ar[u]\ar[ur] & & p_{-\beta}\ar[u]\ar[ul] & \\
}}
\]
If we choose as extremal functions
\begin{itemize}
    \item[(i)] $x_\gamma$ for the orbit $p_\gamma$,
    \item[(ii)] $x_\gamma x_\delta$ for the orbit closure $\ell_{\gamma,\delta}$,
    \item[(iii)] $x_0$ for $X$
\end{itemize}
then we get a Seshadri stratification with all bonds equal to $1$.

\begin{rem}
A normal projective toric variety admits a Seshadri stratification where the subvarieties are orbit closures arising from the torus action. We will return to this example in a separate work.
\end{rem}

\subsection{A compactification of $\mathrm{PSL}_2(\mathbb C)$}\label{subsection_compactification_group}

The De Concini--Procesi compactification $X$ of $G = \mathrm{SL}_2(\mathbb{C})$ is the projective space of all $2\times2$ matrices (see \cite{DCP2}). The group $G\times G$ acts on $X$ by left and right multiplication
\[
(g,h)\cdot[A] = [gAh^{-1}].
\]
If $B$ is the Borel subgroup of upper triangular matrices of $G$, then the $B\times B$--orbit closures in $X$ are: $X_0 = \left\{\left[\begin{array}{cc}0 & 1 \\ 0 & 0\\ \end{array}\right]\right\}$, a unique point, $X_1 = \left\{\left[\begin{array}{cc}0 & * \\ 0 & *\\ \end{array}\right]\right\}$ and $X_2 = \left\{\left[\begin{array}{cc}* & * \\ 0 & 0\\ \end{array}\right]\right\}$, two projective lines,
 $X_3 = \left\{\left[\begin{array}{cc}* & * \\ * & *\\ \end{array}\right]\textrm{ of rank one }\right\}$, a smooth quadric, $X_4 = \left\{\left[\begin{array}{cc}* & * \\ 0 & *\\ \end{array}\right]\right\}$, a projective plane, and $X_5 = X \simeq \mathbb{P}^3(\mathbb{C})$.

Let $\left(\begin{array}{cc}x & y\\ z & t\end{array}\right)$ be the coordinates on the space of matrices and define the following functions: $f_0 = y$, $f_1 = t$, $f_2 = x$, $f_3 = z$, $f_4 = xt - yz$ and $f_5 = z(xt - yz)$.
 
The inclusion relations for the subvarieties and the associated functions are
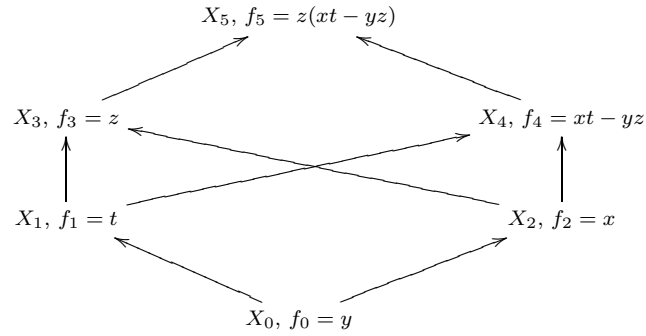
\[{\tiny
\xymatrix{
	& X_5,\,f_5 = z(xt - yz)\ar@{<-}[dl]\ar@{<-}[dr]\\
X_3,\,f_3 = z\ar@{<-}[d]\ar@{<-}[drr] & & X_4,\, f_4 = xt-yz\ar@{<-}[d]\ar@{<-}[dll]\\
X_1,\,f_1 = t\ar@{<-}[dr] & & X_2,\,f_2 = x\ar@{<-}[dl]\\
 & X_0,\,f_0 = y\\
}}
\]

It is easy to check that the orbit closures $X_0, \ldots, X_5$ with the corresponding extremal functions $f_0, \ldots, f_5$ give a Seshadri stratification of $X$.

\subsection{A family of quadrics}\label{subsection_quadric_example}

We see a Seshadri stratification for certain quadrics. Let $[x:y:z:t]$ be the homogeneous coordinates on $\mathbb{P}^3(\mathbb{C})$, let $h = ax + by + cz$ be a linear polynomial where $(a,b)\neq(0,0)$ and $(a,c)\neq(0,0)$ and consider the quadric $X$ defined by $yz - th = 0$ in $\mathbb{P}^3(\mathbb{C})$.

The quadric is smooth if $a\neq 0$ and has $[1:0:0:0]$ as unique singular point if $a=0$. So it is smooth in codimension one.

We define the strata as follows: $X_3 := X$, $X_2 := \{[x:y:z:t]\in X\mid z = t = 0\}$, $X_1 := \{[x:y:z:t]\in X\mid y = t = 0\}$ and $X_0 := \{[1:0:0:0]\}$. It is clear that $X_2$ is the projective line in $\mathbb{P}^3$ defined by $z = t = 0$ and $X_1$ is the projective line in $\mathbb{P}^3$ defined by $y = t = 0$.
We take as extremal functions: $f_3 := t$, $f_2 := y$, $f_1 := z$ and $f_0 := x$. One easily verifies that these data give a Seshadri stratification for $X$.
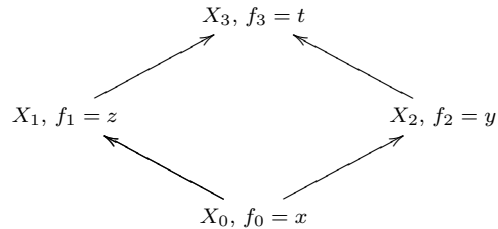
\[{\tiny
\xymatrix{
	& X_3,\,f_3 = t\ar@{<-}[dl]\ar@{<-}[dr]\\
X_1,\,f_1 = z\ar@{<-}[dr]\ar@{<-}[dr] & & X_2,\,f_2 =y\ar@{<-}[dl]\\
& X_0,\,f_0 = x\\
}}
\]

\begin{example}\label{CounterExamplePositive}
We consider the special case $a=0$, $b=c=1$: the quadric $X$ is defined by the homogeneous equation $yz-t(y+z)$. The poset $A=\{p_3,p_2,p_1,p_0\}$ with $X_{p_i}:=X_i$ for $i=0,1,2,3$. Let $\mathfrak{C}=(p_3,p_1,p_0)$ be the maximal chain to the left of the above diagram. 

We claim that the valuation $\mathcal{V}_\mathfrak{C}(y)=(1,0,0)$. Indeed, since all bonds are $1$, we set $N=1$. The function $y$ vanishes on $X_1$ with order $1$; and the function $\frac{y}{t}=1+\frac{y}{z}$, once restricted to $X_1$, is the constant function $1$. 

The function $y$ is positive along the chain $\mathfrak{C}$, but by Lemma \ref{extremalvaluationII}, it is not standard along the chain $\mathfrak{C}$. In fact, we have the following relation in $\mathbb{K}[U_\mathfrak{C}]$:
$$y=t+\frac{ty}{z},\ \text{ with } \mathcal{V}_\mathfrak{C}\left(\frac{ty}{z}\right)=(2,-1,0)>(1,0,0),$$
which explains the equality $\mathcal{V}_\mathfrak{C}(y)=\mathcal{V}_\mathfrak{C}(t)$.
\end{example}

\subsection{An elliptic curve}
Let $X$ be the elliptic curve defined by the homogeneous equation $y^2z-x^3+xz^2$ in $\mathbb{P}^2$. Let $R:=\mathbb{K}[x,y,z]/(y^2z-x^3+xz^2)$ be its homogeneous coordinate ring. We abuse the notation and use the same variables $x,y,z$ for their classes in $R$.

We present two Seshadri stratifications on $X$.
\begin{enumerate}
\item We consider the Seshadri stratification on $X$ defined as follows. The subvarieties are $X_1:=X$ and $X_0:=\{[0:1:0]\}$. As extremal functions we choose $f_1=z$ and $f_0=y$. The Hasse graph is a chain  $X_1\leftarrow X_0$. The vanishing order of $z$ at $X_0$ is $3$, which is the bond between $X_1$ and $X_0$. The quasi-valuation $\mathcal{V}$ is the valuation associated to this unique chain; it is determined by $\mathcal{V}(x)=(1/3,2/3)$, $\mathcal{V}(y)=(0,1)$ and $\mathcal{V}(z)=(1,0)$. Having different valuations, the monomials $x^py^qz^r$, where $p=0,1,2$, $q,r\in\mathbb{N}$ are linearly independent, hence form a basis of $R$. It follows that the valuation monoid is generated by $\mathcal{V}(x)$, $\mathcal{V}(y)$ and $\mathcal{V}(z)$. Theorem \ref{Thm:Degeneration} produces a toric degeneration of $X$ to the toric variety $\mathrm{Proj}(\mathbb{K}[x,y,z]/(y^2z-x^3))$. The Newton-Okounkov complex is in this case the segment connecting $(1,0)$ and $(0,1)$.

This toric degeneration is obtained in \cite{KMM} using a filtration arising from symbolic powers.
\item We give another Seshadri stratification on $X$. The subvarieties are $X_1:=X$, $X_{0,1}:=\{[0:1:0]\}$
and $X_{0,2}:=\{[0:0:1]\}$. The extremal functions are $f_1=x$, $f_{0,1}:=y$ and $f_{0,2}:=z$. The Hasse graph with bonds is depicted below:
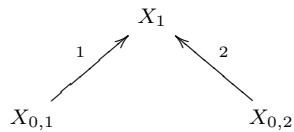
\[{\tiny
\xymatrix{
& X_1 &\\
X_{0,1}\ar[ur]^1 & & X_{0,2}\ar[ul]_2
}}
\]
We leave to the reader to verify that this is indeed a Seshadri stratification, and the degeneration predicted by Theorem \ref{Thm:Degeneration} is the union of two toric varieties, one of the two is a $\mathbb{P}^1$, and the normalization of the other one is a twisted cubic. The Newton-Okounkov complex is a union of two segments at a point. The degree formula reads: $3=1+2$.
\end{enumerate}

\subsection{Flag varieties and Schubert varieties}\label{Schubert}
An important motivation for introducing Seshadri stratifications was the aim to understand
standard monomial theory on Schubert varieties and the associated combinatorics (as developed by Lakshmibai, Musili,
Seshadri and others \cite{LR,LSII,LSIII,LSIV,LS,LS2,L1,Ses,Ses2,SMT,S2}) in the framework of Newton-Okounkov theory.

We present in this subsection just an announcement, the detailed proofs are published in a separate article \cite{CFL}. A different approach without using quantum groups at roots of unity is carried out in \cite{CFL4}. We stick in the
following for simplicity to the case of flag varieties $G/B$ for $G$ a simple simply connected algebraic group,
though the results hold, with an appropriate reformulation, also for Schubert varieties (more generally, unions thereof) contained in a 
partial flag variety in the symmetrizable Kac-Moody case (see \emph{loc.cit}). 

We fix a Borel subgroup $B$ of $G$, a maximal torus $T\subseteq B$ and a regular dominant weight $\lambda$.
The associated line bundle $\mathcal L_\lambda$ on $X:=G/B$ is ample, and we have a corresponding
embedding $G/B\hookrightarrow \mathbb P(V(\lambda))$, where $V(\lambda)$ is the Weyl module 
of highest weight $\lambda$. Let $W$ be the Weyl group of $G$, endowed with the partial order
given by its structure as Coxeter group.

\subsubsection{Seshadri stratification of $G/B$}\label{Seshadri_stratification} 
The Bruhat decomposition $G=\bigcup_{w\in W} BwB$ of $G$ implies: 
$G/B$ has a decomposition into cells $C(w):=BwB/B$, called Schubert cells. The closure
of a cell is called a Schubert variety $X(w):=\overline{C(w)}$. These varieties have an induced 
cellular decomposition: $X(w)=\bigcup_{u\le w}{C(u)}$. Schubert varieties
are known to be smooth in codimension one, see for example \cite[Corollary 4.4.5]{SMT}.

Let $\Lambda:=\Lambda(T)$ be the character group of $T$ and let $\Lambda_\mathbb R$ be the Euclidean 
vector space $\Lambda\otimes_{\mathbb Z}\mathbb R$, endowed with the Killing form as a $W$-invariant scalar product.
Let $P_\lambda\subseteq \Lambda_\mathbb R$ be the polytope obtained as the 
convex hull of all characters with a non-zero weight space in $V(\lambda)$. The vertices of this polytope
are the characters $\{\sigma(\lambda)\mid \sigma\in W\}$.

Fix a highest weight vector $v_\lambda\in V(\lambda)$ and for $\sigma\in W$ set $v_\sigma=n_\sigma(v_\lambda)$,
where $n_\sigma\in N_G(T)$ is a representative of $\sigma$ in the normalizer in $G$ of $T$. 
Then $v_\sigma$ is a $T$-eigenvector for the character $\sigma(\lambda)$. 
The cell $C(\sigma)\subseteq G/B\subseteq \mathbb P(V(\lambda))$ can be identified with the $B$-orbit $B.[v_\sigma]$.

Fix a lowest weight vector $\ell_{-\lambda}\in V(\lambda)^*$. For $\sigma\in W$
let $f_\sigma=n_\sigma(\ell_{-\lambda})\in V(\lambda)^*$, then $f_\sigma(v_\tau)\not=0$ if and only if $\tau=\sigma$.
Using the cellular decomposition, one sees that $\mathcal{H}_{f_\sigma}\cap X(\sigma)$
is the union of all codimension one Schubert varieties in $X(\sigma)$. And the property 
of the weights $\sigma(\lambda)$ to be vertices in $P_\lambda$ can be used to show
that $f_\sigma\vert_{C(\tau)}\not\equiv 0$ implies $[v_\sigma]\in X(\tau)$. 
It follows $X(\sigma)=\overline{B.[v_\sigma]}\subseteq X(\tau)$, and hence $\tau\ge \sigma$.
So for any $\sigma\in W$, $f_\sigma$ vanishes on $X(\tau)$ if $\sigma\not\le \tau$.

Therefore the collection of subvarieties $X(\tau)$ and functions $f_\tau$, $\tau\in W$,
satisfies the three axioms for a Seshadri stratification of  $X=G/B$.

It is clearly evident that to establish a standard monomial theory in such a general setting, one needs good candidates: for the monoid
$\Gamma_{\mathfrak C}$ we present a candidate $L^+_{\mathfrak C,\lambda}$, which is a reincarnation of the path model given by the 
Lakshmibai-Seshadri paths (\emph{LS-paths} for short) of shape $\lambda$. And for each element $\pi\in L^+_{\mathfrak C,\lambda}$ of degree one, 
we present a linear function $p_\pi\in V(\lambda)^*$, a candidate for a representative of a leaf of the quasi-valuation $\mathcal V$. 
The points to prove are: $L^+_{\mathfrak C,\lambda}=\Gamma_{\mathfrak C}$,
$\mathcal V(p_\pi)=\pi$, and the standard monomials (as in Section~\ref{section_standard_monomial_theory}) form a 
basis of the homogeneous coordinate ring $\mathbb K[\hat{X}]$.

\subsubsection{LS-paths} 

We recall first the notion of an LS-path $\pi$ of shape $\lambda$ (see \cite{LS2,L2}).
Given a maximal chain $\sigma=\sigma_t>\ldots >\sigma_0=\tau$ joining two elements $\sigma,\tau\in W$, $\sigma>\tau$, there exist
positive roots $\beta_1,\ldots,\beta_r$  such that $s_{\beta_i}\sigma_{i-1}=\sigma_i$ and $\ell(\sigma_{i})=\ell(\sigma_{i-1})+1$, $i=1,\ldots,t$. 
Given a rational number $a$, the chain is called an {\it $(a,\lambda)$-chain} joining $\sigma$ and $\tau$ if in addition
$a \langle \sigma_i(\lambda),\beta^\vee_i \rangle\in\mathbb Z$ for all $i=1,\ldots,t$. It has been shown in \cite{Dehy}
that  if one maximal chain between $\sigma$ and $\tau$ has the property of being an $(a,\lambda)$-chain,
then all maximal chains joining $\sigma$ and $\tau$ are $(a,\lambda)$-chains.

\begin{definition}\label{LSpath}
An \emph{LS-path} $\pi=(\underline{\sigma}:\sigma_p>\sigma_{p-1}>\ldots>\sigma_1;\underline{a}:0<a_p<\ldots<a_1=m)$ of shape $\lambda$ and degree $m\geq 1$ is a pair of sequences of linearly ordered elements in $W$ and rational numbers such that for all $i=2,\ldots,p$ there exists an $(a_i,\lambda)$-chain joining 
$\sigma_{i-1}$ and $\sigma_{i}$. Let $\mathrm{LS}(\lambda)$ denote the set of LS path of shape $\lambda$ and arbitrary degree $m\geq 0$.

An LS-path $\pi=(\sigma_p>\sigma_{p-1}>\ldots>\sigma_1;\underline{a})\in \mathrm{LS}(\lambda)$ is said to be \textit{supported in $\mathfrak{C}$}, where $\mathfrak C$ is a maximal chain in $W$, if  $\{\sigma_1,\ldots,\sigma_p\}\subseteq\mathfrak C$. 

The LS-path $\pi$ is said to be of shape $\lambda$ \textit{on a Schubert variety $X(\tau)$}
if in addition $\tau\ge \sigma_p$. We call $\sigma_p$  the \textit{initial direction} of the path, and 
$\sigma_1$ the \textit{final direction} of $\pi$. 
\end{definition}

For convenience we add the empty path $()$ as an LS-path of shape $\lambda$ and degree $0$.
The weight of an LS-path  is defined as 
$\pi(1):=\sum_{j=1}^p (a_j-a_{j+1}){\sigma_j(\lambda)}$, where $a_{p+1}:=0$. The weight of the empty path is set to be $0$.
The weight $\pi(1)$ is a weight appearing in the
Weyl module $V(m\lambda)$, and the character of this representation is given
by the following sum, running over all LS-paths of shape $\lambda$ and degree $m$ \cite{L2,L3}:
\begin{equation}\label{characterformula}
\text{Char\,}V(m\lambda)=\sum_{\pi} e^{\pi(1)}.
\end{equation}

\begin{rem}\label{upperbound}
The character formula will play a role later because it gives a dimension bound
for the graded part $\mathbb K[G/B]_m$ of the homogeneous coordinate ring.
We have a natural map from the homogeneous coordinate ring $\mathbb K[G/B]$ to the ring 
of sections $\bigoplus_{m\ge 0} \mathrm{H}^0(G/B,\mathcal L_{\lambda}^{\otimes m})$. The graded 
part of degree $m$: $\mathrm{H}^0(G/B,\mathcal L_{\lambda}^{\otimes m})$ is, as representation,
the dual $V(m\lambda)^*$ of the Weyl module $V(m\lambda)$. So the dimension of $V(m\lambda)$
is an upper bound for the dimension of $\mathbb K[G/B]_m$.
\end{rem}

\subsubsection{LS-paths as fan of monoids} 

We translate the definition of LS-paths into the language of lattices and Hasse diagrams with bonds: 
the indexing set for the Schubert varieties in $X$ is the Weyl group $W$,
endowed with the Bruhat order as partial order. The bonds are given by the Pieri-Chevalley formula:
if $\sigma$ covers $\tau$ and $\tau=s_{\beta}\sigma$ for a positive root $\beta$, then the bond is
$b_{\sigma,\tau}=\langle \tau(\lambda),\beta^\vee\rangle$.
Fix a maximal chain $\mathfrak C:\ w_0=\sigma_r>\ldots>\sigma_0=\mathrm{id}$, where $w_0$ is the maximal element in $W$.  
As before, for such a fixed maximal chain, we simplify the notation by writing $b_i$ instead of $b_{\sigma_i,\sigma_{i-1}}$ for the bonds. Note that 
all extremal functions have degree $1$ and hence $b_0=1$. Let
$L_{\mathfrak C,\lambda}\subseteq \mathbb Q^{\mathfrak C}$ be the lattice
\begin{equation}\label{EqLSLattice}
L_{\mathfrak C,\lambda}=\left\{u=\left(\begin{array}{c}u_r \\ \vdots \\ u_0\end{array}\right)\in \mathbb Q^{\mathfrak C}\left\vert 
{\scriptsize
\begin{array}{r}
b_{r}u_r\in\mathbb Z\\
b_{r-1}(u_r+u_{r-1})\in\mathbb Z\\
\ldots\\
b_{1}(u_r+u_{r-1}+\ldots+u_1)\in\mathbb Z\\ 
u_0 +u_1+\ldots+u_r \in \mathbb Z\\
\end{array} }
 \right.\right\}.
\end{equation}
We call the sum $u_0 +u_1+\ldots+u_r$ the degree of $u$. 

As before, we view $\mathbb Q^{\mathfrak C}$ as a subspace of $\mathbb Q^W$, and we consider
the union of the lattices $L_{\lambda}:=\bigcup_{\mathfrak C}L_{\mathfrak C,\lambda}\subseteq \bigcup_{\mathfrak C}\mathbb Q^{\mathfrak C}$
as a subset of $\mathbb Q^{W}$. The set $L_{\lambda}^+=L_{\lambda}\cap\mathbb Q_{\ge 0}^W$ is our candidate for the fan of monoids $\Gamma$. 
Indeed, the results in \cite{Dehy} mentioned above show that the set $L_{\lambda}^+$ is a fan of monoids.

It remains to give an explicit bijection between the set of LS-paths $\mathrm{LS}(\lambda)$ and the the fan of monoids $L_{\lambda}^+$. The proof of the following lemma will be given in \cite{CFL}, see \cite[Proposition 1]{Chi} for a proof in a slightly different language. In the following formula we set $a_{p+1}=0$.

\begin{lemma}\label{LSlattice}
The map $\nu$ defined below:
$$
\begin{array}{rlcl}
\nu:&\mathrm{LS}(\lambda)&\longrightarrow &L_{\lambda}^+=L_{\lambda}\cap \mathbb Q_{\ge 0}^{W};\\
&\pi=(\sigma_p,\ldots,\sigma_1; 0,a_p,\ldots,a_1=m)&\mapsto&
\nu(\pi):=\sum_{j=1}^p (a_j-a_{j+1})e_{\sigma_j}.\\
\end{array}
$$
induces a bijection between the set of LS-paths of shape $\lambda$ and degree $m\ge 0$, 
and the elements in $L_{\lambda}^+$ of degree $m$.
\end{lemma}
It is understood that the empty path is mapped to $0$. 

Consider again  the lattice $L_{\mathfrak C,\lambda}$ and the intersection $L^+_{\mathfrak C,\lambda}:=L_{\mathfrak C,\lambda}\cap \mathbb Q^{\mathfrak{C}}_{\ge 0}$.
By the above Lemma, we can identify this monoid with the set of all LS-paths of shape $\lambda$ and supported in $\mathfrak C$. This monoid is normal,
and hence by Proposition~\ref{Prop:UniqueDec}, every path in $L^+_{\mathfrak C,\lambda}$ can be decomposed in a unique way into a sum of indecomposable
paths. One can show (see \cite{CFL} or \cite[Proposition 3]{Chi}) that the only indecomposable paths are those of degree $1$.

Following Section~\ref{SMT:normal:leaf}, one can view a path $\pi$ of degree $m>1$ as 
a tuple $\pi=(\pi_1,\ldots,\pi_m)$ of $m$  LS-paths of shape $\lambda$ and degree $1$, satisfying  the additional condition: 
for all $i=1,\ldots,m-1$, the final  direction of $\pi_i$ is larger or equal to the initial direction of $\pi_{i+1}$.

\begin{example}
Let $G=\mathrm{SL}_3(\mathbb K)$ and $\lambda=\omega_1+\omega_2$. 

The pair $\pi=(s_1s_2s_1,s_2s_1,s_1; 0,1,\frac{3}{2},2,3)$ is an LS-paths of shape $\lambda$ and degree $3$.
Then $\pi=(\pi_1,\pi_2,\pi_3)$, where $\pi_1=(s_1s_2s_1;0,1)$, $\pi_2=(s_2s_1,s_1;0,\frac{1}{2},1)$
and $\pi_3=(s_1; 0,1)$ is the decomposition of $\pi$ into LS-paths of shape $\lambda$ and degree $1$.
\end{example}

\subsubsection{The candidates for the leaves}
It was shown in \cite{L1} that one can associate to every  LS-path of shape $\lambda$ and degree $1$
a linear function $p_\pi\in V(\lambda)^*$, called a path vector. Roughly speaking, one fixes for
$\pi=(\sigma_p,\ldots,\sigma_1; 0,a_p,\ldots,a_1=1)$ a natural number $n$ such that $na_j\in \mathbb N$ for all
$j=1,\ldots,p$. The function $p_\pi$ is then intuitively defined as $\sqrt[n]{f_{\sigma_p}^{na_p}\cdots f_{\sigma_1}^{na_1}}$,
i.e. an $n$-th root of this product of extremal functions. Indeed, we need a representation-theoretic trick: Lusztig's quantum Frobenius map
at an appropriate root of unity makes it possible to define analogues of the $n$-th roots of these kind of functions.
 
By construction, these linear functions have the following property:
$p_\pi$ vanishes identically on a Schubert variety $X(\tau)\subseteq X$ if and only if:
\begin{equation}\label{vanishingpropertypathvector}
\textrm{
$\pi=(\sigma_p,\sigma_{p-1},\ldots,\sigma_1; 0,a_p,\ldots,a_1=1)$
is such that $\sigma_p\not\le \tau$.}
\end{equation}

Let $\mathcal V:\mathbb K[\hat{X}]\rightarrow \mathbb Q^{W}$ be the 
quasi-valuation provided by the Seshadri stratification described above. The proof of the following theorem can be found in \cite{CFL}:

\begin{theorem}\label{valuationofpathvector}
For all LS-paths $\pi$ of shape $\lambda$, degree $1$ and support in $\mathfrak C$ one has 
$$\mathcal V(p_\pi)=\mathcal V_{\mathfrak C}(p_\pi)=\nu(\pi).$$ 
The value $\mathcal V(p_\pi)$ is independent of the choice of the total order in the construction of $\mathcal V$.
\end{theorem}

One has as immediate consequences of  Theorem~\ref{valuationofpathvector}:

\begin{coro}
\begin{enumerate}
\item[i)] The fan of monoids $L_{\lambda}^+$ of LS-paths of shape $\lambda$ and degree $m\ge 1$ is contained in the fan of monoids $\Gamma$.
\item[ii)] The set $\mathbb B=\{p_{\pi} \mid \pi \text{\ LS-paths of shape $\lambda$ and degree $1$}\}$ is a basis for $V(\lambda)^*$.
\end{enumerate}
\end{coro}

\begin{proof}
The first claim follows from the fact that the monoids in the fan of monoids $L_{\lambda}^+$ are generated by the degree $1$ elements.
The elements in $\mathbb B$ are linearly independent because vectors with different quasi-valuations are linearly independent \cite{KK}, 
and the character formula in \eqref{characterformula} implies  that $\mathbb B$ is a basis for $V(\lambda)^*$.
\end{proof}

Since $L_{\lambda}^+\subseteq \Gamma$, we know by Theorem~\ref{valuationofpathvector} that a path vector $p_\pi\in V(\lambda)^*$
is a representative for the leaf associated to $\nu(\pi)$. Since these elements are all of degree 1, they are indecomposable.
So we can talk about standard monomials in the sense of Definition~\ref{def:standard:monomial}. Reformulated into the language
of LS-paths this gives: a monomial of degree $m$ of path vectors : $p_{\pi_1}\cdots p_{\pi_m}$  is called a \emph{standard monomial},
if the tuple $\pi=(\pi_1,\ldots,\pi_m)$ is an LS-path of shape $\lambda$ and degree $m$.

\begin{theorem}\label{flag:standard:monomial}
\begin{enumerate}
\item[i)] The above Seshadri stratification of $G/B$ is normal and balanced.
\item[ii)] The standard monomials in the path vectors form a basis of the homogeneous coordinate ring of the embedded
flag variety $G/B\hookrightarrow \mathbb P(V(\lambda))$.
\item[iii)] This basis is compatible with the quasi-valuation $\mathcal V$.
\item[iv)] We have $\Gamma=L_{\lambda}^+$.
\end{enumerate}
\end{theorem}

\begin{proof}
Since the paths occurring in a standard monomial $p_\pi=p_{\pi_1}\cdots p_{\pi_m}$ have support in a common maximal chain, we get
by Proposition~\ref{quasivaluationA} and Theorem~\ref{valuationofpathvector}: 
$$
\mathcal V(p_\pi)=\sum_{i=1}^m\nu(\pi_j),
$$  
and the image is independent of the choice of the total order in the definition of $\mathcal V$.

The uniqueness of the decomposition of a LS-path of degree $m$ into paths of degree $1$
implies the linear independence of the set of standard monomials of degree $m$. Indeed, vectors with different quasi-valuations are linearly independent \cite{KK}.
The character formula in \eqref{characterformula} implies  that $\dim \mathbb K[X]_m\ge \dim V(m\lambda)^*$,
which by Remark~\ref{upperbound} implies that the standard monomials of degree $m$ in the path vectors form a basis 
of $\mathbb K[X]_m$. By construction, the basis is compatible with the leaves of the quasi-valuation, which implies the claim
in the theorem.
\end{proof}

As an immediate consequence, 
Proposition~\ref{proposition:straightening:relation} provides straightening relations expressing a non-standard monomial 
as a linear combination of standard monomials. 

\subsubsection{Schubert varieties}\label{Sec:Schubert}
By Remark~\ref{induction}, each stratum $X(\tau)$ in the Seshadri
stratification of $X=G/B$ is naturally endowed with a Seshadri stratification.
Since the Seshadri stratification is normal and balanced,
by Theorem~\ref{prop:SMT:for:subvarieties}, the induced Seshadri stratification
of a Schubert variety is normal and balanced. In particular,
we get a standard monomial theory for each Schubert variety which is compatible
in the sense of Section~\ref{subvarieties} with the standard monomial theory on $G/B$. 

As consequences we recover the following known results: Schubert varieties are projectively normal (Theorem~\ref{Thm:ProjNormal} and \cite{BW}; see also  \cite{RR} for a different proof using Frobenius splitting);
the degenerate variety is a union of normal toric varieties; the scheme theoretic intersection of Schubert varieties is reduced
(Theorem~\ref{prop:SMT:for:subvarieties}). In this case, further results like 
(1) vanishing theorems for higher cohomology; (2) the identification of the Newton-Okounkov simplicial complex
to the polytopes equipped with an integral structure in \cite{Dehy}; (3) the connection between the flat degeneration
in \cite{Chi} and that in Section~\ref{flatdegen}; will be discussed in \cite{CFL}.
 
We conclude this section with a degree formula, which is an application of Theorem~\ref{volumetheorem1},
but now reformulated into the setting of Schubert varieties. Denote by $V(\lambda)_\tau\subseteq V(\lambda)$
the subspace generated by the affine cone $\hat X(\tau)$ over the Schubert variety  $X(\tau)\subseteq G/B \subseteq \mathbb P(V(\lambda))$.
This subspace is called the Demazure submodule of $V(\lambda)$ associated to $\tau$.
The following formula can also be found in \cite{Chi}, and in \cite{Kn99}
in a symplectic context:

\begin{proposition}\label{Prop:degreebonds}
The degree of the embedded Schubert variety $X(\tau)\subseteq \mathbb P(V(\lambda)_\tau)$ is equal to the sum
$\sum_{\mathfrak C}\prod_{j=1}^sb_{j,\mathfrak{C}}$ running over all maximal chains $\mathfrak C$ in $A_\tau$,
and $\prod_{j=1}^sb_{j,\mathfrak{C}}$ is the product of over all bonds along a maximal chain $\mathfrak C$.
\end{proposition}

\printnoidxglossary[type=symbols,style=mcolindex,title={List of Notations}]


\begin{thebibliography}{99}

\bibitem{AK}
V. Alexeev, A. Knutson,
\emph{Complete moduli spaces of branchvarieties},
J. reine angew. Math. 639 (2010), 39--71.

\bibitem{AH}
K. Altmann, J. Hausen,
\emph{Polyhedral divisors and algebraic torus actions},
Math. Ann. 334 (2006), no. 3, 557--607.

\bibitem{And} 
D. Anderson, \emph{Okounkov bodies and toric degenerations}, Math. Ann. 356 (2013),
no. 3, 1183--1202.

\bibitem{AKL}
D. Anderson,  A. K\"uronya, V. Lozovanu, \emph{Okounkov bodies of finitely generated divisors},
International Mathematics Research Notices, Volume 2014, Issue 9, 2014, Pages 2343–2355.

\bibitem{BZ01}
A. Berenstein, A. Zelevinsky, 
\emph{Tensor product multiplicities, canonical bases and totally positive varieties}, 
Invent. Math., 143(1):77--128, 2001.

\bibitem{Bj}
A. Bj\"orner,
\emph{Shellable and Cohen-Macaulay partially ordered sets,}
Trans. Amer. Math. Soc., 260 (1) 159--183, 1980.

\bibitem{BW}
A. Bj\"orner, M. Wachs, 
\emph{Bruhat order of Coxeter groups and shellability}.
Adv. in Math. 43 (1982), no. 1, 87--100. 

\bibitem{B}
N. Bourbaki, \emph{\'El\'ements de math\'ematique. Fasc. XXXI. Alg\`ebre commutative. Chapitre 7: Diviseurs.} Actualit\'es Scientifiques et Industrielles, No. 1314 Hermann, Paris 1965.


\bibitem{Brion}
M. Brion, \emph{Lectures on the geometry of flag varieties}, Topics in cohomological studies of algebraic varieties, 33--85, 
Trends Math., Birkh\"auser, Basel, 2005. 

\bibitem{BG} 
W. Bruns and J. Gubeladze, \emph{Polytopes, rings, and K-theory}. 
Springer Monographs in Mathematics. Springer, Dordrecht, 2009.

\bibitem{Buch}
B. Buchberger, 
\emph{Ein Algorithmus zum Auffinden der Basiselemente des Restklassenringes nach einem nulldimensionalen Polynomideal}. 
Ph.D. Thesis, University of Innsbruck, 1965,

\bibitem{Cal}
P. Caldero, 
\emph{Toric degenerations of Schubert varieties},
Transform. Groups 7 (1) (2002) 51--60.

\bibitem{Chi}
R. Chiriv\`i,
\emph{LS algebras and application to Schubert varieties}, Transform. Groups 5 (2000), no. 3, 245--264.

\bibitem{Chi2}
R. Chiriv\`i,
\emph{On some properties of LS algebras}, Communication in Contemporary Mathematics, Vol. 22, No. 02, 1850085 (2020).

\bibitem{CFL}
R. Chiriv\`i, X. Fang, P. Littelmann,
\emph{Seshadri stratification for Schubert varieties 
and standard monomial theory}, Proc Math Sci 132, 74 (2022). Special Issue in Memory of Professor C S Seshadri.

\bibitem{CFL2}
R. Chiriv\`i, X. Fang, P. Littelmann, 
\emph{LS algebras, valuations and Schubert varieties}, Atti Accad. Naz. Lincei Cl. Sci. Fis. Mat. Natur. 33 (2022), no. 4, pp. 925--957.

\bibitem{CFL3}
R. Chiriv\`i, X. Fang, P. Littelmann, 
\emph{On normal Seshadri stratifications}, arXiv:2206.13171, to appear in Pure Appl. Math. Q., special volume for Claudio Procesi.

\bibitem{CFL4}
R. Chiriv\`i, X. Fang, P. Littelmann, 
\emph{Seshadri stratifications and Schubert varieties: a geometric construction of a standard monomial theory}, arXiv:2207.08904, to appear in Pure Appl. Math. Q., special volume for Corrado De Concini.

\bibitem{CLS}
D. Cox, J. Little and H. Schenck, 
\emph{Toric Varieties}. Graduate Studies in Mathematics. vol. 124. American Mathematical Society, Providence (2011).

\bibitem{DEP}
C. De Concini, D. Eisenbud and C. Procesi, 
\emph{Hodge algebras}, Ast\'erisque, 91, Soci\'et\'e Math\'ematique de France, Paris 1982.

\bibitem{DL}
C. De Concini and V. Lakshmibai,
\emph{Arithmetic Cohen-Macaulayness and Arithmetic Normality for Schubert Varieties},
American Journal of Mathematics Vol. 103, No. 5 (Oct., 1981), pp. 835--850.

\bibitem{DCP}
C. De Concini and C. Procesi,
\emph{A characteristic free approach to invariant theory},
Advances in Mathematics, Volume 21, Issue 3, September 1976, 330--354.

\bibitem{DCP2}
C. De Concini and C. Procesi, \emph{Complete symmetric varieties}. In: Gherardelli F. (eds) Invariant Theory, 1-44. Lecture Notes in Mathematics, vol 996 (1983). Springer, Berlin, Heidelberg.


\bibitem{Dehy} 
R. Dehy, \emph{Combinatorial Results on Demazure Modules}, J. Algebra {\bf 205}, (1998), pp.~505--524.

\bibitem{Dick} 
 L. E. Dickson, \emph{Finiteness of the odd perfect and primitive abundant numbers with
n distinct prime factors}, Amer. J. Math. 35 (1913), 413–422.

\bibitem{Ehr}
E. Ehrhart, \emph{Sur un probl\`eme de g\'eom\'etrie diophantienne lin\'eaire II}, 
J. Reine Angew. Math. 227 (1967), 25--49.

\bibitem{Eis}
D. Eisenbud, \emph{Introduction to algebras with straightening laws}, 
Ring Theory and Algebra III, Proc. of the third Oklahoma Conf., Lect. Notes in Pure and Appl. Math. No. 55, Dekker, 1980, 243--268.

\bibitem{E}
D. Eisenbud, \emph{Commutative algebra. With a view toward algebraic geometry}. Graduate Texts in Mathematics, 150. Springer-Verlag, New York, 1995. 

\bibitem{FaFoL}
X.~Fang, G.~Fourier, and P.~Littelmann, \emph{On toric degenerations of flag varieties}, in Representation Theory -- Current Trends and Perspectives, edited by H. Krause \textit{et al}, Series of Congress Reports, EMS, 2017.

\bibitem{FL}
X. Fang, P. Littelmann, \emph{From standard monomial theory to semi-toric degenerations via Newton-Okounkov bodies}, Trans. Moscow Math. Soc. 2017, 275--297.

\bibitem{GS}
M. Gross, B. Siebert, \emph{Affine manifolds, log structures, and mirror symmetry}, Turkish J. Math. 27 (2003), no. 1, 33-60.

\bibitem{Hibi1}
T. Hibi, 
\emph{Every affine graded ring has a Hodge algebra structure}, 
Rend. Sem. Mat. Univers. Politecn. Torino 44 (1986), 277--286.

\bibitem{Hibi}
T. Hibi,
\emph{Distributive lattices, affine semigroup rings and algebras with straightening laws},
Commutative algebra and combinatorics, 93--109, Adv. Stud. Pure Math., 11, North Holland,
Amsterdam, 1987.

\bibitem{Hoch}
M. Hochster,
\emph{Grassmannians and their Schubert subvarieties are arithmetically Cohen-Macaulay},
J. Algebra 25 (1973), 40--57.

\bibitem{Hodge}
W.~Hodge, \emph{Some enumerative results in the theory of forms}, Proc. Camb. Phil. Soc. 39, 22--30 (1943).

\bibitem{Igu}
J. Igusa,
\emph{On the arithmetic normality of the Grassmann variety},
Proc. Nat. Acad. Sci. U.S.A. 40 (1954), 309--313. 

\bibitem{Jou}
J.-P. Jouanolou, 
\emph{Th\'eor\`emes de Bertini et applications}, 
Progress in Mathematics, vol. 42, Birkh\"auser Boston Inc., Boston, MA, 1983.

\bibitem{K1} 
K.~Kaveh, \emph{Crystal bases and Newton-Okounkov bodies}, Duke Math. J. 164 (2015), 2461--2506.

\bibitem{KK} 
K.~Kaveh, A. G.~Khovanskii, \emph{Newton-Okounkov bodies, semigroups of integral points, graded algebras
and intersection theory}, Ann. of Math. (2) 176 (2012), no. 2, 925--978.

\bibitem{KM}
K.~Kaveh, C.~Manon, \emph{Khovanskii bases, higher rank valuations, and tropical geometry}, SIAM Journal on Applied Algebra and Geometry 3 (2), 292--336.

\bibitem{KMM}
K. Kaveh, C. Manon, T. Murata, \emph{On degenerations of projective varieties to complexity-one T-varieties}, arXiv:1708.02698, preprint.


\bibitem{Kn99}
A.~Knutson, \emph{A Littelmann-type formula for Duistermaat-Heckman measures.} Invent. Math. 135 (1999), no. 1, 185--200.

\bibitem{Kn06}
A.~Knutson, \emph{Balanced normal cones and Fulton-MacPherson's intersection theory}. Pure Appl. Math. Q. 2 (2006), no. 4, Special Issue: In honor of Robert D. MacPherson. Part 2, 1103--1130.

\bibitem{KLS}
A.~Knutson, T.~Lam and D.E.~Speyer, \emph{Positroid varieties: juggling and geometry},
Compositio Math. 149 (2013), 1710--1752.

\bibitem{LLM}
V. Lakshmibai, P. Littelmann, P. Magyar, 
\emph{Standard Monomial Theory and applications}, ``Representation Theories and Algebraic Geometry'' (A. Broer, ed.), Kluwer Academic Publishers (1998), pp 319--364.

\bibitem{LR} 
V. Lakshmibai, K. N. Raghavan, 
\emph{Standard Monomial Theory},
Encyclopaedia of Mathematical Sciences, 137. Invariant Theory and Algebraic Transformation Groups, 8. 
Springer-Verlag, Berlin, 2008

\bibitem{LSII}
V. Lakshmibai, C.S. Seshadri,  
\emph{Geometry of $G/P$. II. The work of De Concini and Procesi and the basic conjectures}, 
Proc. Indian Acad. Sci. Sect. A 87 (1978), no. 2, 1--54.

\bibitem{LSIII}
V. Lakshmibai, C. Musili, C.S. Seshadri,  
\emph{Geometry of $G/P$. III. Standard monomial theory for a quasi-minuscule P}, 
Proc. Indian Acad. Sci. Sect. A Math. Sci. 88 (1979), no. 3, 93--177.

\bibitem{LSIV}
V. Lakshmibai, C. Musili, C.S. Seshadri,  
\emph{Geometry of $G/P$. IV. Standard monomial theory for classical types}, 
Proc. Indian Acad. Sci. Sect. A Math. Sci. 88 (1979), no. 4, 279--362.

\bibitem{LS}
V. Lakshmibai, C.S. Seshadri, 
\emph{Geometry of $G/P$ - V}, J. Algebra 100 (1986), 462--557.

\bibitem{LS2}
V. Lakshmibai, C. S. Seshadri, 
\emph{Standard monomial theory}, in Proceedings of the Hyderabad Conference on Algebraic Groups, Manoj Prakashan, 1991.

\bibitem{Lak}
D. Laksov,
\emph{The arithmetic Cohen-Macaulay character of Schubert schemes}, 
Acta Math. 129 (1972), no. 1-2, 1--9. 

\bibitem{LM}
R. Lazarsfeld and M. Musta\c{t}$\breve{\rm a}$,
\emph{Convex bodies associated to linear series},
Ann. Sci. \'Ec. Norm. Sup\'er. (4) 42 (2009), no. 5, 783--835. 

\bibitem{Lit98}
P. Littelmann, \emph{Cones, crystals, and patterns}, Transformation groups, 3(2):145--179, 1998.

\bibitem{L1} P. Littelmann, \emph{Contracting modules and standard monomial theory for symmetrizable Kac-Moody algebras}, Jour. Amer. Math. Soc. 11 (1998), no. 3, 551--567.

\bibitem{L2} P. Littelmann, \emph{A Littlewood-Richardson formula for symmetrizable Kac-Moody algebras}, Invent. Math. 116 (1994), 329--346.

\bibitem{L3} P. Littelmann, \emph{Paths and root operators in representation theory}, Ann. Math. 142 (1995), 499--525.

\bibitem{Mac}
F. S. Macaulay,
\emph{Some Properties of Enumeration in the Theory of Modular Systems},
Proc. London. Math. Soc. Volume s2--26, Issue 1, 1927, Pages 531--555.

\bibitem{MS}
D. Maclagan, B. Sturmfels, \emph{Introduction to tropical geometry.}
Graduate Studies in Mathematics, 161. American Mathematical Society, Providence, RI, 2015. xii+363 pp.

\bibitem{Mus}
C. Musili,
\emph{Postulation formula for Schubert varieties},
J. Indian Math. Soc. (N.S.) 36 (1972), 143--171. 

\bibitem{N}
M. Nagata, \emph{Note on a paper of Samuel concerning asymptotic properties of ideals}, Mem.
Coll. Sci. Univ. Kyoto. Ser. A. Math. 30 (1957), 165--175.

\bibitem{O1}
A. Okounkov, \emph{Brunn-Minkowski inequality for multiplicities}. Invent. Math. 125 (1996), no. 3, 405--411.

\bibitem{O}
A. Okounkov, \emph{Multiplicities and Newton polytopes}. Kirillov's seminar on representation theory, 231--244,
Amer. Math. Soc. Transl. Ser. 2, 181, Amer. Math. Soc., Providence, RI, (1998).

\bibitem{O2}
A. Okounkov, \emph{Why would multiplicities be log-concave?} The orbit method in geometry and physics (Marseille, 2000), 329--347, Progr. Math., 213, Birkh\"auser Boston, Boston, MA, 2003. 

\bibitem{RR}
S. Ramanan and A. Ramanathan,
\emph{Projective normality of flag varieties and Schubert varieties},
Invent. Math. 79 (1985), no. 2, 217--224.

\bibitem{R1}
D. Rees, \emph{Lectures on the asymptotic theory of ideals}, 
London Mathematical Society Lecture Note Series, 113, Cambridge University Press, Cambridge, 1988.

\bibitem{R2}
D. Rees, \emph{Valuations associated with ideals (II)},
J. London. Math. Soc. 31 (1956) pp. 221--228.

\bibitem{Sam}
P. Samuel, \emph{Some asymptotic properties of ideals},
Ann. Math. 56 (1952) pp. 11--21.

\bibitem{Ses}
C. S. Seshadri, \emph{Geometry of $G/P$--I, theory of standard monomials for minuscule representations},
C. P. Ramanujan---a tribute, Tata Institute of Fundamental Research Studies in Mathematics, 8, Berlin, New York: Springer-Verlag, pp. 207--239.

\bibitem{Ses2}
C. S. Seshadri, 
\emph{The work of P. Littelmann and standard monomial theory},
Recent trends in Mathematics and Physics: A tribute to Harish Chandra, Narosa Publishing House, 1995, 178--197.

\bibitem{SMT}
C. S. Seshadri, \emph{Introduction to the theory of standard monomials}, 
Second edition, Texts and readings in Mathematics, 46, Hindustan Book Agency.

\bibitem{S2}
C. S. Seshadri, \emph{Standard monomials--A historical account}, 
Collected papers of C. S. Seshadri, Volume 2, Hindustan Book Agency, 2012.

\end{thebibliography}
\end{document}